\documentclass[reqno]{amsart}

\usepackage{amssymb,a4wide,color,eucal,enumerate,mathrsfs,ifthen}
\numberwithin{equation}{section}

\newtheorem{theorem}{Theorem}[section]

\newtheorem{corollary}[theorem]{Corollary}
\newtheorem{lemma}[theorem]{Lemma}
\newtheorem{proposition}[theorem]{Proposition}

\theoremstyle{definition}
\newtheorem{definition}[theorem]{Definition}

\newtheorem{example}[theorem]{Example}

\theoremstyle{remark}
\newtheorem{remark}[theorem]{Remark}

\numberwithin{equation}{section}

\newcommand{\E}{\mathbb{E}}

\newcommand{\N}{\mathbb{N}}

\newcommand{\R}{\mathbb{R}}

\newcommand{\EE}{\mathscr{E}}

\newcommand{\cE}{{\ensuremath{\mathcal E}}}

\newcommand{\cG}{{\ensuremath{\mathcal G}}}

\newcommand{\dd}{{\mbox{\boldmath$d$}}}

\newcommand{\mmu}{{\mbox{\boldmath$\mu$}}}

\newcommand{\beq}{\begin{equation}}
\newcommand{\beql}[1]{\begin{equation}\label{#1}}
\newcommand{\eeq}{\end{equation}}

\newcommand{\ba}{\begin{array}} \newcommand{\ea}{\end{array}}

\newcommand{\eps}{\varepsilon}

\def\d{\mathrm d}
\def\dd{\;\!\mathrm{d}} 

\newcommand{\PROOF}{{\bfseries Proof:} }
\newcommand{\QED}{\mbox{}\hfill\rule{5pt}{5pt}\medskip\par}

\def\bbM{{\mathbb M}}

 \def\bbW{{\mathbb W}} 
 
  \def\calF{{\mathcal F}}
  \def\calI{{\mathcal I}}

  \def\calU{{\mathcal U}}
  
 \def\calZ{{\mathcal Z}}
  
\def\rmd{{\mathrm d}}

 \def\rmw{{\mathrm w}} 
 
  \def\rmC{{\mathrm C}}
\def\rmD{{\mathrm D}}

\renewcommand{\div}{\mathop{\rm div}\nolimits}
\newcommand{\weaksto}{{\rightharpoonup^*}}
\newcommand{\weakto}{\rightharpoonup}
\newcommand{\restr}[1]{\lower3pt \mbox{$|_{#1}$}}
\newcommand{\Restr}[1]{\lower2pt\hbox{$|_{#1}$}}

\newcommand{\graph}{\mathop{\rm graph}}

\newcommand{\argmin}{\mathop\mathrm{Argmin}}
\newcommand{\trepar}{{|\kern-1truept|\kern-1truept|}}
\newcommand{\ltt}{{(\kern-2truept(}}
\newcommand{\rtt}{{)\kern-2truept)}}
\newcommand{\la}{\langle}
\newcommand{\ra}{\rangle}

\newcommand{\pairing}[4]{ \sideset{_{#1 }}{_{ #2}}  {\mathop{\langle #3 , #4  \rangle}}}
\newcommand{\nchi}{{\raise.4ex\hbox{$\chi$}}}
\newcommand{\down}{\downarrow}
\newcommand{\up}{\uparrow}

\newcommand{\Vnorm}[1]{\|#1\|}

\newcommand{\Vanach}{V}

\newcommand{\V}{{V}}

\newcommand{\forae}{\text{for a.a.}}
\newcommand{\foraa}{\text{for a.a.}}

\newcommand{\DDist}[2]{{\vartriangle_{#1}\kern-3pt(#2)}}

\newcommand{\DDistdue}[3]{{\vartriangle_{#1}^{#2}\kern-3pt(#3)}}

\newcommand{\ene}[2]{\cE_{#1}(#2)}

\newcommand{\diff}[2]{\diffname_{#1}(#2)}
\newcommand{\diffi}[3]{\diffname^{#1}_{#2}(#3)}

\newcommand{\cg}[1]{\cG(#1)}
\newcommand{\cgn}[1]{\cG^n(#1)}
\newcommand{\diffname}{\mathrm{F}}
\newcommand{\domainenergy}{D}

\newcommand{\dom}{\mathrm{dom}}

\newcommand{\wlim}{w\kern-2pt-\kern-6pt\mathop{\rm lim}\limits}

\newcommand{\AC}{\mathrm{AC}}

\newcommand{\Utau}{U_\tau}

\newcommand{\UU}{U}
\newcommand{\piecewiseConstant}[2]{\overline{#1}_{\kern-1pt#2}}
\newcommand{\pwC}{\piecewiseConstant}
\newcommand{\underpiecewiseConstant}[2]{\underline{#1}_{#2}}
\newcommand{\upwC}{\underpiecewiseConstant}

\newcommand{\piecewiseLinear}[2]{#1_{\kern-1pt#2}}
\newcommand{\pwL}{\piecewiseLinear}
\newcommand{\pwM}[2]{\widetilde{#1}_{\kern-1pt#2}}
 \def\trait #1 #2 #3 {\vrule width #1pt height #2pt depth #3pt}

 \newcommand{\pwCv}[3]{\overline{#1}_{\kern-1pt#2}^{#3}}
 \newcommand{\upwCv}[3]{\underline{#1}_{\kern-1pt#2}^{#3}}
 \newcommand{\pwLv}[3]{#1_{\kern-1pt#2}^{#3}}
 \newcommand{\pwMv}[3]{\widetilde{#1}_{\kern-1pt#2}^{#3}}

\renewcommand{\E}{E}

\newcommand{\Pt}[3]{\mathrm{P}_{#1}(#2,#3)}
\newcommand{\Pti}[4]{\mathrm{P}^{#1}_{#2}(#3,#4)}
\newcommand{\Ptname}{\mathrm{P}}
\newcommand{\MIN}[3]{\mathrm{A}_{#1,#2}(#3)}
\newcommand{\domaindiff}{\dom(\diffname)}
\newcommand{\enei}[3]{\cE^{#1}_{#2}(#3)}
\newcommand{\margvar}{\eta}

\newcommand{\domainenergyii}[1]{\dom(\cE^{#1})}
\newcommand{\domaindiffi}[1]{\dom(\diffname^{#1})}
\newcommand{\IMIN}[3]{\mathcal{I}_{#1,#2}(#3)}
\newcommand{\energy}{\mathscr{\E}}

\newcommand{\marginal}[3]{I_{#1}(#3,#2)}
\newcommand{\marginali}[4]{I^{#1}_{#2}(#4,#3)}
\newcommand{\marginale}{I}

\newcommand{\dire}{\delta}
\newcommand{\frsub}{\partial}
\newcommand{\lmsub}{\partial_{\mathrm{lim}}}
\newcommand{\margsub}{\widehat{\partial}}
\newcommand{\frsup}{\widehat{\partial}}
\newcommand{\sublev}[1]{\domainenergy_{#1}}
\newcommand{\phimin}{\varphi}
\newcommand{\phiminv}[1]{\phimin_{#1}}
\newcommand{\phiminn}{\phimin_{n}}

\newcommand{\phiminvv}[1]{\tilde{\phimin}({#1})}

\newcommand{\DDD}[3]{\begin{array}[t]{c}#1\vspace*{-1em}\\_{#2}\vspace*{-.5em}\\_{#3}\end{array}}
\newcommand{\ddd}[3]{\DDD{\begin{array}[t]{c}\underbrace{#1}\vspace*{.6em}\end{array}}{\text{\footnotesize #2}}{\text{\footnotesize #3}}}

\definecolor{ddmagenta}{rgb}{0.7,0,1.0}
\definecolor{ddcyan}{rgb}{0,0.1,1.0}

\definecolor{ddmagenta}{rgb}{0.7,0,1.0}
\definecolor{ddcyan}{rgb}{0,0.1,1.0}
\definecolor{dred}{rgb}{.8,0,0}

\newenvironment{RCOMM}{\color{dred} \textsf{R:}\,}{\color{black}}

\begin{document}

\title{Nonsmooth analysis of doubly nonlinear evolution equations}
\date{\today}

\author{Alexander Mielke}
\address{Weierstra\ss-Institut,
  Mohrenstra\ss{}e 39, 10117 D--Berlin and Institut f\"ur
  Mathematik, Humboldt-Universit\"at zu
  Berlin, Rudower Chaussee 25, D--12489 Berlin (Adlershof), Germany.}
\email{\ttfamily mielke\,@\,wias-berlin.de}

\author{Riccarda Rossi}
\address{Dipartimento di Matematica, Universit\`a di
  Brescia, via Valotti 9, I--25133 Brescia, Italy.}
\email{\ttfamily riccarda.rossi\,@\,ing.unibs.it}

\author{Giuseppe Savar\'e}

\address{Dipartimento di Matematica ``F.\
  Casorati'', Universit\`a di Pavia.
  Via Ferrata, 1 -- 27100 Pavia, Italy.}
\email{\ttfamily giuseppe.savare\,@\,unipv.it}
\thanks{R.R. and G.S. have  been partially supported by  a MIUR-PRIN 2008 grant for the  project "Optimal mass
transportation, geometric and functional inequalities and
applications"; A.M. was partially supported by DFG via
FOR787MicroPlast under Mie459/5-2. R.R. acknowledges the kind
hospitality of the Weierstra\ss-Institut, where part of this work
was carried out.}
\begin{abstract}
In this paper we analyze a broad class of abstract doubly nonlinear evolution equations
in Banach spaces, driven by nonsmooth and nonconvex energies. We  provide
some general sufficient conditions, on the dissipation potential and the energy functional,
 for existence of solutions to the related Cauchy problem. We    prove our main existence result by
 passing to the limit in a time-discretization scheme with variational techniques. Finally, we
 discuss an application to a material model in finite-strain elasticity.


\end{abstract}

\maketitle

\section{Introduction}
In this paper we investigate  (the Cauchy problem for) the
\emph{doubly nonlinear evolution equation}
\begin{equation}
  \label{eq:1}
  \partial\Psi(u'(t))+\diff {t}{u(t)}\ni0\quad\text{in }\V^*\quad \text{for a.a.\ }t\in (0,T).
\end{equation}
Here,
\[
\text{
$\V$ is a (separable) reflexive Banach space,}
\]
and
  \begin{equation}
  \label{general-diss}
  \Psi:\V\to [0,+\infty)\quad\text{is a convex potential with}
  \quad
  \Psi(0)=0,\quad
  \lim_{\|v\|\up+\infty} \frac{\Psi(v)}{\|v\|}=+\infty,
  \end{equation}
$\partial\Psi:\V\rightrightarrows \V^*$ is its usual (convex
analysis) subdifferential, and $\diffname:[0,T]\times
\V\rightrightarrows \V^*$
is a time-dependent family of multivalued maps which are induced by
a suitable ``(sub)differential'' (with respect to the variable $u$),
of a lower semicontinuous time-dependent
\[
\text{
 energy functional } \quad
\cE:(t,u)\in [0,T] \times \V \mapsto \ene t u\in (-\infty,+\infty].
\]
The quadruple $(\V,\cE,\Psi,\diffname)$ indeed generates what will
be later on referred to as \emph{generalized gradient system}.  The
aim of this paper is to study existence, stability and approximation
results for solutions to \emph{generalized gradient systems}, for a
large family of quadruples $(\V,\cE,\Psi,\diffname)$. Beside the
generality of the convex dissipation potential $\Psi$ (our main
assumption is that it  exhibits superlinear growth at infinity), we
aim to tackle a class of    multivalued operators $\diffname$ as
broad as possible. Furthermore, we  consider  a general dependence
of the energy functional $\cE$ on time (we refer in particular to
the properties of the map $t\mapsto \ene tu$ and the related notion
$\partial_t \cE$ of derivative with respect to time).

To highlight these issues, let us consider some motivating examples,
in an increasing order of generality.
\begin{enumerate}[\bf 1.]
\item \textbf{\emph{Finite dimensional ODE's.}}
  The simplest example of   gradient system
   is provided by a finite-dimensional space $\V=\R^d$ and
  an energy functional $\cE\in C^1([0,T]\times \V)$; in this case we take
  \[
  \diff tu = \mathrm{D}\ene tu, \ \text{with $\rmD \cE_t$
  the
  standard differential of the energy } \ u \mapsto \ene tu,
  \]
      and \eqref{eq:1} reads
  \begin{equation}
    \label{eq:3}
    \partial\Psi(u'(t))+\rmD\cE_t({u(t)})\ni0 \quad \foraa\, t\in (0,T).
  \end{equation}
  In the quadratic case $\Psi(v):=\frac 12 |v|^2$, $|\cdot|$ being the usual Euclidean norm on $\R^d$,
  \eqref{eq:3} is the Gradient Flow generated by $\cE$
  \begin{equation}
    \label{eq:4}
    u'(t)+\rmD\cE_t(u(t))=0 \quad \foraa\, t\in (0,T).
  \end{equation}
\item \textbf{\emph{$\lambda$-convex functionals, $C^1$-perturbation.}} More generally,
one can consider
energies of the form
 \begin{equation}
    \label{eq:6}
    \cE_t(v):=\E(v)-\langle \ell(t),v\rangle\quad  \text{with }\ell:[0,T]\to \V^* \
     \ \text{an external loading,}
     \end{equation}
     and  $\E : \V \to (-\infty,+\infty]$ a $\lambda$-convex functional
for some $\lambda \in \R$, i.e.\ it satisfies
  \begin{equation}
    \label{eq:7}
    \E((1-\theta)u_0+\theta u_1)\le (1-\theta)\E(u_0)+\theta\E(u_1)-\frac\lambda2\theta(1-\theta) \|u_0-u_1\|_\V^2
    \quad\text{for all } u_0,u_1\in \domainenergy.
  \end{equation}
  In this case, $\diffname$ admits the  representation
  \begin{equation}
    \label{eq:2}
    \diff t u=\diff{}u-\ell(t),\quad
    \diff{}u=\partial \E(v)\quad
    \text{for all } u\in \V,
  \end{equation}
   with $\partial\E$
  the \emph{Fr\'echet} subdifferential of $\E$, defined at $u \in \domainenergy:= \mathrm{dom}(\E)$
  by
  \begin{equation}
    \label{eq:8}
    \xi\in \diff{}u=\partial\E(u)\subset \V^*\quad \Leftrightarrow\quad
   \E(v)-\E(u) -  \la\xi,v-u\ra\geq o(\|v-u\|)\quad\text{as }v\to u\text{ in }\V.
  \end{equation}
   It is well known that for all $u \in \domainenergy$ the
(possibly empty) set $\frsub \E(u) \subset \V^*$ is weakly$^*$
closed, and  it reduces to the singleton $\{\mathrm{D}\E(u) \}$
if the functional $\E$ is G\^ateaux-differentiable at $u$.
Furthermore, if $\E$ is convex, then $\frsub \E(u)$
coincides with the subdifferential of $\E$ in the sense of convex
analysis. In such a framework,
  existence and approximation results for the generalized gradient
  system $(\V,\cE,\Psi,\partial \cE)$, with $\V$  a  reflexive
  space and $\Psi$ a general dissipation potential as in
  \eqref{general-diss},
  have  been proved by \cite{Colli-Visintin90, Colli92}, while in
  \cite{Segatti06} the long-time behavior of the solutions to $(\V,\cE,\Psi ,\partial \cE)$
  has been addressed.
  Notice that, when $\E$ is a $C^1$-perturbation of a convex functional $\E_0$, i.e.
  \begin{equation}
    \label{eq:9}
    \E(u):=\E_0(u)+\E_1(u),\quad \text{$\E_0$ convex, $\E_1$ of class $C^1$,}
  \end{equation}
 then one has the natural decomposition
  \begin{equation}
    \label{eq:10}
   \partial\E(u)= \partial\E_0(u)+\rmD\E_1(u),
  \end{equation}
which has been exploited in \cite{Schimperna-Segatti-Stefanelli07}
to prove   well-posedness (for the Cauchy problem) for the
 gradient system $(\V,\cE,\Psi,\partial \cE)$,
 and
existence of the global attractor for the related dynamical system.
\end{enumerate}
\noindent It is worthwhile mentioning that, in cases 1--2,
  the pair $(E,\diffname)$ satisfies a crucial \emph{closedness}
  property: the graph of the multivalued map $u \mapsto (E(u),\diffname(u)),$
  i.e.\  the set $\{ (u,E(u),\xi)\, :  \ u \in \domainenergy, \ \xi \in \diffname(u)
  \}  \subset \V  \times \R \times
  \V^*$,
  is strongly-strongly-weakly closed, meaning that,
  if sequences $u_n\in V, \EE_n\in \R, \xi_n\in \V^*$ are given,
  then
  \begin{equation}
    \label{eq:11}
    \big(\xi_n\in \diff{}{u_n},\quad \EE_n=\E(u_n),\quad
    u_n\to u,\quad
    \EE_n\to \EE,\quad
    \xi_n\weakto \xi\big) \ \Rightarrow \ \EE=\E(u),\quad \xi\in \diff{}u.
  \end{equation}

  Let us also emphasize that, in   cases 1--2,  under standard conditions on  the external loading $\ell$ as a
  function of time,
 the energy functional  $\ene tv= E(v) - \langle \ell(t),v\rangle $ fulfills the following \emph{chain rule}:
for all $u \in \AC ([0,T];V)$ and $\xi \in L^1 (0,T;V^*)$ with
  $\xi(t)\in \diff t{u(t)}$ for almost all $t \in (0,T)$ (where $\AC ([0,T];V)$ denotes
  the space of absolutely continuous functions on $[0,T]$ with values in $\V$), and  such that
  $\int_0^T \Psi(u'(t)) \dd t <+\infty$, $\int_0^T \Psi^*(-\xi(t)) \dd t <+\infty$, and
  $\sup_{t \in [0,T]} |\ene t{u(t)}|<+\infty$,
  then
 \begin{equation}
 \label{chain-intro}
 \begin{gathered}
 \text{the map $t \mapsto \ene t{u(t)}$  is absolutely continuous and}
 \\
         \frac {\mathrm{d}}{\mathrm{d}t}\ene t{u(t)}   =\la
         \xi(t),u'(t)\ra + \partial_t \ene t{u(t)}
      \quad\text{for a.a.}\, t \in  (0,T)\,.
    \end{gathered}
 \end{equation}
\begin{enumerate}[1]
    \setcounter{enumi}{2}
\item \textbf{\emph{Marginal functionals.}}
  There are examples when the Fr\'echet subdifferential does not satisfy the closedness
  property~\eqref{eq:11},
   see also~\cite{Rossi-Savare06}.
   A typical one, which we
 analyze in Section~\ref{s:2new} in more detail, is given by the so-called \emph{marginal functions}, which are
  defined via an infimum operation. Let us still consider a finite-dimensional case $\V=\R^d$ and a functional
  \begin{equation}
    \label{eq:13}
    \cE_t(u)=\min_{\margvar\in \mathcal{C}} \marginal{t}{u}{\margvar},
  \end{equation}
  where $\mathcal{C}$ is a compact topological space and $\marginale \in C^0([0,T]\times  \mathcal{C} \times \V;\R)$
   is such that
the functional $(t,u)\mapsto\marginal{t}{u}{\margvar}$
  is of class $C^1$ for every $\margvar \in \mathcal{C}$.
  Being $\mathcal{C}$ compact, for every $(t,u)\in [0,T]\times \V$ the set
  \begin{equation}
    \label{eq:14}
    M(t,u):=\argmin \marginal{t}{u}{\cdot}=\Big\{\margvar\in \mathcal{C}: \cE_t(u)=\marginal{t}{u}{\margvar}\Big\}
  \end{equation}
  is not empty.  If in addition
  the map $(t,\margvar,u) \mapsto \rmD_u\marginal{t}{u}{\margvar}$
   is continuous on $[0,T]\times \mathcal{C} \times \V$,
  it is not
   difficult to check that,
   if $\xi$ belongs to the Fr\'echet subdifferential $\partial\cE_t(u)$, then
  \begin{equation}
    \label{eq:15}
    \xi=\rmD_u\marginal{t}{u}{\margvar}\quad \text{for \emph{all}  }\margvar\in M(t,u).
  \end{equation}
  On the other hand, simple examples show that a limit
   $ \xi=\lim_{n\to\infty}\xi_n $
  of sequences $\xi_n$ satisfying \eqref{eq:15} will only obey the relaxed property
  \begin{equation}
    \label{eq:16}
    \xi=\rmD_u\marginal{t}{u}{\margvar}\quad \text{for \emph{some} }\margvar\in M(t,u).
  \end{equation}

   In view of property \eqref{eq:16}, it appears that,
   for
   reduced functionals of the type \eqref{eq:13},
    the  appropriate subdifferential is
    \begin{equation}
    \label{margsubdiff-def}
    \margsub \ene tu:=  \{ \rmD_u \marginal t
u{\margvar}\, : \ \margvar \in M(t,u) \},
    \end{equation}
    which will be  referred to
    as the \emph{marginal subdifferential} of $\cE$.
   We  examine this notion
with some detail in Section \ref{s:2new},  with the help of
 significant examples.
The latter also highlight that, in the case of marginal energies $\cE$ like \eqref{eq:13},
smoothness of the function $t \mapsto \ene tu$ for $u \in \V$ fixed is  no longer to be expected.
That is why,  one has to recur to
  a surrogate
  for the partial derivative $\partial_t \cE$,
  tailored to the marginal case \eqref{eq:13}.
  In Examples~\ref{ex:3.1} and
\ref{ex:3.2}, we  develop some
 heuristics for such a generalization of $\partial_t \cE$, and
 motivate the fact that this object should be also \emph{conditioned} to
 the (marginal) subdifferential of the energy with respect to the variable $u$, and therefore
 depend on the additional variable $\xi \in \margsub \ene tu$.
 This leads to a \emph{generalized derivative with respect to time}
 $\Ptname= \Pt tu{\xi}$, where $\xi \in  \margsub \ene tu$. For the marginal functional
 in \eqref{eq:13}, $\Ptname$ is defined by
  \[
\Pt t{u}{\xi}:= \sup\left\{ \partial_t \marginal{t}{u}{\margvar}\, : \  \margvar \in M(t,u), \  \xi =
\rmD_u \marginal t u{\margvar}
\right\}.
\]
  \item \textbf{\emph{General nonhomogeneous dissipation potentials.}}
  Last but not least, we emphasize that,
beside tackling  the above-mentioned nonsmoothness and nonconvexity
of the energy, at the same time we  treat  \emph{general} convex
dissipation potentials.

First of all,
 we  extend the existence
results of \cite{Rossi-Mielke-Savare08},
which also addressed  doubly nonlinear evolution equations driven by nonconvex energies.
 Moving from the analysis of
 gradient systems in a metric setting, the latter
paper examines the case of  nonconvex energy
functionals, albeit
 smoothly depending on time, but with
dissipation potentials of the form $\Psi(v) =
  \psi(\|v\|)$, where $\psi: [0,+\infty) \to [0,+\infty)$  is convex, l.s.c., and
  with superlinear growth at infinity. However,
in view of  applications it is also
interesting to deal with  dissipation potentials like
  \begin{equation}
  \label{meaningful-dissipation}
\Psi(v)= c_1 |v|^{p_1} + c_2 \|v\|^{p_2},\quad \text{with $p_1 \in
[1,\infty)$, $p_2 \in (1,\infty)$,}
  \end{equation}
  with $|\cdot|$  a second norm on $\V$.
  In particular, dissipations of the type
  \eqref{meaningful-dissipation} arise in the vanishing viscosity
  approximation of rate-independent evolutions described by the
  doubly nonlinear equation
  \begin{equation}
  \label{rate-independent}
\partial \Psi_{\mathrm{hom}}(u'(t))+\diff {t}{u(t)}\ni0\quad\text{in }\V^*\quad \text{for a.a.\ }t\in
(0,T),
  \end{equation}
 featuring the $1$-positively homogeneous dissipation potential $\Psi_{\mathrm{hom}}(v)=
 |v|$. The \emph{natural}
  viscous approximation of \eqref{rate-independent} is
 indeed the gradient system
\begin{equation}
\label{viscous-approximation}
  \partial\Psi_\eps(u'(t))+\diff {t}{u(t)}\ni0\quad\text{in }\V^*\quad \text{for a.a.\ }t\in
  (0,T), \quad \text{with }  \Psi_\eps(v) = |v| + \frac\eps2
  \|v\|^2.
\end{equation}
We mention that the vanishing viscosity limit of
\eqref{viscous-approximation} as $\eps \searrow 0$ has been studied
in \cite{mielke-rossi-savare2010} in the case of a
\emph{finite-dimensional} ambient space $V$. Moving from the
existence results for \emph{viscous} doubly nonlinear equations of
the present paper, we are going to  address the vanishing viscosity
analysis of \eqref{viscous-approximation} in an infinite-dimensional
context in the forthcoming paper \cite{mielke-rossi-savare2011}.

Secondly,  we consider dissipation potentials $\Psi= \Psi_u(v)$ also
depending on the state variable $u$,   hence address the doubly
nonlinear  equation
\begin{equation}
  \label{eq:1-bis}
  \partial\Psi_{u(t)}(u'(t))+\diff {t}{u(t)}\ni0\quad\text{in }\V^*\quad \text{for a.a.\ }t\in (0,T).
\end{equation}
\noindent A significant example for  potentials
 of this type
 will be provided in Section \ref{s:4}, focusing on a model
in finite-strain elasticity. In fact, state-dependent dissipations naturally occur in various plasticity models,
see for example  \cite{Miel03EFME,MaiMie08?GERI,baba-francfort-mora}.
\end{enumerate}

\noindent The discussion  developed throughout Examples 1--4
motivates   the analysis of \emph{generalized gradient systems}
$(\V, \cE, \Psi, \diffname, \Ptname)$ which is developed  in this
paper. As a main goal,
 we will prove an
existence and approximation result for the Cauchy problem for
\eqref{eq:1-bis}, under suitable conditions on $\cE, \Psi,
\diffname, \Ptname$. To be more precise, we will call a function $u:
[0,T] \to \V$ \emph{solution} for the generalized gradient system
 $(\V, \cE, \Psi, \diffname, \Ptname)$, if
$u\in \AC ([0,T];\V)$, and
 there exists $\xi \in L^1 (0,T;\V^*)$ such that
 \begin{subequations}
 \label{sol-def}
 \begin{align}
 & \label{sol-defa}
 \partial \Psi_{u(t)}(u'(t)) +\xi(t) \ni 0   &  \quad \foraa\, t \in (0,T),
 \\
  &
  \label{sol-defb}
   \xi(t) \in \diff t{u(t)}   &  \quad \foraa\, t \in (0,T),
   \end{align}
and $(u,\xi)$ fulfill the \emph{energy identity}
  \begin{align}
  \label{sol-defc}
  \int_0^T \Psi_{u(t)}(u'(t)) + \Psi_{u(t)}^*(-\xi(t)) \dd t +\ene {T}{u(T)} = \ene {0}{u(0)}  +\int_0^T \Pt t{u(t)}{\xi(t)} \dd t.
 \end{align}
 \end{subequations}
Let us point out that the energy identity \eqref{sol-defc} is a
crucial item in our definition of solution to $(\V, \cE, \Psi,
\diffname, \Ptname)$. On the one hand, \eqref{sol-defc} is  a
consequence of \eqref{sol-defa}--\eqref{sol-defb} and of the chain
rule \eqref{chain-intro}, as it can be checked by testing
\eqref{sol-defa} by $u'(t)$ and integrating on $(0,T)$. On the other
hand, as mentioned below,  for proving existence of solutions to
\eqref{eq:1-bis}, in fact we are going to  first derive
\eqref{sol-defb} and \eqref{sol-defc}, and then combine them to
obtain \eqref{sol-defa}.

The \textbf{plan of the paper} is as follows:
 in  Section~\ref{s:2}  we
address the analysis of the doubly nonlinear evolution equation
\eqref{eq:1} in a simplified setting: the
  dissipation potential is \emph{independent} of the variable $u$, and
  the energy  $\cE = \ene tu$ is  possibly nonsmooth and nonconvex
 with respect to   $u$, but
  \emph{smoothly} depending on time.
  Thus, throughout Sec.~\ref{s:2}, the multivalued map $\diffname:[0,T] \times \V \rightrightarrows \V^*
 $ is given by the Fr\'echet subdifferential of the energy, i.e.\ $\diff tu=\partial \ene
 tu$ for all $(t,u) \in [0,T] \times \V$, while  $\Pt tu{\xi} $
  reduces to the usual partial time-derivative $\partial_t \ene tu$.
Nonetheless, the analysis of this case still highlights the most
significant difficulties arising for nonconvex energies.
 In such a context, we enucleate
the main conditions on the energy functional for proving existence for \eqref{eq:1}. First,
we require some suitable
\emph{coercivity} property, which amounts to  asking that the sublevels
 of the energy are compact. Second, we impose that the energy $\cE: [0,T] \times \V \to \R$
  and the Fr\'echet subdifferential $\partial \cE : [0,T] \times \V \rightrightarrows \V^*$ fulfill a
 (joint) \emph{closedness property}, cf.\ \eqref{eq:11}. Third, we
 require that  a suitable  form of the chain rule \eqref{chain-intro} holds.

 Then, we state an existence result for the Cauchy problem for equation \eqref{eq:1}, and
 outline the steps of its proof, viz.\
 approximation by
time-discretization, a priori estimates on the approximate solutions, compactness   arguments, and
 the final   passage to the limit
of the time-discrete scheme. In developing the proof, we highlight
the role played
 at each step by
 the  aforementioned conditions on the energy functional. Namely, the incremental minimization
 leads to a discretized version of \eqref{sol-defa}--\eqref{sol-defb} and, using the \emph{variational interpolant}
 of the discrete solutions, we obtain a discrete upper energy estimate, corresponding to the inequality
 $\leq$ in \eqref{sol-defc}. Exploiting lower semicontinuity arguments and the closedness of the graph of
 the map $(t,u) \mapsto \{ (\ene tu, \xi, \partial_t \ene tu)\, : \ \xi \in \diff tu\}$, the passage to the
 limit yields
 \eqref{sol-defb} and \eqref{sol-defc} with $\leq$ instead of $=$.
 Hence, we employ a suitable \emph{lower chain-rule estimate} to conclude that \eqref{sol-defc}
 holds with equality. From this argument, we also have \eqref{sol-defa}.

  In Section~\ref{s:2new} we
 discuss finite-dimensional examples of marginal energy functionals. In this way,  we  motivate and develop some heuristics for
 new notions of subdifferential of the energy with respect to the variables  $t$ and to
 $u$,
  tailored to the case of
marginal functionals, viz.\ the aforementioned marginal
subdifferential $\margsub \cE$ and the generalized partial
time-derivative $\Ptname$. We emphasize that, even in
finite-dimensional cases, the nonsmoothness of $\cE$ forces us to
make $\Pt t{u}{\cdot}$ dependent on $\xi \in \diff t u$.

 From
 Section~\ref{s:3} on, we  examine the generalized gradient system $(\V, \cE, \Psi, \diffname, \Ptname)$, with a
  \emph{state-dependent} dissipation potential
$\Psi=\Psi_u(v)$. In such a context we
  state our main existence and  approximation
  result for \eqref{eq:1-bis}. We also give an upper semicontinuity result for the
  set of solutions to \eqref{eq:1-bis} with respect to perturbations of the dissipation
  potential and of the
  energy functional.

 In Section~\ref{s:4} we  present an application of our existence
 theorem
 to a
 PDE system for  material models
with finite-strain elasticity. Indeed, we consider dissipative
material models (also called \emph{generalized standard materials},
cf.\ \cite{Fremond2002,mie-generalized}) with an internal variable
$z: \Omega \to \mathrm{K} \subset \R^m$, while the elastic
deformation $\varphi: \Omega \to \R^d$ is quasistatically minimized
at each time instant. Thus, we are in the realm of marginal
functionals
\[
\ene tz = \min_{\varphi \in \mathcal{F}} \marginal{t}{z}{\varphi}, \quad \text{where }
 \marginal{t}{z}{\varphi}= \cE^1 (z) +
 \int_{\Omega} W(\nabla \varphi, z)\dd x - \langle
 \ell(t),\varphi\rangle,
\]
with $\cE^1$  a convex functional with compact sublevels.
 Here, $W(\cdot,z)$ is polyconvex to guarantee that the set of minimizers $M(t,z)$
 is compact and nonempty. We use the technical  assumption
 $|\mathrm{D}_z W(F, z)|\leq \kappa_1 (W(F, z)+\kappa_2)^{1/2}$.

  All the proofs of our abstract results are developed in
 Section~\ref{s:5}, relying on some technical tools and auxiliary results collected
 in the \bigskip Appendix.

\paragraph{\textbf{Basic set-up and notation.}}
 Hereafter, we will   set our analysis in the framework of
 \[
  \text{a reflexive separable Banach space $\V$}
\]
with norm $\| \cdot\|$. We denote by $\langle \cdot, \cdot \rangle$
the duality pairing between $\V^*$ and $\V$ and by $\| \cdot\|_*$
the norm on $\V^*$.

Our basic assumption on the energy functional $\cE: [0,T] \times \V
 \to
(-\infty,+\infty]$  is that there exists $\domainenergy \subset \V$
such that
\begin{equation}
\tag{$\mathrm{E}_0$} \label{Ezero}
\begin{gathered}
 \mathrm{dom}(\cE)=
[0,T] \times \domainenergy,   \ \ \text{the map } u\mapsto \ene t{u}
\  \text{is lower semicontinuous for all $t \in [0,T]$,}
\\
\exists\, C_0>0\, \ \forall\, (t,u) \in [0,T] \times \domainenergy
\, : \ \ene tu \geq C_0.
\end{gathered}
\end{equation}
Indeed,  if the functionals $\cE_t$ are  bounded from below by some
constant independent of $t$, up to a translation it is not
restrictive
to assume such a constant to be strictly positive. 

 Hereafter, we will  use the following notation
\begin{equation}
\label{ene-plus} \cg{u}:= \sup_{t \in [0,T]} \ene tu \quad \text{for
every $u \in \domainenergy$.}
\end{equation}
 Furthermore, we will
denote by $\diffname: [0,T] \times \domainenergy \rightrightarrows
\V^* $ a time-dependent family of multivalued maps, such that for $t
\in [0,T]$ the mapping  $\diffname_t$ is a (suitable notion of)
subdifferential of the functional $u\mapsto\cE_t (u)$. We use the
notation
\[
\begin{gathered}
\domaindiff=\left\{(t,u) \in [0,T]\times \domainenergy\, : \ \ \diff
t u \neq \emptyset \right\}, \\  \mathrm{graph}(\diffname)
=\left\{(t,u,\xi) \in [0,T] \times \domainenergy\times \V^*\, : \ \
\xi \in \diff t u \right\}
\end{gathered}
\]
for the domain and the graph of the multivalued mapping $\diffname :
[0,T] \times \domainenergy \rightrightarrows \V^* $, respectively.
The basic measurability requirement on $\diffname$ is that
$\graph(\diffname)$   is a Borel set of $[0,T]\times \V \times
\V^*$.

 In the framework of
 the space
$\R^m$, we will denote by $|\cdot| $ the Euclidean norm and
   by $B_r(0)$ the
ball  centered at $0$ and of radius $r$. The symbol $\weakto$ will
indicate  weak convergence both  in $\V$ and in $\V^*$. Finally,
throughout the paper we will use the symbols $C$ and $C'$ for
 various positive constants depending only on known quantities.
\section{Analysis in a simplified setting}
 \label{s:2}
In this section,  we deal with  a single dissipation potential
$\Psi$, independent of the
 state variable,  and an  energy functional  $\cE: [0,T]\times \V \to (-\infty,+\infty]$
 as in~\eqref{Ezero}
 with a  \emph{smooth} time-dependence (see for instance \cite[\S
3]{Mielke05}, \cite{Mielke-Rossi-Savare08} for analogous assumptions
within the analysis of abstract
 doubly nonlinear and rate-independent problems). The focus
 of this section is on the
 nonsmoothness  and nonconvexity of the map $u\mapsto \ene tu$.
We leave the questions arising from nonsmooth time-dependence and
state-dependent dissipation potentials to later sections.

 In  the present framework, it is natural to work with
the Fr\'echet subdifferential of  the functionals $\cE_t: \V \to
(-\infty,+\infty]$,  defined in \eqref{eq:8}. Hence, we address the
Cauchy problem
\begin{equation}
\label{Cauchy-special-form}
 \partial\Psi(u'(t))+\frsub \ene t{u(t)}\ni 0 \quad\text{in }\V^*\quad \text{for a.a.\ }t\in (0,T);\qquad
  u(0)=u_0,
\end{equation}
which is a particular case of \eqref{eq:1}, with $\diff t{u}=\frsub
\ene tu $
 for all $(t,u) \in [0,T]\times \domainenergy$.
 In Section \ref{ss:2.1} we
 enucleate the abstract assumptions for Theorem \ref{th:0} below
that yield existence for problem~\eqref{Cauchy-special-form}.
 Next, we
  give an outline of its proof, highlighting the role played by
 the aforementioned  assumptions.
Then, in Section~\ref{ss:2.2} we discuss sufficient conditions for
the latter. We conclude  with some  PDE applications in
Section~\ref{ss:2.3}.
 \noindent
 \subsection{An existence result}
 \label{ss:2.1}
\paragraph{\textbf{Assumptions on the dissipation potential
$\Psi$.}} Throughout this section we will suppose that 
\begin{equation}
\label{eq:psi-sum-1} \tag{$2.{\Psi_1}$}
 \Psi:\V\to [0,+\infty) \ \  \text{is  l.s.c. and  convex,}
  \end{equation}
\begin{equation}
\label{eq:psi-sum-1-bis} \tag{$2.{\Psi_2}$}
  \Psi(0)=0,\quad
  \lim_{\|v\|\up+\infty} \frac{\Psi(v)}{\|v\|}=+\infty, \lim_{\|\xi\|_*\up+\infty}
\frac{\Psi^*(\xi)}{\|\xi\|_*}=+\infty, \quad \text{and}
\end{equation}
\begin{equation}
\label{eq:psi-sum-2} \tag{$2.{\Psi_3}$}
 \forall\, w_1,\, w_2 \in \partial \Psi(v)\,: \qquad
 \Psi^*(w_1) = \Psi^*(w_2),
  \end{equation}
where $\Psi^*$ denotes the Fenchel-Moreau conjugate of $\Psi$.
Hereafter, we will call any $\Psi:\V \to [0,+\infty)$ complying
with~\eqref{eq:psi-sum-1}--\eqref{eq:psi-sum-2} an \emph{admissible
dissipation potential}.

We emphasize that, in this paper we  only consider dissipation
potentials $\Psi$ with  $\mathrm{dom}(\Psi)=\V$. From this it
follows (see, e.g.,\ \cite[Chap.\,I, Cor.\,2.5]{Ekeland-Temam74}),
that $\Psi$ is continuous on $\V$, and that $\Psi^*$ has superlinear
growth at infinity. Hence, the third of \eqref{eq:psi-sum-1-bis}
could be omitted, and has been stated here just for the sake of
analogy with condition \eqref{eq:41.1} later on, for
\emph{state-dependent} dissipation potentials $\Psi=\Psi_u(v)$.

 In fact,
our analysis can  be   extended to the case in which $\mathrm{dom}(\Psi)$ is an
open subset of $\V$ (i.e., it contains one continuity point).
However, this rules out dissipation potentials enforcing irreversible evolution, like for
example in  damage models, see  e.g.\ \cite{MieRou06RIDP}.
\begin{remark}
\label{rmk:useful-properties} \upshape
\begin{enumerate}
\item
 We point out that, since  $ \Psi(0)=0$, we have
 \begin{equation}
 \label{psipos}
\Psi^{*}(\xi) \geq 0 \qquad \text{for all $\xi \in \V^*$.}
\end{equation}
Furthermore, it follows from the superlinear growth of $\Psi$ and
$\Psi^*$ that
\begin{equation}
\label{e:bounded-operator} \partial \Psi: \V \rightrightarrows \V^*
\text{ is a bounded operator, and $\partial \Psi(v) \neq \emptyset$
for all $v \in \V$.}
\end{equation}
\item
 Let us now get some further insight into
condition~\eqref{eq:psi-sum-2}:
 a lower semicontinuous and convex
potential  $\Psi:\V\to [0,+\infty)$ satisfies~\eqref{eq:psi-sum-2}
if and only if
\begin{equation}
\label{differentiability}
 \text{the mapping } \lambda \mapsto\Psi
(\lambda v) \ \text{is differentiable at $\lambda=1$.}
\end{equation}
Indeed, let $v \in \V$ be such that $\partial \Psi(v) \neq
\emptyset$. The convexity of $\Psi$ gives
\[
\liminf_{\lambda \down 0} \frac{\Psi(v +\lambda v)-\Psi(v)}{\lambda}
\geq \langle w, v \rangle \geq \limsup_{\lambda \uparrow 0}
\frac{\Psi(v +\lambda v)-\Psi(v)}{\lambda}
\]
  for all  $w \in
\partial \Psi(v).$ Hence, \eqref{differentiability} holds if
 and only if $\langle w_1,v\rangle = \langle w_2,v\rangle $ for all
$w_1,\, w_2 \in \partial \Psi(v)$, which is obviously equivalent
to~\eqref{eq:psi-sum-2}.

 Therefore,
condition~\eqref{eq:psi-sum-2} is satisfied for example when $\Psi$
is a linear combination of (positively) homogeneous, or
differentiable, convex potentials.
\end{enumerate}
\end{remark}
\paragraph{\textbf{Assumptions on the energy functional $\cE$.\,}}
\begin{description}
\item[Coercivity]
\begin{equation}
  \label{eq:17}
\tag{$2.{\mathrm{E}_1}$}
 \exists\,\tau_o>0\, \ \ \forall\, t \in [0,T]\, : \quad
   u\mapsto \ene tu+\tau_o\Psi(u/\tau_o)\quad\text{has compact sublevels.}
\end{equation}
\item[Variational sum rule]
  If for some $u_o\in \V$ and $\tau>0$  the point $\bar u$ is a minimizer of
  $u\mapsto \ene tu+\tau\Psi((u-u_o)/\tau)$, then
    $\bar u$ satisfies the Euler-Lagrange equation
   $\partial \Psi((\bar u-u_o)/\tau) +\frsub \ene t{\bar u} \ni
  0$, viz.
  \begin{equation}
    \label{2.eq:42-bis}
    \tag{$2.{\mathrm{E}_2}$}
    \exists \,  \xi\in \frsub \ene t{\bar u} \,:\quad -\xi\in \partial \Psi((\bar u-u_o)/\tau).
  \end{equation}
\item[Time-dependence]
\begin{equation}
    \label{hyp:en3}
\tag{$2.{\mathrm{E}_3}$}
\begin{aligned}
&
    \forall \, u \in \domainenergy\, : \
    (t\mapsto \ene tu)  \ \text{ is
      differentiable on $ (0,T)$, \  with derivative } \
    \partial_t\ene tu;
 \\
  & \exists \,C_1>0 \ \ \forall\, (t,u)\in [0,T]\times \domainenergy
    \, : \ \ |\partial_t \ene tu|\leq C_1
    \ene tu.
    \end{aligned}
  \end{equation}
  \item[Chain rule] For every $u\in \AC([0,T];\V)$ and $\xi\in
    L^1(0,T;\V^*)$ with
    \begin{equation}
    \label{conditions-for-chain-rule}
    \begin{gathered}
 \sup_{t \in [0,T]} \left|\ene{t}{u(t)}
\right|<+\infty, \ \
       \xi(t)\in\frsub \ene t{u(t)} \text{ for a.a.   $\, t \in (0,T)$, and}\\
      \int_0^T\Psi(u'(t))\,\d t<+\infty,\quad
       \int_0^T\Psi^*(- \xi(t))\,\d t<+\infty,
       \end{gathered}
    \end{equation}
    the map $t \mapsto \ene t{u(t)}$  is absolutely continuous and
    \begin{equation}
    \label{eq:45tris}
    \tag{$2.{\mathrm{E}_4}$}
    \begin{gathered}
         \frac {\mathrm{d}}{\mathrm{d}t}\ene t{u(t)}   =\la
         \xi(t),u'(t)\ra + \partial_t \ene t{u(t)}
      \quad\text{for a.a.}\, t \in  (0,T)\,.
    \end{gathered}
    \end{equation}
 \item[Weak closedness of $(\cE,\frsub \cE)$]
   For all $t \in [0,T]$ and for all  sequences $(u_n) \subset \V$ and $(\xi_n) \subset \V^* $
   we have the following condition:
    \begin{equation}
      \label{eq:45}
      \tag{$2.{\mathrm{E}_5}$}
       \begin{gathered}
       \text{if $u_n\to u $ in $\V,$ } \
        \xi_n \in \frsub \ene t{u_n}, \
       \text{$\xi_n\weakto \xi$  in $\V^*$,} \
        \partial_t \ene t{u_n}  \to p, \
 \text{$\ene t{u_n}\to \EE$ in $\R$,}
  \\
    \text{then} \quad
      \xi \in \frsub \ene tu,\quad p \leq \partial_t \ene tu, \quad
      \EE=\ene tu.
      \end{gathered}
    \end{equation}
\end{description}
\noindent
 A few comments on the above
abstract conditions are in order:
 \begin{enumerate}
 \item
 in Proposition~\ref{l:variational-sum-rule} we are going to show that
  the variational sum rule
\eqref{2.eq:42-bis} is indeed a consequence of the closedness
property \eqref{eq:45};
\item
 in Section~\ref{ss:2.2} we
discuss sufficient conditions for  \eqref{eq:45} and the chain rule
\eqref{eq:45tris}, showing in particular that they are valid if the
functionals $\ene t{\cdot}$ are $\lambda$-convex;
\item
\eqref{hyp:en3}  and the Gronwall Lemma yield the
  following estimate
  \begin{equation}
    \label{deri-en-3}
    \ene tu\leq  \exp(C_1|t-s|)\ene su \quad \text{for all $t,s \in [0,T],\ u \in
    \domainenergy$,}
  \end{equation}
  whence, in particular,
  \begin{equation}
  \label{eq2.6}
\cg{u} \leq \exp(C_1 T)\ene tu \quad \text{for all $t \in [0,T]$.}
  \end{equation}
  \end{enumerate}
\paragraph{\textbf{Existence theorem and outline of the proof.}}
We are now in the position to state the main result of this section.
\begin{theorem}
\label{th:0}  Let us assume that $(\V,\cE,\Psi,\frsub \cE,\partial_t
\cE)$ comply with~\eqref{eq:psi-sum-1}--\eqref{eq:psi-sum-2}
and~\eqref{Ezero}, \eqref{eq:17}--\eqref{eq:45tris}. Then, for every
$u_0 \in \domainenergy$  there exists  a curve $ u \in \AC
([0,T];\V)$ solving the Cauchy problem~\eqref{Cauchy-special-form}.
In fact, there exists a function $\xi \in L^1 (0,T;\V^*)$ fulfilling
\[
\xi(t) \in \frsub \ene t{u(t)} \cap (-\partial \Psi(u'(t))) \quad
\foraa\, t \in (0,T),
\]
and the \emph{energy identity} for all $0 \leq s \leq t \leq T$
  \begin{equation}
    \label{eq:52bis}
    \int_{s}^{t}
    \Big(\Psi(u'(r)){+}\Psi^*(-\xi(r)) \Big)\,\d r + \ene {t}{u(t)}=\ene
    {s}{u(s)} + \int_s^t \partial_t \ene{r}{u(r)} \, \d r.
  \end{equation}
\end{theorem}
\noindent Theorem~\ref{th:0} is a direct consequence of the more
general Theorem~\ref{thm:viscous2}, which is proved in
Section~\ref{s:5}. Nonetheless, in order to provide some heuristics
for conditions~\eqref{eq:17}--\eqref{eq:45tris},  in the following
lines we  enucleate the main steps   of the proof,  heavily
simplifying most of the technical points and referring to
Section~\ref{s:5}
for all \bigskip details.

\paragraph{\textbf{Sketch of the proof.}} We split the proof in
four \smallskip steps. 

\paragraph{\textbf{Time-discretization.}} Following a well-established routine for
gradient flows  (cf.,\ e.g.,\ \cite{Crandall-Liggett71, Brezis73,
Baiocchi89, Ambrosio95, Rulla96, Savare96, Nochetto-Savare-Verdi00,
Ambrosio-Gigli-Savare08, Rossi-Savare06}),
 and in
general doubly nonlinear equations \cite{Colli-Visintin90, Colli92,
Segatti06, Schimperna-Segatti-Stefanelli07, Mielke-Rossi-Savare08},
we approximate~\eqref{Cauchy-special-form} with the \emph{implicit
Euler scheme}
\begin{equation}
  \label{euscheme}
  \Utau^0:=u_0,\quad
  \partial \Psi\left(\frac{\Utau^n-\Utau^{n-1}}{\tau}\right)+\frsub \ene{t_n}{\Utau^n}\ni 0
  \quad
  n=1,\ldots,N,
\end{equation}
where $\tau=T/N$ is the time step, inducing  a partition of $[0,T]$
with nodes $(t_n:=n\tau)_{n=0}^{N}$. Since \eqref{euscheme} is the
Euler-Lagrange equation for the minimum problem
\[
  U^n_\tau\in \argmin_{U \in \domainenergy}
  \left\{\tau\Psi\left(\frac{U-U^{n-1}_\tau}\tau\right)+\ene{t_n}{U} \right\}
  \quad n=1,\cdots, N,
\]
we look for $(U^n_\tau)_{n=1}^{N}$ solving the above family of
variational problems. Assumption~\eqref{eq:17} yields, via  the
\emph{direct method} in the Calculus of Variations, the existence of
solutions $(U^n_\tau)_{n=1}^{N}$.
 The variational sum
rule~\eqref{2.eq:42-bis} ensures that for every $U^n_\tau$
fulfills~\eqref{euscheme} for all $n=1,\ldots,N$. Hence, we
construct approximate solutions to~\eqref{Cauchy-special-form} by
introducing the  left-continuous piecewise constant $(\pwC
U{\tau})_\tau$ and the piecewise linear $(\pwL U{\tau})_\tau$
interpolants of the discrete solutions $(U^n_\tau)_{n=1}^{N}$ (cf.\
Sec.~\ref{ss:3.1}), which clearly fulfill the approximate equation
\begin{equation}
\label{discrete-eq}
\partial\Psi\left(\pwL
{{\UU}}{\tau}'(t)\right) +\frsub \ene{\pwC
{\mathsf{t}}{\tau}(t)}{\pwC {\UU}{\tau}(t)} \ni 0 \quad \foraa\,  t
\in (0,T),
\end{equation}
where $\pwC {\mathsf{t}}{\tau}$ is the left-continuous piecewise
constant interpolant associated with the partition
$(n\tau)_{n=0}^{N}$
\smallskip
of~$(0,T)$. 

\paragraph{\textbf{Approximate energy inequality and a priori estimates}} In the present
nonconvex setting,   \eqref{discrete-eq} does not yield sufficient
information to pass to the limit and conclude existence for
\eqref{Cauchy-special-form}.  One needs the finer information
provided by the \emph{approximate energy inequality} involving the
Fenchel-Moreau conjugate $\Psi^*$ of $\Psi$, namely
\begin{equation}
\label{dg}
\begin{aligned}
\int_0^T \Psi\left(\pwL {{U}}{\tau}'(t) \right)  \dd t    + \int_0^T
\Psi^*\big({-}\pwM \xi{\tau}(t)\big) \dd t +\ene{T}{\pwC
{\UU}{\tau}(T)} \leq \ene{0}{u_0} + \int_0^T \partial_t \ene{t}{
\pwM U{\tau}(t)}\dd t, \\ \text{with} \ \ \pwM \xi{\tau}(t) \in
\frsub \ene t{\pwM U{\tau}(t))}  \ \ \foraa\, t \in (0,T).
\end{aligned}
\end{equation}
 Here $\pwM U{\tau}$  is the so-called  \emph{De Giorgi variational
 interpolant}
  of the discrete solutions $(U^n_\tau)_{n=1}^{N}$  (see
\eqref{interpmin} for its definition, and Lemma~\ref{lemma:degiorgi}
for its properties). Relying on the positivity of $\Psi$ and of
$\Psi^*$ (cf.\ \eqref{psipos}), on
\eqref{eq:psi-sum-1}--\eqref{eq:psi-sum-2}, on~\eqref{eq:17}, and on
estimate~\eqref{hyp:en3},  from inequality~\eqref{dg} it is possible
to deduce suitable a priori estimates on the sequences $( \pwC
{{U}}{\tau})$, $(\pwL {{U}}{\tau})$, $(\pwM {{U}}{\tau})$, and
$(\pwM {{\xi}}{\tau})$. Then, one infers that,
there exist $u \in \AC ([0,T];\V)$ and $\xi \in L^1(0,T;\V^*)$ such
that, up to a subsequence,
\[
\begin{array}{lll}
& \pwC {{U}}{\tau},\,\pwL {{U}}{\tau},\, \pwM {{U}}{\tau} \to u &
\text{in } L^\infty (0,T;\V),
\\
& \pwL {{U}}{\tau}' \weakto u' & \text{in } L^1 (0,T;\V),
\\
& \pwM {{\xi}}{\tau} \weakto \xi & \text{in } L^1 (0,T;\V^*).
\end{array}
\smallskip
\]

\paragraph{\textbf{Passage to the limit and proof of the upper energy estimate.}} Using
that $\Psi$ and $\Psi^*$ are convex, it is possible to pass to the
limit by lower semicontinuity in \eqref{dg} and conclude that the
functions $u$ and $\xi$ fulfill the \emph{upper energy estimate}
\begin{equation}
\label{dg-ineq}
\begin{aligned}
\int_0^T \Psi(u'(t))  \dd t  + \int_0^T \Psi^*({-}\xi(t)) \dd t
+\ene{T}{u(T)} \leq \ene{0}{u_0} +  \int_0^T \partial_t \ene t{u(t)}
\dd t\,.
\end{aligned}
\end{equation}
Furthermore, the closedness  property~\eqref{eq:45} and an argument
combining Young measures and measurable selection tools yield that
there exists $\xi \in L^1 (0,T;\V^*)$ such that
\begin{equation}
\label{e:conseq-closure} \xi(t) \in \frsub \ene t{u(t)} \ \ \foraa\,
t \in (0,T).
\smallskip
\end{equation}

\paragraph{\textbf{Proof of the  energy identity and conclusion.}} The chain rule
\eqref{eq:45tris} entails
\[
\ene{0}{u_0} +  \int_0^T \partial_t \ene t{u(t)} \dd t =
\ene{T}{u(T)} + \int_0^T \langle -\xi(t), u'(t) \rangle \dd t.
\]
Combining this with~\eqref{dg-ineq}, we ultimately deduce
\[
\int_0^T  \big( \Psi(u'(t))   +  \Psi^*({-}\xi(t)) -\langle
{-}\xi(t), u'(t) \rangle \big)\,  \dd t \leq 0.
\]
Using the Fenchel-Young inequality $\Psi(v) +\Psi^* (w) \geq \langle
w,v \rangle$, we arrive at
\begin{equation}
\label{final-step} \Psi(u'(t))   +  \Psi^*({-}\xi(t)) -\langle
{-}\xi(t), u'(t) \rangle =0 \quad \foraa\, t \in (0,T),
\end{equation}
 whence $-\xi(t) \in \partial \Psi(u'(t)) $
 for almost all $t \in (0,T)$. Combining this with~\eqref{e:conseq-closure},
 we conclude that $u$ is a solution of~\eqref{Cauchy-special-form}.
\QED 
\begin{remark}[Chain-rule inequality]
\label{rem:towards-ineq} \upshape  In view of the discussion
developed in Section \ref{s:3}, let us anticipate that the
chain-rule condition \eqref{eq:45tris} could be weakened. In fact,
 a close perusal at the  argument
for the proof of Theorem~\ref{th:0} reveals that, it is sufficient
to require  the chain-rule inequality
\begin{equation}
\label{inequality-sufficient}
 \frac {\mathrm{d}}{\mathrm{d}t}\ene t{u(t)}  \rangle   \geq \la
 \xi(t),u'(t)\ra  + \partial_t \ene t{u(t)}
      \quad\text{for a.a. }\, t \in  (0,T)\,.
\end{equation}
Indeed, \eqref{inequality-sufficient} yields the \emph{lower energy
estimate}
\[
\ene{T}{u(T)} - \ene{0}{u_0}  \geq
 \int_0^T \langle \xi(t), u'(t) \rangle \dd t  + \int_0^T \partial_t \ene
t{u(t)} \dd t,
\]
which, combined with the \emph{upper energy estimate}
\eqref{dg-ineq}, ultimately yields \eqref{final-step}. In turn, the
\emph{upper energy estimate} \eqref{dg-ineq} is a consequence of the
time-discretization scheme (in particular, of the approximate energy
inequality \eqref{dg}), and of classical lower semicontinuity
results.
\end{remark}
\subsection{Sufficient conditions for closedness and chain rule}
\label{ss:2.2}
 In this section we  revisit the abstract assumptions of
Theorem~\ref{th:0}, and in particular we provide  sufficient
conditions  of $\lambda$-convexity type on the energy functional
$\cE$ for the closedness property~\eqref{eq:45}
 and for the chain rule \eqref{eq:45tris}.

Throughout the following discussion, we will suppose that $\cE$
complies with the time-dependence condition~\eqref{hyp:en3}.
\paragraph{\textbf{Uniformly (Fr\'echet-)subdifferentiable functionals.}} A first sufficient
condition for \eqref{eq:45} and  \eqref{eq:45tris} is some sort of
\emph{uniform subdifferentiability} of the functionals $\cE_t$ on
their sublevels, cf.\ \eqref{unform-modulus} below.

 For every $R>0$ we set
\[
\sublev{R} = \left\{ u \in \domainenergy\, : \ \cg{u} \leq R
\right\}\,.
\]
In view of \eqref{eq2.6}, every  $u \in \domainenergy$ satisfies $u
\in \sublev R$ for some $R>0$.
\begin{proposition}
\label{l:closedness}  Let $\cE : [0,T]\times \V \to
(-\infty,+\infty]$  be a family of time-dependent functionals as
in~\eqref{Ezero}, complying with~\eqref{hyp:en3}. Moreover assume
that   for all $R>0$ there exists a \emph{modulus of
subdifferentiability} $\omega^R:[0,T]\times \sublev{R}\times
\sublev{R} \to [0,+\infty)$   such that for all $t \in [0,T]$:
\begin{equation}
\label{unform-modulus}
\begin{gathered}
 \text{$\omega_t^R(u,u)=0$ for every $u\in
D_R$,}
\\
  \text{the map $(t,u,v) \mapsto \omega_t^R(u,v)$ is upper
semicontinuous, and}
\\
  \ene tv- \ene tu - \la \xi,v-u\ra \geq
  -\omega^R_t(u,v)\Vnorm{v-u} \quad \text{for all $u,\,v \in \sublev R$
and $\xi \in \frsub \ene tu$.}
  \end{gathered}
\end{equation}
Then,   $\cE$ complies with the closedness condition~\eqref{eq:45}
and with the chain rule \eqref{eq:45tris}.
\end{proposition}
\PROOF \textbf{Ad \eqref{eq:45}.} Let $v \in \domainenergy$ be
fixed, and let $(u_n)$, $u$, $(\xi_n)$ and $\xi $ be like in
\eqref{eq:45}.
 It follows from estimate \eqref{eq2.6} that
  there exists some $R>0$ such that $v, \, u_n \in \sublev{R}$ for all $n \in \N$.
Thanks to~\eqref{unform-modulus}, we have
\begin{equation}
\label{forall-n} \ene tv -\ene t{u_n} - \langle \xi_n, v-u_n \rangle
 \geq - \omega_t^{R}(u_n,v)\|v-u_n\| \quad \text{for all $n
\in \N$.}
\end{equation}
Since the functionals $\cE_t$ and  $\omega_t^{R}$ are respectively
lower and upper semicontinuous, we can pass to the limit
in~\eqref{forall-n}, obtaining
\[
\ene tv -\ene t{u} -  \langle \xi, v-u \rangle \geq -
\omega_t^{R}(u,v)\|v-u\|\,.
\]
It is not difficult to check that this inequality yields $\xi \in
\frsub \ene tu$: indeed, notice that, in the definition \eqref{eq:8}
of Fr\'echet subdifferential, it is not restrictive to consider
sequences $(v_k)_k$ converging to $u$, such that $\limsup_{k \to
\infty} \ene t{v_k} \leq \ene tu$.
 Furthermore, choosing $v=u$ in \eqref{forall-n}
(notice that $u \in \sublev{R}$ by lower semicontinuity),
 and exploiting the properties of $\omega_t^{R}$,  we
have the following chain of inequalities
\[
\begin{aligned} 0 \leq \limsup_{n \to \infty}  \left(\ene t{u_n}  {-} \ene tu
\right)   \leq \lim_{n \to \infty} \langle \xi_n, u_n{-}u \rangle +
\limsup_{n \to \infty} \omega_t^{R}(u_n,u)\|u-u_n\| =0,
\end{aligned}
\]
whence $\ene t{u_n} \to \ene tu$, which concludes the proof
of~\eqref{eq:45}.
\\
\textbf{Ad \eqref{eq:45tris}.} Let $u \in \AC ([0,T];\V)$ and $\xi
\in L^1 (0,T;\V^*)$ fulfill \eqref{conditions-for-chain-rule}. Up to
a suitable reparametrization (cf.\
\cite[Lemma\,1.1.4]{Ambrosio-Gigli-Savare08}), it is possible to
assume that the curve $u$ is $1$-Lipschitz. Furthermore, due to
$\sup_{t \in [0,T]} \ene t{u(t)}<+\infty$  there
exists $R>0$ such that  $u(t) \in \sublev{R}$ for all $t \in [0,T]$.
In order to show the absolute continuity of the map $t \mapsto \ene
t{u(t)}$, we estimate for $0 \leq s \leq  t \leq T$ the difference
\[
\ene t{u(t)} - \ene s{u(s)} = \ene t{u(t)} - \ene s{u(t)} + \ene
s{u(t)} - \ene s{u(s)}.
\]
First, due to \eqref{hyp:en3}, for all $0 \leq s \leq t\leq T$ there
holds
\begin{equation}
\label{first-term} \ene t{u(t)} - \ene s{u(t)} = \int_s^t \partial_t
\ene r{u(t)} \, \d r \leq C_1 \cg{u(t)} (t-s) \leq C_1 R (t-s).
\end{equation}
Second,  in view of \eqref{unform-modulus} we have
\begin{equation}
\label{second-term} \ene s{u(t)} - \ene s{u(s)} \geq \langle \xi(s),
u(t) - u(s) \rangle - \omega^R_s(u(s),u(t)) \| u(t) - u(s)\|.
\end{equation}
Combining \eqref{first-term} and \eqref{second-term}, inverting the
role of $s$ and $t$, and using the $1$-Lipschitz continuity of $u$,
we conclude
\begin{equation}
\label{as-abs-cont} |\ene t{u(t)} - \ene s{u(s)}| \leq (C_1 R +
\|\xi(t)\|_* + \|\xi(s)\|_* + \omega^R_t(u(s),u(t)) +
\omega^R_s(u(s),u(t)))|t-s|.
\end{equation}
Now, the  upper semicontinuity of   $(t,u,v) \mapsto \omega_t^R(u,v)$, joint with the fact that $u
\in \mathrm{C}^0 ([0,T];\V)$, yield that there exists $C>0$ such that
  there holds
 $\omega^R_t(u(s),u(t)), \, \omega^R_s(u(s),u(t)) \leq C$ for all $s,\,t \in [0,T]  $. Therefore,
arguing as in \cite[Thm.\,1.2.5]{Ambrosio-Gigli-Savare08}, we
conclude that $t \mapsto \ene t{u(t)}$ is absolutely continuous. Let
us now fix a point $t \in (0,T)$ such that $u'(t)$, $\frac{\d}{\d t
} \ene t{u(t)}$ exist: arguing throughout
\eqref{first-term}--\eqref{second-term},   it is not difficult to
deduce that
\[
\ene{t+h}{u(t+h)} - \ene t{u(t)} \geq \int_t^{t+h} \partial_t \ene
r{u(t)} \, \d r + \langle \xi(t), u(t+h) - u(t) \rangle -
\omega^R_t(u(t),u(t+h)) \| u(t+h) - u(t)\|.
\]
Dividing by $h>0$ and $h<0$ and taking the limit as $h \down 0$ and
$h \up 0$, we prove  the chain rule \eqref{eq:45tris}.
 \QED
 \begin{remark}[$\lambda$-convex functionals]
\upshape A sufficient condition yielding~\eqref{unform-modulus} in
the case of an energy functional $\cE: [0,T]\times \V \to
(-\infty,+\infty]$ complying with \eqref{Ezero} and \eqref{hyp:en3},
is that the map $u \mapsto \cE_t(u)$ is \emph{$\lambda$-convex}
uniformly in $t \in [0,T]$, namely
\begin{equation}
\label{e:V-convexity}
\begin{aligned}
\exists\, \lambda \in \R \,  \quad &  \forall\, t \in [0,T] \
\forall\,u_0,\, u_1 \in \domainenergy\ \forall\, \theta \in [0,1]\, : \\
& \ene t{(1-\theta)u_0+\theta u_1}\le (1-\theta)\ene t{u_0}+
  \theta\ene t{u_1}-\frac\lambda2\theta(1-\theta)\|u_0-u_1\|^2\,.
  \end{aligned}
\end{equation}
Indeed,  given $u,\, v \in \domainenergy$, \eqref{e:V-convexity} and
the very definition \eqref{eq:8} of Fr\'echet subdifferential yield
 for any $\xi \in \frsub \ene{t}{u}$ and $\theta \down 0$
\[
\begin{aligned}
\theta \left( \ene{t}{v} -\ene{t}{u}\right)& \geq \frac{\lambda}2
\theta (1-\theta) \Vnorm{v-u}^2 +  \ene t{(1-\theta)u+\theta v} -
\ene{t}{u}\\ & \geq \theta \left( \frac{\lambda}2  (1-\theta)
\Vnorm{v-u}^2 + \langle \xi,  v-u\rangle + o (1) \right)\,.
\end{aligned}
\]
Upon diving by $\theta$, we conclude inequality
 \eqref{unform-modulus} with the choice $\omega_t(u,v):= \frac
{\lambda^-}2 \|u-v\|$.
\end{remark}
\begin{remark}[Perturbations of $\lambda$-convex functionals]
\upshape
 In
\cite{Rossi-Mielke-Savare08} a broad family of time-dependent
energies, which for instance encompasses $\lambda$-convex
functionals, was tackled.  However,  as hinted in the Introduction,
\cite{Rossi-Mielke-Savare08}  focuses on the analysis (from a metric
viewpoint) of  doubly nonlinear equations driven by a less general
class of dissipation potentials than those considered in the present
paper.
 While referring to \cite{Rossi-Mielke-Savare08} for details, here we  recall
that the energies therein  considered are given by the sum of two
time-dependent functionals $\cE^1,\, \cE^2 : [0,T]\times \V \to
(-\infty,+\infty]$, such that the functionals $\cE^1_t$ are
$\lambda$-convex, uniformly with respect to $t \in [0,T]$, and the
functionals $\cE^2_t$ are \emph{dominated concave perturbations} of
 $\cE^1_t$ (cf.\ \cite[Definitions
5.4, 5.10]{Rossi-Mielke-Savare08}, as well as \cite{Rossi-Savare06}
for an analogous class of functionals). In \cite[Propositions 5.6,
5.10]{Rossi-Mielke-Savare08} it was shown that, under the above
conditions, the energy $\cE=\cE^1+ \cE^2 $ fulfills the closedness
condition \eqref{eq:45} and the chain rule \eqref{eq:45tris}.
\end{remark}
%
\subsection{Examples and applications}
\label{ss:2.3}
\begin{example}[A model in ferro-magnetism]
\label{ex:ferro-magnetism} \upshape We take $\Omega \subset \R^3$ a
bounded sufficiently smooth domain,  and let
\[
V=L^2 (\Omega;\R^3) \qquad \text{and} \ \ \Psi(v)=\int_{\Omega}
|v|+\frac12 |v|^2 \dd x =\|v\|_{L^1(\Omega;\R^3)} +  \frac12
\|v\|_{L^2(\Omega;\R^3)}^2.
\]
We  consider a simplified model for ferro-magnetism, in which the
interplay between  the elastic and the magnetic effects is neglected
(see~\cite[Sec.~7.4]{Mielke05} for a rate-independent model
accounting for both features). In this framework,  the relevant
energy functional $\cE: [0,T] \times L^2 (\Omega;\R^3) \to
(-\infty,+\infty]$  is given by
\begin{equation}
\label{e:magnetization} \begin{aligned} &\ene{t}{m} = \left\{
\ba{ll} \!\!\!  \int_{\Omega}\Big( \frac{\alpha}2|\nabla m|^2 + W(m)
\Big)\dd x +\int_{\R^3}|\nabla \Phi_m|^2\dd x -\pairing{}{H^1}{
H_{\textrm{ext}}(t)}{ {m}}  & \text{if $m \in H^1(\Omega;\R^3), $}
\\
 \!\!\! +\infty & \text{otherwise.} \ea \right.
\end{aligned}
\end{equation}
Here,
$\pairing{}{H^1}{\cdot}{\cdot}$ is a short-hand notation for the
duality pairing between $H^1(\Omega;\R^3)^*$ and
$H^1(\Omega;\R^3)$,
  $m : \Omega \to \R^3$ is the magnetization,
\begin{equation}
\label{e:pot-w}
\begin{gathered}
\text{ $W \in \rmC^1 (\R^3;\R)$ a $\lambda_W$-convex potential for
some $\lambda_W \in \R$, such that}
\\
\exists\, c_W, \, C_W>0 \ \forall\,m \in \R^3 \, : \ \ W(r) \geq c_W
|m|^2 - C_W
\end{gathered}
\end{equation}
(e.g.,\ $W(m)= (1-|m|^2)^2)$, and  the external magnetic field
$H_{\textrm{ext}}$ fulfills
\begin{equation}
\label{ext-magnetic-field} H_{\textrm{ext}} \in \rmC^1 ([0,T];H^1
(\Omega;\R^3)^*), \qquad \div(H_{\textrm{ext}}(t))\equiv0,
\end{equation}
 the latter equation meaning that
 $\int_{\Omega} H_{\textrm{ext}}(t) \cdot \nabla v \dd x =0$ for all $v \in
 H^1(\Omega;\R^3)$.
In \eqref{e:magnetization}, the potential $\Phi_m$ describes the
field induced by the magnetization inside the body. Hence, the
magnetic flux is $J= (H_{\textrm{ext}}-\nabla \Phi_m +
\mathbb{E}_\Omega (m))$, $\mathbb{E}_\Omega (m)$ denoting the
trivial extension of $m$ to all of $\Omega$ by $0$. Thus, $\div J=0$
and \eqref{ext-magnetic-field} yield that $\Phi_m$ is  the solution
of
\[ \div\left(-\nabla\Phi_m + \mathbb{E}_\Omega (m) \right)=0 \quad
\text{in $\R^3$.}
\]
Note that the operator $\mathcal{J}: L^2 (\Omega;\R^3) \to L^2
(\R^3;\R^3)$ mapping $m \mapsto \nabla\Phi_m$
 is bounded and linear; it was proved in~\cite{DeSimone93} that
 $m \mapsto \mathcal{J}(m)\Restr{\Omega}$ is an orthogonal projection on
 $ L^2 (\Omega;\R^3)$, satisfying
\begin{equation}
\label{orthogonal} \int_{\R^3}|\nabla \Phi_m|^2 \dd x= \int_{\Omega}
m \cdot \mathcal{J}(m)\Restr{\Omega} \dd x \qquad \text{for all } m
\in L^2 (\Omega;\R^3)\,.
\end{equation}

One can see that for all $ (t,m)  \in [0,T]\times H^1 (\Omega;\R^3)
$ such that $\frsub \ene tm \neq \emptyset$ there holds
\[
\frsub \ene tm =\left\{ {-}\Delta m + \rmD W(m) +
\mathcal{J}(m)\Restr{\Omega} \right\}\,.
\]
Therefore, with the present choices of $\Psi$ and $\cE$, the Cauchy
problem  \eqref{Cauchy-special-form} translates into
\begin{equation}
\label{magnetization-cauchy} \mathrm{Sign}(\dot{m}) + \dot{m}
-\alpha \Delta m + \rmD W(m) + \mathcal{J}(m)\Restr{\Omega} =
H_{\textrm{ext}} \quad \text{a.e.\ in $\Omega \times (0,T)$;} \ \
m(0)=m_0,
\end{equation} with variational boundary conditions;
here
\[
\mathrm{Sign} : \R^3 \rightrightarrows \R^3 \ \ \text{ is given by }
\ \ \mathrm{Sign}(v) = \left\{
\begin{array}{lll}
 \displaystyle \frac v{|v|}  & \text{if $v \neq 0$,}
\medskip
\\
  \overline{B_1(0)}  & \text{if $v = 0$.}
\end{array}
 \right.
\]

Now,
 combining~\eqref{e:pot-w}--\eqref{ext-magnetic-field}, it is
easy to see that $\cE$  complies with~\eqref{eq:17} and
\eqref{hyp:en3}. Also using~\eqref{orthogonal} and arguing as
in~\cite[Sec.~7.2]{Rossi-Mielke-Savare08}, we further check that for
some $\lambda \in \R$ the energy  $\cE$ is $\lambda$-convex (uniformly in $t
\in [0,T]$), with respect to the $L^2 (\Omega;\R^3)$-norm. Therefore, the
closedness property \eqref{eq:45} and, a fortiori, the variational
sum rule~\eqref{2.eq:42-bis} (cf.\
Proposition~\ref{l:variational-sum-rule}) hold, as well as the chain
rule \eqref{eq:45tris}.
 Thus, it follows from Theorem~\ref{th:0}
that the Cauchy problem \eqref{magnetization-cauchy}
admits a solution. 
\end{example}
\begin{example}[Doubly nonlinear evolutions of Allen-Cahn type]
\upshape Let us consider the  following class of evolution equations
\begin{equation}
\label{allen-cahn} \varrho \mathrm{Sign} (\dot{u}) +
|\dot{u}|^{p-2}\dot{u}
 -\div(\beta(\nabla u))
+W'(u) = \ell \quad \text{in $\Omega \times (0,T),$}
\end{equation}
with $\varrho >0$, $1<p<\infty$,  and $u: [0,T]\times \Omega \to
\R$, where
 $\Omega \subset \R^d$, $d \geq 1$, is
a  sufficiently smooth bounded domain.  In \eqref{allen-cahn},
$\beta: \R^d \to \R^d$ is the gradient of some smooth function $j$
on $\R^d$, $W: \R \to \R$ a differentiable function and $\ell:
\Omega \times (0,T) \to \R$ some source term. To fix ideas (cf.\
\cite[Sec.~8.2]{Rossi-Mielke-Savare08} for the precise statement of the
assumptions on $j$ and $W$), we may think of the case in which
$j(\zeta)=\frac1q |\zeta|^q$ for some $q>1$  (hence $\beta(\zeta)=
|\zeta|^{q-2} \zeta$ and the elliptic operator in~\eqref{allen-cahn}
is indeed the $q$-Laplacian), and    $W$ is given by the sum of a
convex function, perturbed by a nonconvex nonlinearity which
complies with suitable growth conditions (like for instance in the
classical, double-well potential case $W(u):= (u^2 -1)^2/4$).

We supplement equation~\eqref{allen-cahn} with homogeneous Dirichlet
boundary conditions, and notice that this boundary-value problem can
be written
 in the abstract
form~\eqref{Cauchy-special-form},  in the framework of the ambient
space
\begin{equation}
\label{ambient-lp}
\text{$\V = L^p(\Omega)$, with the dissipation potential } \
\Psi(v) = \varrho \,\| v \|_{L^1 (\Omega)} + \frac1p \|v \|_{L^p
(\Omega)}^p
\end{equation}
and the driving energy functional
\begin{equation}
\label{driving-energy} \ene tu  = \left\{ \begin{array}{ll}
\int_{\Omega} \left( j(\nabla u(x)) + W (u(x)) \right) \, {\rm d}x
-\pairing{}{W_0^{1,q}}{\ell(t)}{u} & \text{if $u \in
W_0^{1,q}(\Omega)$, $ W(u) \in L^1 (\Omega)$,}
\\
+\infty & \text{otherwise,}
\end{array}
\right.
\end{equation}
 where $\pairing{}{W_0^{1,q}}{\cdot}{\cdot}$ stands for the duality pairing between $W^{-1,q'}(\Omega)$
 and $W_0^{1,q}(\Omega)$,  with $q'=q/(q-1)$.

 In~\cite[Sec.~8.2]{Rossi-Mielke-Savare08}, under suitable
conditions on $j$ and $W$
 the existence of a solution  for the initial-boundary value
problem for~\eqref{allen-cahn}  was proved
 for
 $\varrho=0$. Namely, in
 \cite{Rossi-Mielke-Savare08} only the case of a
  dissipation potential $\Psi$
 induced by the \emph{single norm} $ \|\cdot \|_{L^p(\Omega)}$ was considered,
which does not include the more physical
   form \eqref{ambient-lp}.

  Relying on the analysis
 of~\cite{Rossi-Mielke-Savare08}, it can be checked that,
if $\ell \in \mathrm{C}^1
 ([0,T];W^{-1,q'}(\Omega))$, then
  the energy functional
 $\cE$~\eqref{driving-energy} complies
 with~\eqref{eq:17}--\eqref{eq:45}. Hence,
 Theorem~\ref{th:0} applies,  yielding
the existence of a solution $ u \in L^\infty (0,T;W_0^{1,q}(\Omega))
\cap W^{1,p}(0,T;L^p(\Omega))$ to the initial-boundary value problem
for \eqref{allen-cahn}.
\end{example}
\section{Motivating examples for marginal subdifferentials}
\label{s:2new}
 In this section,
we restrict  to a finite-dimensional setting and give an outlook to
a twofold generalization of the set-up considered in
Section~\ref{s:2}. Such an extension is motivated by the analysis of
abstract evolutionary systems of the form
\begin{equation}
\label{general-coupled} \left.
\begin{array}{lll}
&
\partial\Psi (u'(t)) + \frsub \cE^1 (u(t)) +  \rmD_u \marginal t{u(t)}{\margvar(t)}
\ni \ell(t) & \text{ in $V^*$}
\\
& \rmD_\margvar \marginal t{u(t)}{\margvar(t)}=0 & \text{ in $X^*$}
 \end{array}
 \right\} \quad \foraa\, t \in (0,T),
\end{equation}
where $\cE^1 : \V \to (-\infty,+\infty]$ is a convex energy,
perturbed by some smooth functional $I : [0,T]\times X \times \V \to
\R$ (where $X$ is a second Banach space), and
 $\ell
:[0,T] \to \V^*$  is the external loading.
 Couplings like \eqref{general-coupled}
arise in the modeling of physical systems described in terms of two
variables $(\margvar,u)$, such that energy   dissipation  only
occurs through the   \emph{internal variable} $u$, and $\margvar$
fulfills some stationary law. PDE systems of the type
\eqref{general-coupled} typically arise in connection with
rate-independent behavior (cf.\ \cite{Mielke05} and the references
therein). Nonetheless,  they can also occur in the modeling of
\emph{rate-dependent} evolutions, like for instance in the case of
\emph{quasistationary phase-field models}, cf.\
\cite{Luckhaus90,Schatzle00,Visintin96,Rossi-Savare06}.
 In  Section~\ref{s:4} later on, we analyze
a PDE system of the type \eqref{general-coupled} in finite-strain
elasticity.

Let us observe that the  second stationary relation in
\eqref{general-coupled} is the Euler-Lagrange equation for the
minimum problem $\inf_{\margvar \in X} I_t (u(t),\margvar)$,  and
suppose that
\[
M(t,u) : = \argmin_{\margvar \in X} \marginal t u{\margvar} \neq
\emptyset \quad \text{for all $(t,u) \in [0,T]\times \V$}.
\]
Hence,
 we introduce the
\emph{reduced energy} functionals $\cE^2, \, \cE: [0,T]\times \V\to
(-\infty,+\infty] $
\begin{equation}
\label{marginal-energies} \left\{
\begin{array}{ll}
\enei 2tu:= \min_{\margvar \in X} \marginal t u{\margvar},
\\
\ene tu:= \cE^1 (u)  + \enei 2tu - \langle \ell(t), u \rangle =
\min_{\margvar \in X} \left(\cE^1 (u) +  \marginal t u{\margvar} -
\langle \ell(t), u \rangle \right).
\end{array}
\right.
\end{equation}

In this setting, it is natural to introduce the following
subdifferential notion for energy $\cE^2$, tailored to its reduced
form.
\begin{definition}[Marginal subdifferential]
\label{def-margsubdif} The \emph{marginal subdifferential} (with respect to
the variable $u$) of the reduced functional $\cE^2: [0,T]\times \V
\to (-\infty,+\infty]$ at $(t,u) \in [0,T]\times \V$ is
\begin{equation}
\label{tailored_notion} \margsub \enei 2tu:=  \{ \rmD_u \marginal t
u{\margvar}\, : \ \margvar \in M(t,u) \}.
\end{equation}
\end{definition}
Hereafter, we will address  the doubly nonlinear evolution equation
\begin{equation}
\label{tailored-equation}
\partial \Psi (u'(t)) + \frsub \cE^1 (u(t)) + \margsub \enei 2t{u(t)}
\ni \ell(t) \ \ \text{in $\V^*$} \quad \foraa\, t \in (0,T).
\end{equation}
Clearly, solutions to \eqref{tailored-equation} are in fact
solutions to the quasi-stationary evolution system
~\eqref{general-coupled}.

Notice that \eqref{tailored-equation} may be viewed as a
generalization of the doubly nonlinear equations, featuring the
Fr\'echet subdifferential, examined in Section \ref{s:2}. Indeed,
 under quite standard assumptions  there holds
\begin{equation}
\label{frechet-subdiff-inclusion} \frsub \ene t{u} \subset \frsub
\cE^1 (u) + \margsub \enei 2t{u} - \ell(t) \quad \text{for all }
(t,u)\in [0,T]\times \V,
\end{equation}
while the converse inclusion is not true, in general.

This fact is illustrated in Example~\ref{ex:3.0}: for a specific
choice of the functional $\cE^1$ and for a \emph{time-independent}
marginal functional $\cE^2$, it is shown that the Fr\'echet
subdifferential  of the energy $\ene tu = \cE^1 (u) + \cE^2 (u)
-\langle \ell(t), u \rangle $ in \eqref{marginal-energies} does not
 comply with the closedness property \eqref{eq:45}. Furthermore, the
closure of $\partial \cE$ in the sense of graphs coincides for all
$(t,u) \in [0,T] \times \V$ with the set  $\frsub \cE^1 (u) +
\margsub \cE^2 (u) -\ell(t)$.

Another important feature which sets aside reduced energy
functionals from the class of energies examined in Section \ref{s:2}
is that,  even if the function $t\mapsto \marginal t u{\margvar}$ is
smooth, the resulting reduced functionals $\cE^2$ and $\cE$ (cf.\
\eqref{marginal-energies}) may be nonsmooth with respect to time,
see  Examples \ref{ex:3.1} and~\ref{ex:3.2}. Therein, we suggest the
usage of a   \emph{generalized time-derivative}, defined in such a
way as to comply with a suitable chain-rule inequality.

\begin{example}
\label{ex:3.0} \upshape For simplicity we restrict to the
one-dimensional case $\V=X=\R$, and to gradient flows, hence taking
$\Psi(v) = \frac12 |v|^2$. As in \cite[Ex.~2]{Rossi-Savare06}, we
choose
\[
\cE^1 (u) = \frac12 |u|^2, \ \ \marginal t u{\margvar}= \marginal {}
u{\margvar} = \frac12 |\margvar|^2 - u\margvar +
\mathcal{W}(\margvar),
\]
where $\mathcal{W} : \R \to \R$ is the (piecewise quadratic) double
well potential
 \begin{equation}
      \mathcal{W}(\margvar):=
      \left\{
      \begin{array}{lll}
      (\margvar+1)^2
      & \margvar<- \frac{1}{2} ,\\
      -\margvar^2+ \frac{1}{2}
      & |\margvar| \leq \frac{1}{2},
      \\
      (\margvar-1)^2
    & \margvar > \frac{1}{2},
  \end{array}
  \right.
  \ \
  \text{with derivative}\ \
    \mathcal{W}'(\margvar)=
    \left\{
  \begin{array}{lll}
    2(\margvar+1)
    & \margvar <- \frac{1}{2} ,\\
    -2 \margvar    & |\margvar| < \frac{1}{2},
    \\
    2(\margvar-1)
    & \margvar > \frac{1}{2}.
 \end{array}
  \right.
  \label{eq:dwpot}
  \end{equation}
In this setting, given some smooth external loading $\ell: [0,T]\to
\R$, the coupled system \eqref{general-coupled} reads
\[
\left\{
\begin{array}{ll}
 u'(t) + u(t)-\margvar(t)=\ell(t),\\
 \mathcal{W}'(\margvar(t)) + \margvar(t) =u(t)
          \end{array}
          \right. \quad \foraa\, t \in (0,T),
\]
which may be viewed as  the one-dimensional caricature of the \emph{quasistationary
phase-field system} (cf.\ \cite{Luckhaus90, Schatzle00, Visintin96}).

It was observed in \cite[Sec.~2.1]{Rossi-Savare06} that the
Fr\'echet subdifferential of the energy functional $\cE$ defined
in~\eqref{marginal-energies} satisfies
\eqref{frechet-subdiff-inclusion}. More precisely,
\begin{equation}
\label{first-subdiff}
 \frsub \ene tu \neq \emptyset \ \Rightarrow \ \ \left\{
\begin{array}{ll}
   \text{$\frsub \ene
tu $ and $M(t,u) $ reduce to a  singleton, and } \\ \frsub \ene tu = \frsub  \cE^1 (u) + \margsub
\enei 2{}{u} - \ell(t) =
\frsub \cE^1 (u) - M(t,u) - \ell(t).
\end{array} \right.
\end{equation}
Now, since the subdifferential mapping $\frsub \cE_t : \R
\rightrightarrows \R$ is not closed in the sense of graphs, it is
natural to introduce its closure, i.e.\ the \emph{limiting subdifferential}
 (cf.\ \cite{Mordukhovich84,
Mordukh06}, and \cite{Rossi-Savare06, Rossi-Segatti-Stefanelli08,
Rossi-Segatti-Stefanelli09} for some analysis of gradient flow and
doubly nonlinear equations featuring such a notion of
subdifferential), defined by
\[
\lmsub \ene tu := \{\xi \in \R\,: \ \exists\, (u_n),\, (\xi_n)
\subset \R, \ u_n \to u, \ \xi_n \to \xi, \ \ene t{u_n} \to \ene tu
\}.
\]
From the closedness of the graph of the  multivalued mapping
$M(t,\cdot) : \R \rightrightarrows \R$ we infer that a weaker form
of \eqref{first-subdiff} passes to the limit, i.e.
\begin{equation}
\label{second-subdiff} \lmsub \ene tu  \subset \frsub  \cE^1 (u) -
M(t,u) - \ell(t) =\frsub  \cE^1 (u) + \margsub \enei 2{}{u} -
\ell(t) \quad \text{for every $(t,u) \in [0,T]\times \R$.}
\end{equation}
In fact, in the case of \eqref{eq:dwpot} we even have $\lmsub \ene
tu  = \frsub  \cE^1 (u) + \margsub \enei 2{}{u} - \ell(t)$.
Relations \eqref{first-subdiff} and \eqref{second-subdiff} suggest
the choice of the subdifferential notion $\diff tu :=\partial \cE^1
(u) + \margsub \enei 2{}{u} - \ell(t)$ for reduced energies of the
type \eqref{marginal-energies}. We explore this viewpoint in
Section~\ref{s:4}.
\end{example}
\begin{example}
\label{ex:3.1}
 We take $V = \R$, $\Psi(v) = \frac12 |v|^2$ and, we set
\[
\ene t {u} = - \alpha |u - \beta t| \qquad \text{for all $u \in \R$,
$t \in (0,T)$, with } \alpha > \beta >0, \ \beta <1\,.
\]
Note that $\cE $ is a marginal function: indeed,
\[
\ene t {u} = \min_{\margvar \in \{ 0,1\}} \marginal t u{\margvar}, \
\text{ with }
 \
\marginal t u{\margvar} = \left\{
 \begin{array}{lll} - \alpha u + \alpha \beta t & \text{if $\margvar=0$,}
\\
\alpha u - \alpha \beta t & \text{if $\margvar=1$.}
\end{array}
\right.
 \]
In this case, $\cE$ does not comply with the smooth time-dependence
condition \eqref{hyp:en3}, and it is only Lipschitz continuous with
respect to both variables $t$ and $u$. It is then natural to
consider the
 \emph{Clarke} subdifferentials of the energy $\cE$ with respect to $u$ and
 $t$, which
are easily calculated:
\begin{equation}
\label{clarke}
\partial_u^{\mathrm{Clarke}}\ene t{u} =
 \left\{
\begin{array}{lll}
\{-\alpha\} & \text{if $u>\beta t$,}
\\
 {[-\alpha, \alpha]} & \text{if $u = \beta t$,}
\\
\{\alpha\} & \text{if $u<\beta t$,}
\end{array}
\right.
 \quad
 \partial_t^{\mathrm{Clarke}}\ene t {u}=
 \left\{
\begin{array}{lll}
\{\alpha\beta\} & \text{if $u>\beta t$,}
\\
 {[-\alpha\beta,
\alpha\beta]} & \text{if $u = \beta t$,}
\\
\{-\alpha \beta\} & \text{if $u<\beta t$.}
\end{array}
\right.
\end{equation}
Notice that  $\frsub \ene tu \subset
\partial_u^{\mathrm{Clarke}}\ene t {u}$ for all $(t,u) \in [0,T]\times
\R$. Furthermore, the multivalued mapping
$\partial_u^{\mathrm{Clarke}}\cE_t: \R \rightrightarrows \R$ is
closed in the sense of graphs. We may choose $\diff tu : =
\partial_u^{\mathrm{Clarke}}\ene tu$ and  consider the gradient flow
\begin{equation}
\label{clarke-gflow} u'(t) + \partial_u^{\mathrm{Clarke}}\ene t
{u(t)} \ni 0 \qquad \foraa\, t \in (0,T)\,.
\end{equation}

 We immediately verify that the
 curve
$\bar{u}:[0,T]\to \R$ defined by
 $ \bar{u}(t) = \beta t $
  is a solution of \eqref{clarke-gflow}.  Now, we aim to get some insight into a possible surrogate notion of
chain rule in this nonsmooth setting.  Imposing that the chain-rule
inequality~\eqref{inequality-sufficient} holds along the curve
$\bar{u}$, with the Clarke subdifferentials~\eqref{clarke},  we
arrive at
\begin{equation}
\label{e:crule} \frac{\mathrm{d}}{\mathrm{d}t} \ene t {\bar{u}(t)}
\geq \xi \bar{u}'(t) + p \ \  \text{for all $\xi \in F_t
(\bar{u}(t))$, $p \in\partial_t^{\mathrm{Clarke}}\ene t
{\bar{u}(t)}$, and} \ \foraa\, t\in (0,T).
\end{equation}
Since $\ene t {\bar{u}(t)} \equiv 0$ and $\bar{u}'(t) \equiv \beta$,
this amounts to checking if
\[
0 \geq \xi \beta + p \qquad \text{for all $\xi \in [-\alpha, \alpha]
$, $p \in [-\alpha\beta, \alpha\beta]$,}
\]
which does not hold.

However, it is true that for a fixed $\xi \in F_t (\bar{u}(t)) $
there exists  a set $\widehat{\mathrm{P}}_t ({\bar{u}(t)},{\xi})$
such that inequality~\eqref{e:crule} holds for $\xi$ and all
elements $p \in \widehat{\mathrm{P}}_t ({\bar{u}(t)},{\xi})$, namely
\begin{equation}
\label{new-pt} \widehat{\mathrm{P}}_t ({\bar{u}(t)},{\xi})=
[-\alpha\beta,-\xi\beta]\,.
\end{equation} Finally, we may observe that, if we ask for an equality sign in~\eqref{e:crule}
for a fixed $\xi \in F_t (\bar{u}(t)) $, then the corresponding set
$\widehat{\mathrm{P}}_t({\bar{u}(t)},{\xi})$ reduces to the
singleton $\{-\xi\beta\}$, cf.\ also Remarks \ref{rmk:example-p} and
\ref{rmk:versus} for further related comments.
\end{example}
\begin{example}
\label{ex:3.2} We reconsider the triple $(\V,\ene tu,\Psi)= (\R,
\min_{\margvar \in \{0,1\}} \marginal t u \margvar, \frac12
|\cdot|^2)$ of Example \ref{ex:3.1}, but choose for  $\diffname$ the
marginal subdifferential of $\cE$ (cf.\ Definition
\ref{def-margsubdif}), viz.
\[
\diff tu = \margsub \ene t u = \left\{
\begin{array}{lll}
\{-\alpha\} & \text{if $u > \beta t$}
\\
 \{-\alpha,\alpha \} & \text{if $u = \beta t$}
\\
\{\alpha\} & \text{if $u < \beta t$.}
\end{array}
 \right.
\]
Notice that, in this case,  the
 curve
$\bar{u}:[0,T]\to \R$ defined by
 $ \bar{u}(t) = \beta t $ is not a solution of the gradient flow $
 u'(t) + \margsub \ene t
{u(t)} \ni 0$ on $(0,T)$. Imposing that the
 chain-rule
inequality~\eqref{inequality-sufficient} holds along the curve
$\bar{u}$, for  the marginal subdifferential with respect to $u$ and the
Clarke subdifferential with respect to $t$, yields
\[
0 \geq \xi \beta + p \qquad \text{for all $\xi \in \{-\alpha,
\alpha\} $, $p \in [-\alpha\beta, \alpha\beta]$.}
\]
Thus, referring to notation \eqref{new-pt}, we conclude that, in this case,
\[
\widehat{\mathrm{P}}_t ({\bar{u}(t)},{\xi})=
[-\alpha\beta,-\xi\beta]= \left\{
 \begin{array}{lll}
 \{ -\alpha\beta \}
& \text{if $\xi=\alpha$}
\\
 {[-\alpha\beta, \alpha\beta]} & \text{if $\xi=-\alpha.$}
\end{array}
 \right.
\]
\end{example}

 Examples \ref{ex:3.1} and \ref{ex:3.2} seem to suggest that,
  to deal with marginal functions, one
should use a notion of time-derivative  $\Ptname$ \emph{conditioned,
via the chain rule, to elements $\xi$ of the  subdifferential}. This
means that, in addition to the $(t,u)$-dependence,  such a notion
$\Ptname$ also depends on the elements $\xi \in \diff tu$. This is
the point of view we are going to adopt in what follows.

\section{Main results}
 \label{s:3}
\subsection{Assumptions}
\label{ss:3.0} We recall that
  $\V$ is  a reflexive separable Banach space.
Below we  enlist our general assumptions on the
\emph{state-dependent} dissipation $\Psi= \Psi_u(v)$,  and
 on the energy functional
$\cE: [0,T] \times \V
 \to
(-\infty,+\infty]$, with domain $[0,T]\times\domainenergy$.  We
emphasize that the conditions on $\cE$ involve \emph{both}     its
subdifferential
   $\diffname: [0,T] \times \domainenergy \rightrightarrows
\V^* $ (with domain  and graph $\domaindiff$ and
$\mathrm{graph}(\diffname)$, respectively), \emph{and} its
\emph{generalized partial time-derivative} $\Ptname = \Pt tu{\xi}$,
for $(t,u,\xi) \in \mathrm{graph}(\diffname)$,  since we encompass a
nonsmooth dependence of the energy $\cE$ on the time
\medskip variable.

\paragraph{\textbf{A (Finsler) family of   dissipation
potentials.}}
 \noindent We consider a family
\begin{equation}
   \label{e:4.1}
   \tag{$4.{\Psi_1}$}
  \begin{gathered}
  \text{$\Psi_u:\V\to [0,+\infty)$, $u \in \domainenergy$, of \emph{admissible dissipation
  potentials}}
  \end{gathered}
\end{equation}
i.e.\ $\Psi_u$ complies with
\eqref{eq:psi-sum-1}--\eqref{eq:psi-sum-2} for all $u \in
\domainenergy$.   We now require that the potentials $(\Psi_u)_{u
\in \domainenergy}$ and $(\Psi_u^*)_{u \in \domainenergy}$ have a
superlinear growth, uniformly with respect to $u$ in sublevels of the energy
$\cE$, viz.
\begin{equation}
  \label{eq:41.1}
  \tag{$4.{\Psi_2}$}
 \forall\, R>0\,: \ \ \left\{
\begin{array}{ll}
\displaystyle \lim_{\|v \|\to +\infty} \frac{1}{\|v\|} \inf_{\cg u
\leq R} \Psi_u(v) =+\infty,
\medskip
\\
\displaystyle \lim_{\|\xi \|_*\to +\infty} \frac{1}{\|\xi\|_*}
\inf_{\cg u \leq R} \Psi_u^*(\xi) =+\infty,
\end{array}
\right.
\end{equation}
where we have used the notation $\cg u = \sup_{t \in [0,T]}\ene
tu$.
 Furthermore, we  require that the dependence $u \mapsto
\Psi_u$ is continuous, on sublevels of the energy,  in the sense of
\textsc{Mosco}-convergence (see, e.g, \cite[\S\,3.3,
p.\,295]{Attouch84}), i.e.
\begin{equation}
\label{mosco} \tag{$4.{\Psi_3}$}
\begin{array}{llllll}
 \! \! \! \! \! \!  \forall\, R>0\,:    & u_n\to
u,  &  \cg {u_n}\leq R,  &  v_n\weakto v \text{ in }\V
  &\Rightarrow\quad
  \liminf_{n\to\infty}\Psi_{u_n}(v_n)\ge \Psi_u(v)
\\
  \! \! \! \! \! \! \forall\, R>0\,:    & u_n\to u, & \cg {u_n}\leq R, & v\in \V &\Rightarrow\quad
 \left\{
 \begin{array}{ll}
    \exists\, v_n\to v:\\
  \lim_{n\to\infty}\Psi_{u_n}(v_n)=\Psi_u(v).
  \end{array}
  \right.
  \end{array}
\end{equation}
 For later use, we explicitly remark that
assumption~\eqref{eq:41.1} means that
\begin{equation}
\label{superlinear}
\begin{aligned}
& \! \! \! \! \! \!  \forall\, R>0, \,  M>0\,
\\ & \! \! \! \! \! \!  \left\{
\begin{array}{llllll}
\exists K>0  &  \forall\, u\in \domainenergy \text{ with } \cg u\leq
R & \forall\, v \in \V\, :  &  \|v \|\geq K  & \Rightarrow &
\Psi_u(v) \geq M\|v\|,
\\
 \exists K^*>0  &  \forall\, u\in \domainenergy \text{ with }
\cg u\leq R &  \forall\, \xi \in \V\, :
 &  \|\xi \|_* \geq K^* &  \Rightarrow & \Psi_u^*(\xi) \geq
M\|\xi\|_*.
\end{array}
\right.
\end{aligned}
\end{equation}
We  also recall an important consequence of assumption \eqref{mosco}
(see \cite[Chap.\,3]{Attouch84}): for all $R>0$
\begin{equation}
\label{e:impo-consequence}
 u_n\to u \ \ \text{ in $V$}, \quad  \cg {u_n}\leq R,  \qquad \xi_n\weakto \xi \ \ \text{ in
 $\V^*$} \ \ \Rightarrow \ \  \liminf_{n\to\infty}\Psi_{u_n}^*(\xi_n)\ge
 \Psi_u^*(\xi).
\end{equation}
Indeed, it has been proved in \cite[Lemma\,4.1]{Stef08?BEPD} that
the first condition in  \eqref{mosco}, combined with
\eqref{e:impo-consequence}, is
 in fact
\emph{equivalent} to \medskip \eqref{mosco}.

\paragraph{\textbf{Assumptions on the  energy functional.}}
We now formulate our assumptions on the functional $\cE$.
 We recall the basic condition
\begin{equation}
\tag{$4.{\mathrm{E}_0}$} \label{Ezero-bis}
\begin{gathered}
u\mapsto \ene t{u} \  \text{is l.s.c. for all $t \in [0,T]$,} \quad
\exists\, C_0>0\, \ \forall\, (t,u) \in [0,T] \times \domainenergy
\, : \ \ene tu \geq C_0 \ \text{ and}
\\
 \graph(\diffname) \text{ is a Borel set
of $[0,T]\times \V \times \V^*$.}
\end{gathered}
\end{equation}
\begin{description}
\item[Coercivity]
For all $t \in [0,T]$
\begin{equation}
  \label{eq:17-bis}
  \tag{$4.{\mathrm{E}_{1}}$}
 \text{the map}\quad
  u \mapsto \ene tu  \quad\text{has compact sublevels.}
\end{equation}
\item[Variational sum rule]
  If for some $u_o\in \V$ and $\tau>0$ the point  $\bar u$ is a minimizer of
  $u\mapsto \ene tu+\tau\Psi_{u_o}((u-u_o)/\tau)$, then  $\bar{u}$
  fulfills the Euler-Lagrange equation
  \begin{equation}
    \label{eq:42-bis}
    \tag{$4.{\mathrm{E}_2}$}
    \exists \,  \xi\in \diff t{\bar u} \,:\quad -\xi\in \partial \Psi_{u_o}((\bar u-u_o)/\tau).
  \end{equation}
\item[Lipschitz  continuity]
\begin{equation} \label{eq:diffclass_a}
\tag{$4.{\mathrm{E}_3}$}
 \exists\, C_1>0 \ \forall\, u \in \domainenergy \
\forall\,  t,s \in [0,T]\,: \quad
  |\ene{t}{u} -\ene{s}{u}| \leq
  C_1 \ene tu |t-s|.
\end{equation}
\item[Conditioned one-sided time-differentiability]
\begin{equation}
\label{e:ass_p_a} \tag{$4.{\mathrm{E}_4}$}
\begin{gathered}
\text{there exists a Borel function} \quad \Ptname:
\mathrm{graph}(\diffname) \to \R \quad \text{and a constant $C_2>0$ such that} \\
\forall\, (t,u,\xi) \in \mathrm{graph}(\diffname)\, : \quad
\liminf_{h \downarrow 0} \frac{\ene{t+h}{u}- \ene{t}{u}}{h} \leq
\Pt{t}{u}{\xi} \leq C_2\cg{u}\,.
\end{gathered}
\end{equation}
\item[Chain-rule inequality]
 $\cE$ satisfies the   \emph{chain-rule inequality with respect
to the triple $(\Psi,\diffname,\Ptname)$}, i.e.\
 for every absolutely continuous curve
   $u \in \AC ([0,T];\V)$
  and for all
  $\xi\in L^1(0,T;\V^*)$
 such that
 \begin{equation}
\label{e:basic-conditions} \sup_{t \in (0,T)} |\ene t
{u(t)}|<+\infty, \qquad  \xi(t)\in \diff t{u(t)}\ \ \text{for a.a.\
}t\in (0,T),
 \end{equation}
    \begin{equation}
  \label{conditions-1}
  \begin{gathered}
  \int_0^T\Psi_{u(t)}(u'(t))\,\d t<+\infty,\qquad \text{and}\qquad
    \int_0^T\Psi_{u(t)}^*(-\xi(t))\,\d t<+\infty,
\end{gathered}
 \end{equation}
\begin{equation}
    \label{eq:48strong}
     \tag{$4.{\mathrm{E}_5}$}
    \begin{gathered}
     \text{the map  $t\mapsto \ene t{u(t)}$ is absolutely continuous
     and}
     \\
    \frac \d{\d t}\ene t{u(t)} \geq \la \xi(t),u'(t)\ra+
    \Pt{t}{u(t)}{\xi(t)}
    \quad
    \text{for a.a.\ }t\in (0,T).
    \end{gathered}
  \end{equation}
\item[Weak closedness of $(\cE,\diffname,\Ptname)$]
For all $t \in [0,T]$ and for all sequences $\{ u_n\} \subset \V$,
     $\, \xi_n\in \diff{t}{u_n},$
     $\,\EE_n=\ene{t}{u_n}, $  $\, p_n
     =\Pt{t}{u_n}{\xi_n}$
    fulfilling
\[
     u_n\to u \ \ \text{in $\V$,}  \quad  \xi_n\weakto \xi  \ \ \text{weakly in
     $\V^*$,}
     \quad
    p_n \to p   \ \ \text{and}  \ \  \EE_n\to \EE \ \ \text{in
    $\R$,}
\]
there holds
    \begin{equation}
      \label{eq:468}
      \tag{$4.{\mathrm{E}_6}$}
      (t,u)\in \dom(F),\quad
      \xi\in \diff tu,\quad
      p\leq \Pt{t}{u}{\xi},\quad
      \EE=\ene tu.
    \end{equation}
\end{description}
For later use, we point out that \eqref{eq:diffclass_a} yields the
following estimate
\begin{equation}
\label{e:useful-later}
\begin{aligned}
 \exists\, C_3>0 \ \
  \forall\, u \in \domainenergy\, : \cg{u} \leq C_3 \inf_{t
\in [0,T]} \ene{t}{u}\,.
\end{aligned}
\end{equation}

Notice that, under the above conditions (cf.\
\eqref{eq:diffclass_a}), for fixed $u \in \domainenergy$ the
function $t \mapsto \ene tu$ is Lipschitz continuous, hence a.e.\
differentiable. Still, it may happen that, along some  curve $u \in
\AC ([0,T];V)$, the energy $ \ene t{u}$ is not differentiable at
$(t,u(t))$, for any   $t \in [0,T]$, cf.\ e.g.\ Example
\ref{ex:3.1}. Hence, one needs to recur to the generalized notion
$\Ptname$.
\begin{example}[Example~\ref{ex:3.1} revisited]
\upshape \label{rmk:example-p} Let us refer to the setting of
Example~\ref{ex:3.1}, and to the subdifferential $\diff tu  =
\partial_u^{\mathrm{Clarke}}\ene tu$, explicitly calculated in~\eqref{clarke} (analogous considerations can
be developed in the case $\diff tu  = \margsub \ene tu$ examined in
Example~\ref{ex:3.2}). Since $\frsub \ene tu \subset
\partial_u^{\mathrm{Clarke}}\ene tu$, in view of the forthcoming
Proposition~\ref{l:variational-sum-rule},  condition
\eqref{eq:42-bis} is satisfied. As for the choice of the function
$\Ptname: \mathrm{graph}(\diffname) \to \R $ in such a way that
 chain-rule inequality holds, it
follows from \eqref{clarke} that
\[
\left\{
\begin{array}{llll}
u >\beta t  & \Rightarrow & \diff tu = \{ -\alpha\}, &  \Pt
t{u}{-\alpha}= \alpha \beta,
\\
u <\beta t  & \Rightarrow & \diff tu = \{ \alpha\}, &  \Pt
t{u}{\alpha}= -\alpha \beta.
\end{array} \right.
\]
Asking for the chain-rule inequality \eqref{e:crule}  along the
curve $\bar{u}(t)=\beta t$   only needs, for every $\xi \in \diff
t{\bar{u}(t)}=[-\alpha,\alpha]$, that  $\Pt t{\bar{u}(t)}{\xi} $ is
a selection in the set $[-\alpha \beta,-\xi\beta]$ of the admissible
$p$'s. However, the closedness condition \eqref{eq:468} is fulfilled
only for the choice $\Pt t{\bar{u}(t)}{\xi} = -\xi \beta$.
\end{example}
  We conclude this section with a result providing sufficient
conditions for the variational sum rule \eqref{eq:42-bis}. As in the
case of sum rules for convex functionals, we use that $\Psi_{u_o}$ is
locally Lipschitz, since its domain is the whole space $\V$.
 Our proof relies on
\cite[Lemma 2.32]{Mordukh06}.
\begin{proposition}
\label{l:variational-sum-rule} Let $\{\Psi_u\}$ be a family of
admissible dissipation potentials on the reflexive space $\V$, and
$\cE : [0,T]\times \V \to (-\infty,+\infty]$ an energy functional
complying with \eqref{Ezero-bis}, with subdifferential $\diffname:
[0,T]\times \V \rightrightarrows \V^*$. Suppose that
\begin{equation}
\label{subdifferential-inclu} \frsub \ene tu \subset \diff tu \quad
\text{for every } (t,u)\in [0,T]\times \domainenergy,
 \end{equation}
and that $(\cE,\diffname)$ comply with the weak closedness condition
\eqref{eq:468}.

Then, the variational sum rule \eqref{eq:42-bis} holds.
\end{proposition}
\PROOF
 Let $\bar{u}$ be a minimizer of
  $u\mapsto \ene tu+\tau\Psi_{u_o}((u-u_o)/\tau)$. It follows from~\cite[Lemma 2.32, p.\ 214]{Mordukh06}
that
\begin{equation}
\label{eta-inequalities} \forall\, \eta >0 \ \ \exists\, u^1_\eta
\in V, \ u^2_\eta \in \domainenergy \, :  \quad \left\{\ba{cl}
\Vnorm{u^1_\eta-\bar{u}} + \Vnorm{u^2_\eta-\bar{u}} &\leq \eta\,,
\\
| \Psi_{u_o} \big(\frac{u^1_\eta - u_o}\tau \big)
-\Psi_{u_o}\big(\frac{\bar{u} - u_o}\tau \big)| & \leq \eta\,,
\\
| \ene{t}{u^2_\eta} - \ene{t}{\bar{u}} | &\leq \eta\,, \ea \right.
\end{equation}
 and
\begin{equation}
\label{eta-equality}
 \exists \, w_\eta \in \partial \Psi_{u_o}
\left(\frac{u^1_\eta - u_o}\tau \right), \ \xi_\eta \in \frsub
\cE_{t}(u^2_\eta), \
 \zeta_\eta \in \Vanach^*\, : \quad  \Vnorm{\zeta_\eta}_* \leq \eta,
\  \ \text{and} \ \ w_\eta + \xi_\eta +  \zeta_\eta=0\,.
\end{equation}
Due to \eqref{subdifferential-inclu}, we have $\xi_\eta \in \diff
t{u^2_\eta}$.  Choosing $\eta=1/n$, we find sequences $(u^1_n)$,
$(u^2_n)$, $(w_n)$, $(\xi_n)$, and $(\zeta_n)$ such that $\zeta_n
\to 0$ in $\Vanach^*$,
 $u^1_n \to
\bar{u}$ and $u^2_n \to \bar{u}$ in $V$, with $ \cE_{t}(u^2_n) \to
\cE_{t}(\bar{u}) $. Since $w_n \in \partial \Psi_{u_o} \left(\frac{u^1_n -
u_o}\tau \right)$  and   $\partial \Psi_{u_o} : \V \rightrightarrows \V^*$ is
a bounded operator, we also deduce that $\sup_n \Vnorm{w_n}_*
<+\infty$. Hence, in view of \eqref{eta-equality},
 we ultimately have that $(\xi_n)$ is bounded in $\Vanach^*$.
Thus,
 there exists $\xi \in \Vanach^*$ such that, up to a (not relabeled)
 subsequence, $\xi_n \weakto \xi$ in $\V^*$. Due to \eqref{eq:468},
 we conclude  that $\xi \in \diff t {\bar{u}}$. On the other hand,
 passing to the limit in \eqref{eta-equality} and using the
 well-known
strong-weak
 closedness property of $\partial \Psi_{u_o}$ gives  $-\xi \in \partial
 \Psi_{u_o}((\bar{u}-u_o)/\tau)$, and \eqref{eq:42-bis} ensues.
\QED
\subsection{Approximation}
\label{ss:3.1}
 For a fixed  initial datum $u_0 \in \domainenergy$ and a time step $\tau>0$,
 we consider a uniform
partition $\mathscr{P}_\tau= \{0=t_0<t_1<t_2<\cdots<t_{N-1}<T\le
t_{N}\}$ with $t_n:=n\tau$, $U^0_\tau =  u_0$,  and construct a
sequence $(U^n_\tau)_{n=1}^{N}$ by recursively solving
\begin{equation}
  \label{e:min-scheme}
  U^n_\tau\in \argmin_{U \in \domainenergy} \left\{\tau\Psi_{U^{n-1}_\tau}\left(\frac{U-U^{n-1}_\tau}\tau\right)+\ene{t_n}U
\right\}
  \quad n=1,\cdots, N.
\end{equation}
Using the \emph{direct method} in the Calculus of Variations and
exploiting assumption~\eqref{eq:17-bis}, one sees (cf.\
Lemma~\ref{lemma:degiorgi}) that
 for all $u_0 \in \domainenergy $
and $\tau \in (0,\tau_o)$ there exists at least one solution
$(U^n_\tau)_{n=1}^{N}$ to the time-incremental minimization
problem~\eqref{e:min-scheme}. We  denote by $\pwC {U}\tau$ and
$\upwC {U}{\:\tau}$, respectively,  the left-continuous and
right-continuous piecewise constant interpolants of the values
$(U^n_\tau)_{n=1}^{N}$, i.e.,
\begin{equation}
\label{e:pwc} \pwC {U}\tau (t):=\Utau^n \quad \text{for } t \in
(t_{n-1},t_n], \quad \quad \upwC {U}{\:\tau}(t):=\Utau^{n-1} \quad
\text{for } t \in [t_{n-1},t_n), \quad n=1,\ldots, N,
\end{equation}
and by $\pwL {U}\tau$  the piecewise linear interpolant
\begin{equation}
\label{e:pwl} \pwL {U}\tau (t):= \frac{t-t_{n-1}}{\tau}\Utau^n +
\frac{t_n-t}{\tau}\Utau^{n-1}  \quad \text{for $ t \in
[t_{n-1},t_n),$} \quad n=1,\ldots, N.
\end{equation}
Thanks to \eqref{eq:42-bis}, for all $n =1,\ldots,N$ there exists
$\xi_\tau^{n} \in \diff{t_n}{\Utau^n} \cap
(-\partial\Psi_{U^{n-1}_\tau} ({\Utau^n - \Utau^{n-1}}/\tau))$. We
denote
 by $\pwC \xi{\tau}$ the (left-continuous) piecewise
constant interpolant of the family $(\xi_\tau^{n})_{n=1}^{N} \subset
\V^*$.

 Furthermore, we  also consider the \emph{variational
interpolant} $ \pwM U\tau$ of the discrete
va\-lues~$(\Utau^n)_{n=1}^{N}$, which was first introduced  by
\textsc{E. De Giorgi} within the \emph{Minimizing Movements} theory
(see~\cite{DeGiorgi-Marino-Tosques80, DeGiorgi93, Ambrosio95}). It
is defined in the following way: the map $t\mapsto \pwM U\tau (t)$
is Lebesgue measurable in $(0,T)$ and  satisfies
\begin{equation}
  \label{interpmin}
  \left\{
    \begin{array}{l}
      \pwM U\tau(0)=u_0, \quad \text{and, for }
      t=t_{n-1} + r \in (t_{n-1}, t_{n}],
      \\ \\
     \pwM U\tau(t)
      \in \argmin_{U \in \domainenergy}\left\{r
      \Psi_{U^{n-1}_\tau}\left(\frac{U-\Utau^{n-1}}{r}\right) +
 \ene{t}{U} \right\},
    \end{array}
  \right.
\end{equation}
The existence of such a measurable selection is ensured
by~\cite[Cor.\,III.3, Thm.\,III.6]{Castaing-Valadier77}, see
also~\cite[Rem.\,3.4]{Rossi-Mielke-Savare08}. When $ t=t_n$, the
minimization problems \eqref{e:min-scheme} and~\eqref{interpmin}
coincide, so that we may assume
\begin{equation}
\label{e:coincidence-nodes}
  \pwC U\tau(t_n)=    \upwC U{\:\tau}(t_n)=
  \pwL U\tau(t_n)= \pwM U\tau(t_n),\quad
  \text{for every } n =1, \ldots,N.
\end{equation}
Then, $\pwM U\tau$ contains all the information on $\pwL U\tau$,
$\pwC U\tau$, and $\upwC U{\:\tau}$.
 Furthermore, again by ~\eqref{eq:42-bis}  and
 the
measurable selection result \cite[Thm.\,8.2.9]{AubFra90SVA},
 with $\pwM \UU\tau $ we can associate a measurable
function $\pwM {\xi}{\tau}:(0,T) \to \V^*$ fulfilling the Euler
equation for the minimization problem \eqref{interpmin}, i.e.
\begin{equation}
\label{interpxi} \pwM {\xi}{\tau}(t) \in \diff{t}{\pwM \UU\tau(t)}
\cap \left( {-}\partial\Psi_{\upwC \UU{\:\tau}(t)}\left(\frac{\pwM
\UU\tau(t) - \upwC \UU{\:\tau}(t)}{t-t_{n-1}} \right) \right) \qquad
\forall\, t \in [t_{n-1}, t_n), \ \ n =1,\ldots, N.
\end{equation} 
For later notational convenience, we also introduce the piecewise
constants interpolants $\pwC {\mathsf{t}}\tau$ and $\upwC
{\mathsf{t}}{\:\tau}$ associated with the partition
$\mathscr{P}_{\tau}$, namely
\begin{equation}
\label{interp-time}
 \pwC {\mathsf{t}}{\tau}(0)=\upwC {\mathsf{t}}{\:\tau}(0):=0, \quad \pwC
{\mathsf{t}}{\tau}(t):=t_k \quad \text{for} \quad  t \in (
t_{k-1},t_{k}],
 \quad \upwC {\mathsf{t}}{\:\tau}(t):=t_{k-1} \quad
\text{for} \quad  t \in [ t_{k-1},t_{k}).
\end{equation}
 Of course, for every $t
\in [0,T]$ we have $\pwC {\mathsf{t}}{\tau} (t) \down t$  and $\upwC
{\mathsf{t}}{\:\tau} (t) \up t$ as $\tau \down 0$.

\subsection{Main existence result} \label{sec:viscous}
Before stating our main existence result, let us first specify the
notion of solution we are interested in.
\begin{definition}
\label{def-energy-pair} Let $\cE : [0,T]\times \V \to
(-\infty,+\infty]$ be an energy functional fulfilling
\eqref{Ezero-bis} and  $\{ \Psi_u\}_{u \in \domainenergy}$ a family
of \emph{admissible dissipation potentials}. Suppose that $\cE$
complies with the chain rule \eqref{eq:48strong}. We say that
$(u,\xi) \in \AC ([0,T];\V) \times L^1 (0,T;\V^*) $ is a
\emph{solution pair} to the generalized gradient system
$(\V,\cE,\Psi,\diffname,\Ptname)$ if
\begin{enumerate}
\item
$(u,\xi) $ fulfills the doubly nonlinear equation
\begin{equation}
\label{xi-equation} \xi(t) \in \diff t{u(t)}, \quad \partial
\Psi_{u(t)}(u'(t)) + \xi(t)\ni 0 \quad \foraa\, t \in (0,T),
\end{equation}
\item  $(u,\xi) $ complies with the energy identity
\begin{equation}
\label{enid1}
 \int_{s}^{t} \big( \Psi_{u(r)}(u'(r))
 {+}\Psi^*_{u(r)}(-\xi(r)) \big)
 \,\d r +   \ene {t}{u(t)} =
\ene {s}{u(s)} + \int_s^t \Pt{r}{u(r)}{\xi(r)}\dd r
\end{equation}
for every $0 \leq s \leq t \leq T$.
\end{enumerate}
We shortly say that $u \in \AC ([0,T];\V)$ is a solution to the
generalized gradient system $(\V,\cE,\Psi,\diffname,\Ptname)$, if
there exists $\xi \in L^1 (0,T;\V^*) $ such that $(u,\xi)$ is  a
solution pair to $(\V,\cE,\diffname,\Ptname,\Psi)$.
\end{definition}
\begin{theorem}[Existence]
  \label{thm:viscous2}
  Assume that
$(\V,\cE,\Psi,\diffname,\Ptname)$ comply
with~\eqref{e:4.1}--\eqref{mosco}
 and with~\eqref{Ezero-bis}--\eqref{eq:468}.

Then, for every $u_0\in \domainenergy$ there exists a solution $u\in
\AC([0,T];\V)$ to the doubly nonlinear equation~\eqref{eq:1-bis},
fulfilling the initial condition $u(0)=u_0$.

 In fact, for any family of approximate solutions $(\pwM{U}{\tau},\pwM{\xi}{\tau})_{\tau>0}$
  there exist
 a sequence
$\tau_k \down 0$ as $k \to \infty$,   and $\xi \in L^1 (0,T;\V^*)$
such that the following convergences
  hold as~$k \to \infty$
\begin{subequations}
\label{e:conv}
\begin{align}
& \label{convu1} \pwC{U}{\tau_k},\ \pwL{U}{\tau_k}, \
\pwM{U}{\tau_k} \to u\quad \text{in $L^\infty (0,T;\V)$,}
\\
 & \label{convu2} \pwL{U}{\tau_k} \weakto u \quad \text{in
$W^{1,1}(0,T;\V)$,}
\\
& \label{convu4} \ene{t}{\pwC U{\tau_k}(t)} \to \ene{t}{u(t)}\quad
\text{for all } t \in [0,T],
\\
& \label{convu5}  \int_{s}^{t} \Psi_{\upwC U{\:\tau_k}(r)}(\pwL
{U'}{\tau_k}(r))\,\d r \to \int_{s}^{t} \Psi_{u(r)}(u'(r))\,\d r
\quad \text{for all }   0 \leq s \leq t \leq T,
\\
&
 \int_{s}^{t} \Psi^*_{\upwC
U{\:\tau_k}(r)}(-\pwM \xi{\tau_k}(r))\,\d r \to \int_{s}^{t}
\Psi^*_{u(r)}(-\xi(r))\,\d r \quad  \text{for all }   0 \leq s \leq t \leq
T, \label{enhanced-middle}
\end{align}
\end{subequations}
and $(u,\xi) $ is   a solution pair to the generalized gradient
system $(\V,\cE,\Psi,\diffname,\Ptname)$.

Furthermore, if
\begin{equation}
\label{strict-convexity} \text{
    $\Psi_u^*$ is  strictly convex for all $u \in \V$,}
\end{equation}
 we have the additional convergence
\begin{align}
 \label{convu3}
 &
 \pwM{\xi}{\tau_k} \weakto \xi \quad \text{in
$L^{1}(0,T;\V^*)$.}
\end{align}
\end{theorem}
\noindent The proof of Theorem~\ref{thm:viscous2} is developed
throughout Section~\ref{s:5}.
\begin{remark}
\label{rmk:versus} The considerations set forth
 in Remark
\ref{rem:towards-ineq} for energies smoothly depending on time
extend to the present setting. Namely,  the proof of Theorem
\ref{thm:viscous2} reveals that the \emph{one-sided} chain-rule
inequality \eqref{eq:48strong} is  sufficient to conclude  the
existence of solutions to the Cauchy problem for \eqref{eq:1-bis},
in that it is combined with the upper energy estimate following from
the discretization scheme.

 Clearly, in order to enforce  the energy identity \eqref{enid1} for   \emph{any} solution to
 \eqref{eq:1-bis}, it would be necessary to impose \eqref{eq:48strong}
 as an \emph{equality}. As shown by Example \ref{ex:3.1}, this may
 lead to restrictions on the admissible functions $\Ptname$.
\end{remark}
\paragraph{\textbf{Weakened assumptions.}}
The two ensuing remarks  explore the possibility of refining our
requirements on the chain rule~\eqref{eq:48strong}, and on the
properties of the dissipation potentials.
\begin{remark}[A weaker chain rule]
\upshape \label{rem:weaker}  Like in the gradient flow case
 (cf.\ \cite[Thm.\,2]{Rossi-Savare06},
 and~\cite[Thm.\,2.3.1]{Ambrosio-Gigli-Savare08} in the metric setting),
  it is possible to state our  existence result
   for the Cauchy problem for~\eqref{eq:1-bis} under  a (slightly)
weaker form of the chain rule~\eqref{eq:48strong}, which requires
that for every absolutely continuous curve $u \in \AC([0,T];\V)$ and
  $\xi\in L^1(0,T;\V^*)$ satisfying~\eqref{e:basic-conditions},
 as well as
  \begin{equation}
  \label{conditions}
  \begin{gathered}
 \int_0^T\Psi_{u(t)}(u'(t))\,\d t<+\infty,  \  \ \  \     \int_0^T\Psi_{u(t)}^*(-\xi(t))\,\d t<+\infty,
    \\
    \text{and such  that  the map $
t\mapsto \ene t{u(t)} $ is a.e.\ equal to a function $\EE$ of
bounded variation,}
\end{gathered}
 \end{equation}
   there holds
\begin{equation}
    \label{eq:48}
    \tag{$\textsc{chain}_\mathrm{weak}$}
    \frac \d{\d t} \EE(t) \geq \la \xi(t),u'(t)\ra+
    \Pt{t}{u(t)}{\xi(t)}
    \quad
    \text{for a.a.\ }t\in (0,T).
  \end{equation}
In this case, suitably adapting the proof
of~\cite[Thm.\,2]{Rossi-Savare06}  to the doubly nonlinear case, one
obtains that there exist $u \in \AC([0,T];\V)$ and $\xi \in L^1
(0,T;\V^*)$ fulfilling the differential
inclusion~\eqref{xi-equation} and the following \emph{energy
inequality} (compare with the energy identity~\eqref{enid1} under
the chain rule~\eqref{eq:48strong})
\begin{equation}
\label{weaker-enid}
\begin{aligned}
   \int_{s}^{t} \big( \Psi_{u(r)}(u'(r))
 {+}\Psi^*_{u(r)}(-\xi(r)) \big)
 \,\d r &  +   \ene {t}{u(t)} \leq
\ene {s}{u(s)} + \int_s^t \Pt{r}{u(r)}{\xi(r)}\dd r \\ &  \text{for
all $t \in [0,T]$ and almost all $s \in (0,t)$.}
\end{aligned}
\end{equation}
\end{remark}
\begin{remark}[Weaker conditions on the dissipation]
In fact, condition~\eqref{eq:psi-sum-2} in the definition of the
dissipation potentials $(\Psi_u)_{u \in \V}$, is only used in the
proof of the forthcoming Lemma~\ref{lemma:degiorgi}, which is the
crucial technical result for  the a priori estimates on the
approximate solutions $(\pwC \UU\tau)$, $(\pwL \UU\tau)$ and $(\pwM
\UU\tau)$. As shown in Remark~\ref{rem:alternative} later on, it is
possible to dispense with~\eqref{eq:psi-sum-2} if the following
condition, involving the energy $\cE: [0,T]\times \V \to
(-\infty,+\infty]$ and its subdifferential mapping $\diffname:
[0,T]\times \V \rightrightarrows \V^*$, holds:
\begin{equation}
\label{structure-condition}
\begin{gathered}
\text{for all $u,\,v \in \V$ the directional derivative }
\dire\ene{t}{u;v} := \lim_{h \down 0} \frac1h\left(
\ene{t}{u+hv}{-}\ene{t}{u}\right) \text{ exists,}
\\
\text{and }\langle \xi,v \rangle \geq \dire\ene{t}{u;v} \quad
\text{for all $\xi \in \diff{t}{u}$ and $v \in \V$.}
\end{gathered}
\end{equation}
Condition~\eqref{structure-condition} has to be coupled with a
strengthened version of the first inequality in~\eqref{e:ass_p_a},
namely
\begin{equation}
\label{strengthened}
\begin{aligned}
 &  \text{for every } (t,u,\xi) \in
\mathrm{graph}(\diffname) \text{ and }
 (u_h)\subset\V  \text{ such
that }  u_h \to u \text{ as $h \down 0$,} \\ &  \text{ there holds }
 \liminf_{h \downarrow 0} \frac{\ene{t+h}{u_h}- \ene{t}{u_h}}{h} \leq
\Pt{t}{u}{\xi}.
\end{aligned}
\end{equation}
Notice that \eqref{structure-condition}   holds for marginal
functionals which are $\lambda$-concave.
\end{remark}
\subsection{Upper semicontinuity of the set of  solutions.}
\label{ss:stability}
We now address the issue of \emph{upper semicontinuity} of the set
of solutions to  the Cauchy problem for \eqref{eq:1-bis},
 with respect to
convergence of the initial data and (a suitable kind of) variational
convergence for the driving energy functionals.

We consider sequences $(\V,\cE^n,\Psi^n,\diffname^n,\Ptname^n)$ of
generalized gradient systems, and impose the following.
 \paragraph{\textbf{Assumption
(H1).}}
 Let $(\cE^n)_{n \in \N}$ be a sequence of lower semicontinuous
 energy
functionals $\cE^n : [0,T]\times \V \to (-\infty,+\infty]$, with
domains $\dom(\cE^n) = [0,T]\times \domainenergy_n$ for some 
$\domainenergy_n \subset \V$,   and with subdifferentials
$\diffname^n: [0,T]\times \domainenergy_n \rightrightarrows \V^*$;
 we use the notation
$\cgn{u}:= \sup_{t \in [0,T]}\enei ntu$ for  $u \in
\domainenergy_n$.  We suppose that the functionals $(\cE^n)_{n \in
\N}$ comply with \eqref{Ezero-bis}, \eqref{eq:17-bis},
\eqref{eq:diffclass_a}, and \eqref{e:ass_p_a}, with constants
\emph{uniform} with respect to $n \in \N$.
We also  require  that there exists  a generalized gradient system
$(\V, \cE,\Psi,\diffname,\Ptname)$, such that the energy
\[
 \cE: [0,T]\times \V \to (-\infty,+\infty] \ \
\text{complies with \eqref{Ezero-bis} and the chain rule
\eqref{eq:48strong},}
\]
and
  the functionals $(\cE^n)_n$ converge to
$\cE $ in the following sense: for all $t \in [0,T]$ and for all
sequences $\{ u_n\} \subset \V$,
     $\, \xi_n\in \diffi n{t}{u_n} $,
    fulfilling
\[
     u_n\to u \ \ \text{in $\V$,}  \quad  \xi_n\weakto \xi  \ \ \text{weakly in
     $\V^*$,}
     \quad
    \Pti n{t}{u_n}{\xi_n} \to p,
\]
 there holds
    \begin{equation}
    \label{closure-n}
    \begin{gathered}
      (t,u)\in \dom(F),\quad
      \xi\in \diff tu,\quad
      p\leq \Pt{t}{u}{\xi},
\\
 \text{and, if $\enei n{t}{u_n}$ converges to some $\EE \in \R$, then $\EE
= \ene tu$.}
\end{gathered}
\end{equation}
 \paragraph{\textbf{Assumption (H2).}} Let  $\{\Psi_u^n
\}_{u \in \domainenergy_n}$  be a family of \emph{admissible}
dissipation potentials, satisfying conditions of superlinear growth
on sublevels of the energies $\cE^n$,  uniformly with respect to $n$
(i.e., \eqref{superlinear} holds for constants \emph{independent} of
$n$). We also suppose that the potentials $(\Psi_u^n)_{u \in
\domainenergy_n}$ Mosco converge on sublevels of the energies to a
family $(\Psi_u)_{u \in \domainenergy}$ of admissible potentials,
viz.
\begin{equation}
\label{ppsi-n}
\begin{array}{llllll}
 \! \! \! \! \! \!  \forall\, R>0\,:    & u_n\to
u,  &  \sup_{n \in \N}\cgn {u_n} \leq R,  &  v_n\weakto v \text{ in
}\V
  &\Rightarrow\quad
  \liminf_{n\to\infty}\Psi_{u_n}^n(v_n)\ge \Psi_u(v)
\\
  \! \! \! \! \! \! \forall\, R>0\,:    & u_n\to u, & \sup_{n \in \N}\cgn {u_n}\leq R, & v\in \V &\Rightarrow\quad
 \left\{
 \begin{array}{ll}
    \exists\, v_n\to v:\\
  \lim_{n\to\infty}\Psi_{u_n}^n(v_n)=\Psi_u(v).
  \end{array}
  \right.
  \end{array}
\end{equation}
\begin{theorem}[Upper semicontinuity]
  \label{thm:viscous-stability}
 Let $(\V, \cE_n,\Psi_n,\diffname_n,\Ptname_n)$  be a family of
  generalized gradient systems
 complying
with \textbf{Assumption (H1)} and  \textbf{Assumption (H2)}.  Let $(
u_0^n )_n$ be a sequence of initial data, with $u_0^n \in
\domainenergy_n$ for all $n\in \N$, such that
\begin{equation}
\label{data-convergence} u_0^n \weakto u_0   \ \  \text{ in $\V$ and }  \ \
 \enei n{0}{u_0^n} \to \ene {0}{u_0},
\end{equation}
and let $(u_n,\xi_n)_{n \in \N}$ a sequence of \emph{solution pairs}
to the Cauchy problems
\begin{equation}
\label{cauchy-n}
\partial\Psi_{u(t)}(u'(t))+\diffi n{t}{u(t)}\ni0\quad\text{in }\V^*,\quad \text{for a.a.\ }t\in (0,T); \qquad
  u(0)=u_0^n
\end{equation}
(in particular, complying with the energy identity \eqref{enid1} for
all $n \in \N$). Then, there exist a subsequence
$(u_{n_k},\xi_{n_k})_{k \in \N}$ and functions  $(u,\xi) \in
\AC([0,T];\V)\times L^1(0,T;\V^*)$ such that $(u,\xi)$ is a solution
pair of the Cauchy problem for \eqref{eq:1-bis}, and the following
convergences hold as $k \to \infty$
\begin{subequations}
\label{e:conv-enne}
\begin{align}
& \label{convu1-n} u_{n_k} \to u  \text{ in $\mathrm{C}^0
([0,T];\V)$,} \quad u_{n_k} \weakto u \text{ in $W^{1,1}(0,T;\V)$,}
\\
& \label{convu2-n}
  \enei {n_{k}}{t}{u_{n_k}(t)} \to \ene{t}{u(t)} \text{ for all $t
\in [0,T],$}
\\
& \label{convu5-n} \left\{
\begin{array}{ll}
 \int_{s}^{t} \Psi_{u_{n_k}(r)}^{n_k}(u_{n_k}'(r))\,\d
r \to \int_{s}^{t} \Psi_{u(r)}(u'(r))\,\d r, \\ \int_{s}^{t}
(\Psi^{n_k}_{u_{n_k}(r)})^*(-\xi_{n_k}(r))\,\d r \to \int_{s}^{t}
\Psi^*_{u(r)}(-\xi(r))\,\d r  \end{array}\right. \
 \text{ for all $
0 \leq s \leq t \leq T$.}
\end{align}
\end{subequations}
\end{theorem}
\noindent The proof of this result is outlined at the end of Section
\ref{s:5}.
\begin{remark}
\upshape Suppose that  the energy functionals $\cE^n$ have the
special form
\[
\begin{aligned}
\enei ntu= E^n(u) - \langle \ell^n(t), u \rangle, \ \  & \text{
 with $E^n: \V \to (-\infty,+\infty]$
convex functionals and}
\\
&
 (\ell^n)  \subset \mathrm{C}^1 ([0,T];\V^*).
\end{aligned}
\]
Hence,  if  the functionals $(E^n)$ Mosco-converge to some convex
functional $E: \V \to (-\infty,+\infty]$, and if the functions
$(\ell_n)$ suitably converge to some $\ell\in \mathrm{C}^1
([0,T];\V^*) $,
 then the energies $(\cE^n)$
converge to $\ene tu := E(u) - \langle \ell(t), u \rangle$ in the
sense specified by \textbf{Assumption (H1)}. Indeed, Theorem
\ref{thm:viscous-stability} might be viewed as an extension, to the
doubly nonlinear case, of the
 result on stability of gradient flows
 (with $V$ a Hilbert space and $\Psi(u)=\frac 12 \|u\|_V^2$),
 with respect to Mosco-convergence of
 the (convex) energies, stated in  \cite[Thm.\,3.74(2), p.\,388]{Attouch84}.
 The reader may also consult \cite{Stef08?BEPD}
and the references therein.
\end{remark}
\section{Application: evolutions driven by marginal functionals in finite-strain elasticity}
\label{s:4}
  In this section we
 examine a mechanical model for finite-strain elasticity, described
in terms of the elastic deformation and of some internal,
dissipative variable $z$. Its analysis  has already been developed
in~\cite{FraMie06ERCR}, in the case of a \emph{rate-independent}
evolution for $z$. Therein, existence of \emph{energetic solutions}
to the (Cauchy problem for the) related PDE system has been proved.
Here, we address the case in which the evolution of  $z$ is driven
by \emph{viscous} dissipation.
\subsection{Problem set-up and existence result}
We consider an elastic body occupying a bounded domain $\Omega
\subset \R^d$, $d\geq 1$, with Lipschitz boundary $\Gamma$. We
denote by $\phi: \Omega \to \R^d$ the elastic deformation field, and
assume that the inelasticity of $\Omega$ is described by an internal
variable $z: \Omega \to \R^m$, $m \geq 1$, which we may envisage as
a mesoscopic averaged phase variable.
\paragraph{\textbf{Energy functional.}}
The stored energy  ${\marginale} =\marginal t z\phi$ has the form
\begin{equation}
\label{stored-energy} \marginal t z {\phi}= \cE^1(z)  + \marginali 2
t z{\phi}+E_0,
\end{equation}
with $E_0\in \R$ to be precised later on (cf. Lemma \ref{l:1}).

 In \eqref{stored-energy}, $\cE^1 : \V \to (-\infty,+\infty]$
is the \emph{convex} functional
\begin{equation}
\label{energy1}
\begin{gathered}
\cE^{1}(z) : = \left\{ \begin{array}{lll}
 \int_{\Omega}\left(\frac1q |\nabla z|^q + I_{\mathrm{K}}(z) \right) \, \d x &
 \text{if $z \in W^{1,q}(\Omega;\R^m),$}
 \\
 +\infty & \text{otherwise,}
 \end{array}
 \right.
 \\
\quad \text{where } q>d,  \text{and  $\mathrm{K} $ is a compact
subset of $\R^m$,}
\end{gathered}
\end{equation}
and
 we consider the Fr\"obenius norm $|\nabla z| = \left( \sum_{j=1}^d \sum_{i=1}^m |\partial_{x_j} z^i|^2\right)^{1/2}$
 of the matrix $\nabla z$.

The \emph{nonconvex} contribution $\marginale^2 : [0,T]\times
W^{1,p}(\Omega;\R^d)\times L^2 (\Omega;\R^m)
 \to (-\infty,+\infty]$
is given by
\[
 \marginali 2 t z \phi := \int_{\Omega} W(\nabla \phi(x),z(x)) \, \d x -
\pairing{}{W^{1,p}}{\ell(t)}{\phi},
\]
where  $p>d$, $\pairing{}{W^{1,p}}{\cdot}{\cdot}$ denotes the
duality pairing between $W^{1,p}(\Omega;\R^d)^*$ and
$W^{1,p}(\Omega;\R^d)$, and we suppose that
\begin{equation}
\label{load} \ell \in \mathrm{C}^{1}([0,T];W^{1,p}(\Omega;\R^d)^*).
\end{equation}
 The
 stored energy density  $W: \R^{d \times d}
\times \mathrm{K}' \to (-\infty,+\infty]$ has domain
$\mathrm{dom}(W)= D_W \times \mathrm{K}'$, where $\mathrm{K}' $ is a
compact subset of $\R^m$ containing $\mathrm{K}$.
 We  neglect the
dependence of $W$ on the variable $x$ for the sake of simplicity and
with no loss of generality.  We impose the  following conditions
on~$W$:
\begin{equation}
\label{w2} \tag{$\mathrm{W}_1$} \exists\, \kappa_1,\, \kappa_2>0 \
\forall\, (F,z)\in \R^{d \times d} \times \mathrm{K}' \, : \ \
W(F,z) \geq \kappa_1 |F|^p - \kappa_2 \quad \text{ with $p>d$};
\end{equation}
\begin{equation}
\label{w1} \tag{$\mathrm{W}_2$}
\begin{aligned}
&
 \exists\, \bbW: \R^{\mu_d} \times \mathrm{K}' \to
(-\infty,+\infty] \ \text{such that}
\\ &
\text{(i)} \ \bbW \ \text{is lower semicontinuous,}
\\
& \text{(ii)} \  \forall\, (F,z) \in \R^{d \times d} \times
\mathrm{K}'\, :
 \ \ W(F,z)= \bbW(\bbM (F),z),
\\
& \text{(iii)} \ \forall z\in \mathrm{K}'\,:  \  \ \bbW(\cdot,z):
\R^{\mu_d} \to (-\infty,+\infty] \text{ is convex;}
\end{aligned}
\end{equation}
\begin{equation}
\label{w3} \tag{$\mathrm{W}_3$}
\begin{aligned}
&
 \text{for all $F \in D_W $ the map $W(F,\cdot)$ is continuous and G\^ateau-differentiable  on $\mathrm{K}'$, and
}
\\ &
\text{(i)} \
 \exists\, \kappa_3,\, \kappa_4>0 \
\forall\, (F,z)\in D_W  \times \mathrm{K}' \, : \ \ |\rmD_z W(F,z)|
\leq \kappa_3 (W(F,z) + \kappa_4)^{1/2};
\\
& \begin{aligned} \text{(ii)} \exists\, \kappa_5,\, \kappa_6>0, \
 & \exists\, \alpha \in (0,1] \
\ \forall\,F \in D_W \ \forall\, z_1,\,z_2 \in \mathrm{K}'\, : \\
& |\rmD_z W(F,z_1) - \rmD_z W(F,z_2)| \leq \kappa_5 |z_1
-z_2|^{\alpha} (W(F,z_1) +\kappa_6)^{1/2}.
\end{aligned}
\end{aligned}
\end{equation}
In \eqref{w1}, we have used the notation $\mu_d = \sum_{s=1}^d
\binom{d}{s}^2 $,  and $\bbM : \R^{d\times d}\to \R^{\mu_d}$ is the
function which maps a matrix to all its minors (subdeterminants).
Hence, \eqref{w1} states that for all $z \in \mathrm{K}'$ the map
$W(\cdot,z)$ is \emph{polyconvex}.  
\paragraph{\textbf{Dissipation.}}
We consider a measurable (dissipation density) function $\psi:
\mathrm{K}\times \R^m \to [0,+\infty)$ (again, we  omit the
dependence of $\psi$ on the variable $x$ with no loss of
generality), fulfilling
\begin{equation}
\label{d1} \tag{$\psi_1$}
\begin{aligned}
 \psi:  K \times \R^m \to [0,+\infty) \ \text{ is
continuous;}
\end{aligned}
\end{equation}
\begin{equation}
\label{d2} \tag{$\psi_2$}
\begin{aligned}
  &  \forall\, z \in
\mathrm{K}\,: \  \   \   \psi(z,\cdot):  \R^m \to [0,+\infty) \
\text{ is convex, with } \ \psi(z,0)=0, \ \text{ and}
\\
& \psi^* (z,w_1) = \psi^* (z,w_2) \ \text{ for all $w_1,\,w_2
\in \partial_v \psi(z,v)$ and all $v \in \R^m$;}
\end{aligned}
\end{equation}
\begin{equation}
\label{d4} \tag{$\psi_3$}
\begin{aligned}
  \exists\,
\kappa_7,\, \kappa_8, \, \kappa_9>0 \ \
 \forall\, z \in \mathrm{K} \ \forall\,v,\,\xi \in \R^m \, : \quad
\left\{
\begin{array}{ll}
 \psi(z,v) \geq \kappa_7 |v|^2 -\kappa_9,
 \\
 \psi^*(z,\xi)  \geq   \kappa_8 |\xi|^2 - \kappa_9.
\end{array}
 \right.
\end{aligned}
\end{equation}
In \eqref{d2} the symbols $\partial_v \psi$ and $\psi^*$
 respectively  denote the subdifferential and the Fenchel-Moreau conjugate
of the function $\psi(z,\cdot)$.
Let us point out that there is a crucial interplay between the exponent $1/2$ in \eqref{w3}(ii),
and the exponents $2$ in \eqref{d4}, see also Remark \ref{rem:interplay} later on.
\paragraph{\textbf{PDE system and existence theorem.}} Within this setting, we address the analysis
of the doubly nonlinear evolution equation
\begin{subequations}
\label{pde}
\begin{equation}
\label{pde1}
\begin{aligned}
\partial_v \psi(z(t,x),\dot{z}(t,x)) -\Delta_q z(t,x) + \partial
I_{\mathrm{K}}(z(t,x))  +  & \rmD_z W(\nabla \phi(t,x),z(t,x)) \ni 0
\\ &  \foraa\, (t,x) \in  (0,T)\times \Omega,
\end{aligned}
\end{equation}
(where $\Delta_q z= \mathrm{div}(|\nabla z|^{q-2} \nabla z)$),
supplemented with homogeneous Neumann boundary conditions, and
coupled with the minimum problem
\begin{equation}
\label{pde2} \phi(t,x) \in \argmin\{ \marginal t {z(t,x)} \phi \, :
\ \phi \in \calF\} \quad \foraa\, (t,x) \in
(0,T)\times \Omega,
\end{equation}
\end{subequations}
where $\calF$ denotes the set of   the
kinematically admissible deformation fields, viz.
\[
\mathcal{F}= \left\{ \phi \in
W^{1,p}(\Omega;\R^d)\, : \ \phi = \phi_{\mathrm{Dir}} \text{ on
$\Gamma_{\mathrm{Dir}}$} \right\},
\]
for some $\Gamma_{\mathrm{Dir}} \subset \Gamma$,
 $\Gamma_{\mathrm{Dir}} \neq \emptyset$ with positive Hausdorff measure, and
\begin{equation}
\label{hyp-phi-dir} \phi_{\mathrm{Dir}} \in W^{1,p}
(\Omega;\R^d),  \ \text{such that the map } x \mapsto  \max_{z \in \mathrm{K}} W(\nabla
\phi_{\mathrm{Dir}}(x),z) \text{ is in $L^1 (\Omega)$.}
\end{equation}
\begin{theorem}
\label{th:exist-marginal-elast} Under assumptions \eqref{load},
\eqref{w2}--\eqref{w3}, \eqref{d1}--\eqref{d4}, and
\eqref{hyp-phi-dir}, for every
\begin{equation}
\label{condition-datum} z_0 \in W^{1,q}(\Omega;\R^m)  \ \text{ with
 } \  z_0(x) \in \mathrm{K}   \ \text{ for all $x \in \Omega$,}
\end{equation}
 there exist
functions $z \in L^\infty (0,T;W^{1,q}(\Omega;\R^m)) \cap H^{1}
(0,T;L^2(\Omega;\R^m))$ and $\phi \in L^\infty
(0,T;W^{1,p}(\Omega;\R^d))$ fulfilling \eqref{pde1},
 supplemented with homogeneous Neumann boundary conditions and the initial condition
 $z(0,x) = z_0(x) $ $  \foraa\, x \in \Omega$,
 and \eqref{pde2}. In particular, there exists $\xi \in L^2
(0,T;L^2(\Omega;\R^m)) $ satisfying, for almost all $(t,x) \in (0,T)
\times \Omega$, the inclusions
\begin{equation}
\label{e:pde-xi}
\left\{
\begin{array}{ll}
\partial_v \psi(z(t,x),\dot{z}(t,x))  + \xi(t,x) \ni0, \\
\xi(t,x) \in -\Delta_q z(t,x) + \partial I_{\mathrm{K}}(z(t,x)) +
\rmD_z W(\nabla \phi(t,x),z(t,x))
\end{array}
\right.
\end{equation}
and such that $(z,\phi,\xi)$ fulfill the energy identity for all $0
\leq s \leq t \leq T$
\begin{equation}
\label{enid-marginal}
\begin{aligned}
 \int_{s}^{t}
\int_{\Omega} \big(\psi(z(r,x),\dot{z}(r,x))  & +
\psi^*(z(r,x),{-}\xi(r,x))  \big) \dd x \dd r + \marginal
t{z(t)}{\phi(t)}\\ &  = \marginal s{z(s)}{\phi(s)} - \int_s^t
\pairing{}{W^{1,p}}{\ell'(r)}{\phi(r)} \dd r.
\end{aligned}
\end{equation}
\end{theorem}
\begin{example}
\label{ex:concrete} In finite-strain elasticity there are two main conditions,
namely $(i)$ frame indifference and $(ii)$ local invertibility:
\[
\begin{aligned}
& i)  & W(R F,z)= W(F,z) & \quad \text{ for all $ R \in SO(d), \ F
\in D_W,\ z \in \mathrm{K}'$,}
\\
& ii)  & W(F,z)= \infty & \quad \text{ for all $ \mathrm{det}(F) \leq 0$.}
\end{aligned}
\]
These conditions are compatible with polyconvexity, e.g.\ by
choosing functions of the type
\[
W(F,z)= C |F|^p +w_\mathrm{co}(F,z) + h(\mathrm{det}(F) ),
\]
where $h : \R \to (-\infty,+\infty]$ is continuous, convex, and satisfies
$h(y) =\infty$ for $y \leq 0$. Thus,
 $D_W = \{F\, : \ \mathrm{det}(F)>0\}\subset \R^{d \times d}$
is the nonconvex domain. Recall that the Fr\"obenius norm $|F|=
(\mathrm{tr}(F^{T}F))^{1/2}$ satisfies $|RF|=|F|$. Conditions
\eqref{w1}--\eqref{w3} can be now satisfied if the coupling energy
$w_\mathrm{co}$ satisfies $w_\mathrm{co} \in \mathrm{C}^1 (\R^{d
\times d}\times \mathrm{K}';\R)$, and
\[
\begin{aligned}
& w_\mathrm{co}(F,z) \geq 0, \quad w_\mathrm{co}(RF,z)=w_\mathrm{co}(F,z),
 \quad w_\mathrm{co}(\cdot,z) \text{ is polyconvex,}
 \\
 &
 |\mathrm{D}_z w_\mathrm{co} (F,z)| \leq \kappa (|F|+1)^{p/2},
\\
&
 |\mathrm{D}_z w_\mathrm{co} (F,z_1)-\mathrm{D}_z w_\mathrm{co} (F,z_2) | \leq \kappa |z_1-z_2|^{\alpha}(|F|+1)^{p/2}
\end{aligned}
\]
for all arguments.

 In magnetism (see \cite{mie-generalized}), $z$
denotes the magnetization (with respect to material coordinates),
and we have $z \in \R^d$ and $\mathrm{K}=\{ z\in \R^d \, : \ |z|
\leq z_{\mathrm{sat}}\} $, where  the subscript ``$\mathrm{sat}$''
stands for \emph{saturation}. A choice for the coupling energy for
$p\geq 4$ is $w_\mathrm{co} (F,z)= C |F z|^2 +\tilde{w}(z)$, where
$\tilde{w} \in \mathrm{C}^2 (\mathrm{K}')$ gives the anisotropy of
magnetization, as well
 as the saturation term $\frac1{4\delta} (|z|^2 - z_{\mathrm{sat}}^2)^2$, with $\delta>0$.

For shape memory alloys, $z$
may denote volume fractions of $m$ different phases, such that
$\mathrm{K}= \{ z = (z_1,z_2,\ldots,z_m) \in \R^m\, : \ z_j \leq 0, \ \sum_{k=1}^m z_k=1 \}$.
Denoting by $\mathrm{cof} F \in \R^{d \times d}$ the cofactor matrix $\mathrm{det}(F) F^{-T}$
(which is contained in $\mathbb{M}(F)$), and by  $C_n(z)$, $n=1,\ldots,N$ the
$z$-dependent, effective
transformation Cauchy strains,
we may use
\[
w_\mathrm{co} (F,z) = \sum_{n=1}^N \alpha_n |F C_n(z)^{-1} -\mathrm{cof} F |^2 + w_{\mathrm{mix}}(z),
\]
with $\alpha_1,\ldots,\alpha_N>0$ and a mixture energy
$w_{\mathrm{mix}} \in \mathrm{C}^2(\mathrm{K}')$, see
\cite{KrMiRo05MMES,GoMiHa02FEMV}. Here we need $p \geq 2d$, because
$|\mathrm{D}_z w_\mathrm{co} (F,z)| \leq C (|F|+1)^d$, as the
highest power $|\mathrm{cof} F |^2 \sim O(|F|^{2d-2})$ is
independent of $z$. Note that we follow the ideas in
\cite{GoMiHa02FEMV,MiThLe02VFRI}, where $W(\cdot,z)$ is considered
to be a polyconvex relaxation, under given volume fractions of the
different phases.
\end{example}

\begin{example}
\label{ex:concrete-diss} Most commonly, the dissipation potentials
$\psi$ are assumed to be independent of the state $z$, i.e.\
$\psi(z,v)=\psi(v)$, which simplifies the analysis considerably.
However, there are cases where $\psi$ must depend on $z$,  like in
finite-strain elasticity where the internal variable is the plastic
tensor $P \in \mathrm{SL}(d)= \{P\in \R^{d \times d}\, : \
\mathrm{det}(P)=1 \}$, and $\psi_{P}(\dot{P})= \hat{\psi}(\dot{P}
P^{-1})$. In the framework of the modeling for magnetization
illustrated in Example \ref{ex:concrete}, we may consider $\psi:
\R^d \times \R^d \to [0,\infty)$ of the form
\[
\psi(z,v) = \psi_{\mathrm{rad}}(z \cdot v) +
\psi_{\mathrm{tang}}\left(\left(\mathrm{I}-\tfrac{1}{|z|^2} z
\otimes z\right)v\right),
\]
to account for different dissipations for enlarging the magnetization or changing its orientation.
\end{example}

\subsection{Proof of Theorem \ref{th:exist-marginal-elast}}
\paragraph{\textbf{Outline of the proof.}}
We  follow an abstract approach to the analysis of \eqref{pde}, by
rephrasing it as a  doubly nonlinear equation of the type
\eqref{eq:1-bis}, generated by the  generalized gradient system
$(V,\cE,\Psi,\diffname, \Ptname)$ specified in the following
\medskip lines.
\paragraph{
\textbf{Space} $\V$:} We choose
\[
\text{the ambient space \  $\V=L^2 (\Omega;\R^m)$.}
\smallskip
\]
\paragraph{\textbf{Energy} ${\cE}$:} We consider
 the \emph{reduced functional} $\cE:
[0,T]\times W^{1,q}(\Omega;\R^m) \to
(-\infty,+\infty]$ obtained by minimizing out the displacements from
$\marginale$, i.e.
\begin{align} \label{reduced-energy}
\begin{gathered}
\ene tz:= \inf\{\marginal t z{\phi}\, : \ \phi\in \calF\}
\
  \text{ with domain $[0,T]\times \domainenergy$, where }   \\
  \domainenergy =\{z \in W^{1,q}(\Omega;\R^m)\, : \ z(x) \in
\mathrm{K}
  \ \foraa\, x \in \Omega\}.
  \end{gathered}
\end{align}
 We  often use the
 decomposition of $\cE$ as a sum of a convex  and of a nonconvex,
 reduced functional
\begin{equation}
\label{decomp-energy}
 \ene tz = \cE^1(z) + \inf\{ \marginali 2 t \phi  \,
: \ \phi\in \calF\} \doteq \cE^1(z) + \enei 2tz \quad \text{for
$(t,z) \in [0,T]\times \domainenergy$.}
\end{equation}
Indeed, in Lemma \ref{l:1} below we  prove that for all $(t,z)\in
[0,T]\times \domainenergy$, the set
\begin{equation}
\label{argmin-not-empty-elast} M(t,z): = \argmin_{\phi\in \calF}\{ \marginal
t z \phi  \}  \ \text{ is nonempty and weakly compact  in
$W^{1,p}(\Omega;\R^d)$.}
\smallskip
\end{equation}
\paragraph{\textbf{Subdifferential} ${\diffname}$:}  Reflecting
\eqref{decomp-energy}, we  use the following subdifferential notion
\begin{equation}
\label{subdiff-marginal} \diff tz := \partial \cE^1 (z) + \margsub
\enei 2tz \quad \text{for all $(t,z) \in [0,T]\times\domainenergy$,}
\end{equation}
where, as in  Section~\ref{s:2new},    $\margsub \enei 2tz $  is the
\emph{marginal subdifferential} of the reduced energy $\cE^2$, viz.
\[
\margsub \enei 2tz = \left\{ \rmD_z \marginali 2t z \phi \, :
\ \phi\in M(t,z) \right\},
\]
with $\rmD_z \marginali 2 t{\cdot}{\phi} $ the G\^ateau derivative
of the functional \smallskip $\marginali 2 t{\cdot}{\phi}$.

\paragraph{\textbf{Generalized time-derivative} ${\Ptname}$:} We set
\[R(t,z,\xi):= \{\phimin \in M(t,z)\, : \ \xi \in \partial \cE^1 (z)  +   \rmD_z \marginali 2t z \phimin \}
\quad \text{for all  $(t,z,\xi) \in
\mathrm{graph}(\diffname)$}
\]
and define
\begin{equation}
\label{def-P} \Pt tz{\xi} :=\max_{\phimin \in R(t,z,\xi)}
\pairing{}{W^{1,p}}{-\ell'(t)}{\phimin}.
\smallskip
\end{equation}
\paragraph{\textbf{Dissipation potential} ${\Psi}$:} We  consider  the Finsler family $(\Psi_z)_{z \in
\domainenergy}$ of dissipation potentials
\begin{equation}
\label{psiz} \Psi_z: V \to [0,+\infty) \  \text{ defined by } \
\Psi_z(v) := \int_{\Omega}\psi(z(x),v(x)) \dd x.
\end{equation}

In what follows, throughout  Lemmas~\ref{l:1}--\ref{l:chain-rule} we
 check that the above generalized gradient system
$(V,\cE,\Psi,\diffname,\Ptname)$
 complies
with the abstract assumptions \eqref{e:4.1}--\eqref{mosco},
\eqref{Ezero-bis}--\eqref{eq:468} of Theorem \ref{thm:viscous2}. The
latter result yields the existence of  a pair $(u,\xi)$ fulfilling
the Cauchy problem for \eqref{xi-equation}, and the energy identity
\eqref{enid1}. The forthcoming calculations show that, in the
present setting, \eqref{xi-equation} and \eqref{enid1} entail
\eqref{pde} and \eqref{enid-marginal}.
\begin{remark}
\label{rem:interplay}
\upshape
As it will be clear from the ensuing calculations,
it is  possible to generalize the theory to the case where,
in place of \eqref{d4},
we have  for some $r \in (1,\infty)$
\[
\psi(z,v) \geq \kappa_7 |v|^r -\kappa_9, \quad \psi(z,v) \geq \kappa_8 |v|^{r'} -\kappa_9
\]
(where $r'=r/(r-1)$ is the conjugate exponent of $r$),
and the growth conditions in \eqref{w3} are replaced by
\[
\begin{aligned}
 & |\mathrm{D}_z W(F,z)| \leq \kappa_3 ( W(F,z) +\kappa_4)^{1-1/r},
\\
&
 |\rmD_z W(F,z_1) - \rmD_z W(F,z_2)| \leq \kappa_5 |z_1
-z_2|^{\alpha} (W(F,z_1) +\kappa_6)^{1-1/r}.
\end{aligned}
\]
Under these assumptions, it is  again possible to develop the
abstract approach of Section \ref{s:3}.   The  natural ambient space
is now $V= L^r(\Omega;\R^m)$ and,   as in Theorem
\ref{th:exist-marginal-elast},  one concludes the existence of a pair
$(z,\xi)$  fulfilling
  $\dot{z} \in L^r (0,T;L^r(\Omega;\R^m))$,  $\xi \in L^{r'} (0,T;L^{r'}(\Omega;\R^m))$,
  and satisfying \eqref{e:pde-xi} and \eqref{enid-marginal}.
  \end{remark} 
\paragraph{\textbf{Coercivity and time-dependence of $\cE$.}}
\begin{lemma}
\label{l:1} Assume \eqref{load}, \eqref{w2}--\eqref{w3}, and
\eqref{hyp-phi-dir}. Then \eqref{argmin-not-empty-elast} holds.
Moreover, there exist positive constants  $c_1,\ldots,c_6>0$ such
that for all $(t,z) \in [0,T]\times D$ and all $ \phimin \in M(t,z)$
we have
\begin{align}
\label{e2-lower-order}
&
c_1 \|\phimin
\|_{W^{1,p}(\Omega;\R^d)}^p -c_2 \leq \enei 2tz \leq c_3,
\\
\label{even-better} &   \int_{\Omega} W(\nabla \phimin(x) ,z(x)) \dd
x \leq c_4,
\\
 \label{stima2}
&
 \ene tz \geq
c_5 \| z \|_{W^{1,q}(\Omega;\R^m)}^q  -c_6\,.
\end{align}
Further, for a sufficiently large constant $E_0$ (cf.
\eqref{stored-energy}),   the energy functional $\cE$ is bounded
from below by a positive constant, it complies with
\eqref{Ezero-bis} and~\eqref{eq:17-bis}, and for every $(z_n),\, z
\subset L^2(\Omega;\R^m)$ we have
\begin{equation}
\label{fortiori} \left(z_n \weakto z \ \text{ in $ L^2(\Omega;\R^m)$
and } \ \sup_n
 \cE^1 (z_n)<+\infty
\right) \ \Rightarrow \
 z_n \to z \quad \text{ in $\mathrm{C}^0
(\overline{\Omega};\R^m)$.}
\end{equation}

 Moreover,
\begin{equation}
\label{estimate-time} \exists\, c_7>0 \ \ \forall \, t,\, s \in
[0,T] \ \  \forall\, z \in \domainenergy \, : \  \ |\ene tz - \ene
sz| \leq c_7|t-s|.
\end{equation}
Hence,
  $\cE$ fulfills~\eqref{eq:diffclass_a}.
\end{lemma}
\PROOF
 We
have  for every $(t,\phi,z) \in [0,T]\times \calF\times
L^2(\Omega;\R^m)$
\begin{equation}
\label{pass-1}
\begin{aligned}
\marginal tz {\phi}& \geq \kappa_1 \int_{\Omega}|\nabla\phi(x)|^p \, \d x
-\kappa_2 |\Omega| - \| \ell(t)\|_{W^{1,p}(\Omega;\R^d)^*} \, \|
 {\phi}\|_{W^{1,p}(\Omega;\R^d)} \\
 & \geq
 c_1 \|\phi\|_{W^{1,p}(\Omega;\R^d)}^p -\kappa_2 |\Omega| -
C' \| \ell \|_{L^\infty (0,T;W^{1,p}(\Omega;\R^d)^*}^{p'},
\end{aligned}
\end{equation}
where the first inequality follows from the positivity of the
functional $\cE^1$ and from
 \eqref{w1}, and the second
one from
 Poincar\'e's and Young's inequalities. Taking into
account \eqref{load}, we deduce the lower estimate in
\eqref{e2-lower-order}. Hence, it is sufficient to choose
$E_0:=2c_2$ in order to have $\cE$ bounded from below by a positive
constant.

 Next, we remark that the functional
$\marginal t z{\cdot}$ is (sequentially) lower semicontinuous
 with respect to the weak topology of
$W^{1,p}(\Omega;\R^d)$. Indeed, let $(\phi_k)_k$ weakly converge
to some $\phi \in W^{1,p}(\Omega;\R^d)$ as $k \to \infty$.
 Then, by the weak continuity of minors of
gradients (cf.\ \cite{Ball76CCET,Resh67SCMM}), $\bbM (\nabla \phi_k) \weakto
 \bbM (\nabla \phi)$ in $L^{p/d}
(\Omega;\R^{\mu_d})$. Taking into account the polyconvexity
assumption \eqref{w1}, we ultimately have $\liminf_{k \to \infty} \marginal  t z{\phi_k}
\geq \marginal
t z{\phi}$.
We combine this weak  lower semicontinuity property with
 the coercivity estimate
\eqref{pass-1},  and thus we  conclude that the set of minimizers
\eqref{argmin-not-empty-elast} is not empty via the direct method of the calculus
of variations.

 Secondly, we observe that
\[
\begin{aligned}
\enei 2tz =  \min_{\phi \in \calF} \marginali 2 t z  {\phi}
 & \leq   \int_{\Omega}
W(\nabla \phi_{\mathrm{Dir}}(x),z(x)) \,\d x - \pairing{}{W^{1,p}}{\ell(t)}{\phi_{\mathrm{Dir}}}
\\  &  \leq \int_{\Omega}
\max_{z \in \mathrm{K}} W(\nabla
\phi_{\mathrm{Dir}}(x),z)  \,\d x + C
\doteq c_3
\end{aligned}
\]
where the last inequality follows from \eqref{load} and \eqref{hyp-phi-dir}. Hence,
the upper estimate in \eqref{e2-lower-order} ensues.
Then, \eqref{even-better} follows from
\[
\begin{aligned}
\int_{\Omega} W(\nabla \phimin,z) \dd x
 & \leq \enei 2tz + \kappa_2
|\Omega| +  C \| \ell \|_{L^\infty(0,T;W^{1,p}(\Omega;\R^d)^*)}
\|\phimin\|_{W^{1,p}(\Omega;\R^d)}\\ &  \leq  c_3 + \kappa_2
|\Omega| + \left(\frac{c_3 + c_2}{c_1}\right)^{1/p} \| \ell
\|_{L^\infty(0,T;W^{1,p}(\Omega;\R^d)^*)},
\end{aligned}
\]
where the first inequality is due to \eqref{pass-1} and \eqref{load}, and the second one to
\eqref{e2-lower-order}.

Next, in view of \eqref{pass-1} we have for all $(t,z) \in
[0,T]\times L^2(\Omega;\R^m) $ and for every $\phimin \in M(t,z)$
\[
\begin{aligned}
\ene tz     \geq \frac1q\|\nabla z \|_{L^q(\Omega)}^q +\int_{\Omega}
I_{\mathrm{K}}(z(x)) \dd x + c_1
\|\phimin\|_{W^{1,p}(\Omega;\R^d)}^p -C.
\end{aligned}
\]
Then, \eqref{stima2} ensues from the  Poincar\'e inequality, and
\eqref{fortiori} follows from \eqref{stima2} and the fact that
$q>d$, hence $W^{1,q}(\Omega;\R^m) \Subset C^0
(\overline{\Omega};\R^m)$.

To prove~\eqref{eq:diffclass_a}, we observe that for all $z\in
\domainenergy$, for every $0 \leq s \leq t \leq T$ and every $\phiminv t \in M(t,z)$ and
$\phiminv s \in M(s,z)$ there holds
\[
\begin{array}{ll} \ene tz - \ene sz  = \enei2tz - \enei2sz &  = \marginali 2tz{\phiminv t} - \marginali 2sz{\phiminv s} \\  &
\leq   \marginali 2tz{\phiminv s} - \marginali 2sz{\phiminv s}
\\ &
\leq
-\pairing{}{W^{1,p}}{\ell(t)-\ell(s)}{\phiminv
s}\\ & \leq C \| \ell(t)
-\ell(s) \|_{W^{1,p}(\Omega;\R^d)^*} \|\phiminv s \|_{W^{1,p}(\Omega;\R^d)} \\
&  \leq C \| \ell'\|_{L^\infty (0,T;W^{1,p}(\Omega;\R^d)^*)}\,
|t-s|\, c_1^{-1/p}(\enei 2sz +c_2)^{1/p} \leq c_7 |t-s|
\end{array}
\]
where  we have  used~\eqref{load}
and~\eqref{e2-lower-order}. Exchanging the roles of $s$ and $t$,
  we
infer~\eqref{estimate-time}.
 \QED

\paragraph{\textbf{Properties of the dissipation potentials $(\Psi_z)_{z \in
\domainenergy}$.}}
\begin{lemma}
\label{lemma:diss} Assume \eqref{load}, \eqref{w2}--\eqref{w3}, and
\eqref{d1}--\eqref{d4}. Then, the dissipation potentials
$(\Psi_z)_{z \in \domainenergy}$ satisfy \eqref{e:4.1},
\eqref{mosco}, and for every $z \in \domainenergy$ we have
\begin{equation}
\label{replacement} \left\{
\begin{array}{ll}
 \Psi_z(v) \geq \kappa_7
\|v\|_{L^2(\Omega;\R^m)}^2  -\kappa_9|\Omega|
\\
\Psi_z^{*} (w) \geq  \kappa_8 \|w\|_{L^2(\Omega;\R^m)}^2
-\kappa_9|\Omega|  \end{array} \right. \quad \text{for every
$v,\,w\in L^2(\Omega;\R^m) $,}
\end{equation}
where $\kappa_7, \, \kappa_8$ and $\kappa_9$ are the same constants as in \eqref{d4}.
 Thus, \eqref{eq:41.1} is fulfilled.
\end{lemma}
 \PROOF It follows from  \cite[Prop.\,2.16,
p.\,47]{Brezis73} that, for every $z \in \domainenergy$ the
subdifferential and conjugate of the potential $\Psi_z$ are given
for all $v, \, w \in L^2 (\Omega;\R^m)$ by
\begin{equation}
\label{subdif-conj}
\begin{aligned}
\left\{ \begin{array}{ll} w \in \partial \Psi_z (v)\
\Leftrightarrow \ w(x) \in
\partial_v \psi(z(x),v(x)) \ \ \foraa\, x \in \Omega, \\  \Psi_z^*
(w) = \int_{\Omega} \psi^* (z(x),w(x)) \dd x.
\end{array}
\right.
\end{aligned}
\end{equation}
Hence, \eqref{d2} yields that $(\Psi_z)_{z \in \domainenergy}$ is a
family of \emph{admissible dissipation potentials} on $L^2
(\Omega;\R^m)$ in the sense of
\eqref{eq:psi-sum-1}--\eqref{eq:psi-sum-2}, and \eqref{d4} obviously
implies \eqref{replacement}. Finally,  exploiting \eqref{d1},
\eqref{d2}, \eqref{fortiori},  and relying on Ioffe's theorem
\cite{Ioff77LSIF}, it is not difficult to check that the first of
\eqref{mosco} and \eqref{e:impo-consequence}  are fulfilled. This
implies \eqref{mosco}.
 \QED

\paragraph{\textbf{Closedness and variational sum rule.}}
We need the following preliminary result.
\begin{lemma}[Subdifferentials]
\label{l:subdifferentials} Assume \eqref{load},
\eqref{w2}--\eqref{w3}, and \eqref{hyp-phi-dir}. Then:
\begin{enumerate}
\item the subdifferential of $\cE^1$ is
\begin{equation}
\label{frsub1}
\begin{gathered}
\frsub\cE^1(z)  =-\Delta_q z + \partial I_{\mathrm{K}}(z) \ \text{
for all } z \in \mathrm{dom}(\frsub \cE^1), \end{gathered}
\end{equation}
with domain described by the following conditions
\begin{equation}
\label{dom-frsub1} z \in \dom(\partial \cE^1) \ \Leftrightarrow \
\left\{
\begin{array}{ll}
z \in \domainenergy \subset W^{1,q}(\Omega;\R^m),
\\
-\Delta_q z \in L^2 (\Omega;\R^m),
\end{array}
 \right.
\end{equation}
hence
\begin{equation}
\label{p-lapl-regularity} \dom(\partial \cE^1) \subset
W^{\nu,q}(\Omega;\R^m) \quad \text{for all $\nu  \in
\left[1,1+\tfrac1q\right)$.}
\end{equation}
\item  There exists a constant $c_8$ such that
for every $(t,z) \in [0,T]\times \domainenergy $ and  $\phimin \in M(t,z)$ we have
 $ \rmD_z W(\nabla \phimin, z)
\in L^2(\Omega;\R^m)$, with
\begin{equation}
\label{even-better2}
\|\rmD_z W(\nabla \phimin,
z)\|_{L^2(\Omega;\R^m)} \leq c_8.
\end{equation}
 Hence the \emph{marginal} subdifferential
\begin{equation}
\label{gateau2} \margsub \enei 2tz = \left\{ \rmD_z W(\nabla
\phimin, z)\,: \ \phimin \in M(t,z)
\right\} \quad \text{ is bounded in $L^2 (\Omega;\R^m)$.}
\end{equation}
\item For all $t \in [0,T]$ and all  $z_1,\, z_2 \in \domainenergy$ with
$z_1(x), \, z_2(x)\in \mathrm{K}'$ for all $x \in \Omega$,  and for
every $\phiminv1 \in M(t,z_1)$ there holds
\begin{equation}
\label{useful-preliminary}
\begin{aligned}
    \enei 2t{z_2}  & - \enei 2t{z_1}   -  \int_{\Omega} \rmD_z W(\nabla
\phiminv1,z_1) (z_2-z_1) \dd x
\\
      &  \quad \leq \kappa_5 \|z_1 -z_2\|_{L^\infty
(\Omega;\R^m)}^{\alpha} \left(\int_{\Omega} W(\nabla
 \phiminv1,z_1) \dd x +
\kappa_6 |\Omega|\right)^{1/2} \|z_1 -z_2\|_{L^2 (\Omega;\R^m)} \\
\end{aligned}
\end{equation}
where $\kappa_5$, $\kappa_6$, and $\alpha$ are the same constants as
in \eqref{w3}.
\item For every $(t,z) \in [0,T]\times \domainenergy$
the Fr\'echet subdifferential $\frsub \cE_t$ satisfies
\begin{equation}
\label{inclu-diff} \frsub \ene tz \subset \diff tz =  -\Delta_q z
+\partial I_{\mathrm{K}}(z) + \{\rmD_z W(\nabla  \phimin, z)\, : \ \phimin \in M(t,z) \}\,.
\end{equation}
\end{enumerate}
\end{lemma}
\PROOF Formulae \eqref{frsub1} and \eqref{dom-frsub1} can be
obtained by  adapting the proof of
\cite[Lemma\,2.4]{Schimperna-Segatti-Stefanelli07}, see also
\cite[Prop.\,2.17]{Brezis73}. Notice that \eqref{p-lapl-regularity}
ensues from \eqref{dom-frsub1} and the regularity results in
\cite{Savare98}, cf.\ also \cite{EbmFre99}. We conclude
\eqref{even-better2} combining  condition \eqref{w3}(i) with
estimate \eqref{even-better}, and then \eqref{gateau2} follows from
trivial calculations.

Estimate \eqref{useful-preliminary} is a consequence of the
following chain of  inequalities
\begin{equation}
\label{ingredient1}
\begin{aligned}
\enei 2t{z_2} - \enei 2t{z_1}  & \leq
\int_{\Omega} \left( W(\nabla
\phiminv1, z_2) - W(\nabla  \phiminv1, z_1) \right)\dd x \\
& = \int_{\Omega}\int_0^1 \rmD_z W(\nabla \phiminv1, (1-\theta)z_1 +\theta z_2) (z_2-z_1) \dd \theta \dd x \\
& \leq I+ \int_{\Omega} \rmD_z W(\nabla
\phiminv1,z_1) (z_2-z_1) \dd x
\end{aligned}
\end{equation}
where, relying on  condition \eqref{w3}(ii), we estimate
\begin{equation}
\label{ingredient2}
\begin{aligned}
I & = \int_0^1 \int_{\Omega} |\rmD_z W(\nabla \phiminv1, (1-\theta)z_1 +\theta z_2) - \rmD_z W(\nabla
 \phiminv1, z_1)| |z_2-z_1| \dd x \dd
\theta
\\ & \leq \int_0^1 \int_{\Omega} \kappa_5 \theta^{\alpha}
|z_2-z_1|^{\alpha} \left(W(\nabla
\phiminv1,z_1) + \kappa_6 \right)^{1/2} |z_2-z_1| \dd x \dd \theta.
\end{aligned}
\end{equation}
Then,  \eqref{useful-preliminary} follows upon using H\"older's
inequality.

Finally, we  prove \eqref{inclu-diff},  in fact in the following
stronger form
\begin{equation}
\nonumber
 \text{if } \frsub \ene tz  \neq \emptyset, \
\text{then } \ \xi - \rmD_z W(\nabla
\phimin,z) \in \frsub \cE^1(z)
 \   \text{for
every $\xi \in \frsub \ene tz$ and $\phimin \in M(t,z)$.}
\end{equation}
which in particular yields \eqref{inclu-diff}. Indeed, we  show that
for every $\xi \in \frsub \ene tz$ and $\phimin \in M(t,z)$,
 and for every $z_n \to z $ in $L^2 (\Omega;\R^m)$, there holds
 \begin{equation}
\label{claim}
\begin{aligned}
  0 \leq \liminf_{n \to \infty} \frac{\cE^1(z_n) - \cE^1({z}) -
  \int_{\Omega} \left(\xi - \rmD_z W(\nabla
\phimin,z)\right)
 ( z_n - z )\, \d x }{\|z_n - z\|_{L^2 (\Omega;\R^m)}} \doteq \Lambda.
 \end{aligned}
\end{equation}
To this aim, we observe that
\begin{equation}
\label{auxiliary-label}
\begin{aligned}
\Lambda \geq \liminf_{n \to \infty}   & \frac{\ene{t}{z_n} - \ene{t}{z}
-  \int_{\Omega}\xi
 (z_n-z) \, \d x}{\|z_n - z\|_{L^2 (\Omega;\R^m)}} \\ & + \liminf_{n \to
 \infty} \frac{\enei 2{t}{z} -   \enei
 2{t}{z_n} + \int_{\Omega}\rmD_z W(\nabla
\phimin,z) ( z_n - z )\, \d x}{\|z_n - z\|_{L^2 (\Omega;\R^m)}}.
 \end{aligned}
 \end{equation}
Since the first summand on the right-hand side of the above
inequality is nonnegative by  definition of the Fr\'echet
subdifferential $\frsub \ene tz$, it remains to prove that the
second term   is nonnegative. Now, it is not restrictive to suppose
for the sequence $(z_n)$ in \eqref{claim} that $\sup
 \cE^1(z_n)<+\infty$. Then,
 $z_n(x) \in \mathrm{K}$ for all $x \in \Omega$ and  $n \in \N$.
 Hence, estimate
\eqref{useful-preliminary} with the choices $z_1=z$ and $z_2=z_n$ yields
\[
\begin{aligned}
 & \liminf_{n \to
 \infty} \frac{\enei 2{t}{z} -   \enei
 2{t}{z_n} + \int_{\Omega}\rmD_z W(\nabla
\phimin,z) ( z_n - z )\, \d x}{\|z_n - z\|_{L^2
 (\Omega;\R^m)}} \\ & \geq -  \kappa_5  \lim_{n \to \infty}
\frac{\|z_n -z\|_{L^\infty (\Omega;\R^m)}^{\alpha}
\left(\int_{\Omega} W(\nabla  \phimin,z)
\dd x + \kappa_6\right)^{1/2} \|z_n -z\|_{L^2 (\Omega;\R^m)} }{\|z_n
- z\|_{L^2
 (\Omega;\R^m)}} =0,
 \end{aligned}
\]
and the last limit follows from \eqref{fortiori} and the bound
\eqref{even-better}. Ultimately, \eqref{claim} ensues.
 \QED
\begin{lemma}[Closedness]
\label{closedness} Assume~\eqref{load}, \eqref{w2}--\eqref{w3}, and
\eqref{hyp-phi-dir}.
 Then, for all $\{ t_n \} \subset [0,T]$, $\{z_n\}\subset
L^2(\Omega;\R^m)$, and $\{ \xi_n\}\subset L^2(\Omega;\R^m) $ with
$\xi_n \in \diff{t_n}{z_n}$ for all $n \in \N$, we have
\begin{equation}
\label{contin-prop}
\begin{aligned}
 & \left( t_n \to t, \ \   z_n \weakto z \
\text{in $L^2(\Omega;\R^m)$,}\ \ \xi_n\weakto \xi \ \text{in
$L^2(\Omega;\R^m)$,}\ \  \ene{t_n}{z_n} \to \EE   \ \ \text{as $n
\to \infty$}\right)
\\ &
\Longrightarrow  \ \xi \in \diff{t}z \ \ \text{and } \ \EE= \ene tz.
\end{aligned}
\end{equation}
 In particular, $\mathrm{graph}(\diffname)$ is a Borel set of
$[0,T]\times L^2(\Omega;\R^m)\times L^2(\Omega;\R^m)$.
\end{lemma}
\PROOF From $ \sup_n \ene{t_n}{z_n} <+\infty$ and from
\eqref{e2-lower-order},  \eqref{stima2}, and \eqref{fortiori},  we
deduce that $z_n \weakto z$ in $W^{1,q}(\Omega;\R^m)$ and $z_n \to
z$ in $\mathrm{C}^0(\overline{\Omega};\R^m)$,
 and that there exist
$\phi \in \calF$ and a (not relabeled) subsequence $(\phiminn)$ such
that $\phiminn \weakto \phi$ in $W^{1,p}(\Omega;\R^d)$ as $n \to
\infty$. Hence, we argue in the same way as in the proof of Lemma
\ref{l:1}:   combining the polyconvexity assumption \eqref{w1} with
the continuity of the map $z\mapsto W(F,z)$, we apply Ioffe's
theorem \cite{Ioff77LSIF}
 to find  $
 \liminf_{n \to \infty}
 \int_\Omega
W(\nabla \phiminn,z_n) \dd x \geq
 \int_\Omega
W(\nabla \phi,z) \dd x.$ Therefore,
in view of \eqref{load} we have
\begin{equation}
\label{liminf-nk}
\begin{aligned}
 \liminf_{n \to \infty} \enei 2{t_n}{z_n}  & = \liminf_{n \to \infty}
 \left(\int_\Omega
W(\nabla \phiminn,z_n) \dd x -
\pairing{}{W^{1,p}}{\ell(t_n)}{ \phiminn}
\right)
\\ & \geq \int_\Omega W(\nabla\phi,z) \dd x -
\pairing{}{W^{1,p}}{\ell(t)}{\phi} \geq \enei 2 tz.
\end{aligned}
\end{equation}
On the other hand, $\phiminn \in M(t_n,z_n)$
gives
\begin{equation}
\label{comparison-n}
\enei 2{t_n}{z_n}  = \marginali 2{t_n}{z_n}{\phiminn} \leq \marginali 2{t_n}{z_n}{\phimin}
\to \marginali 2{t}{z}{\phimin} = \enei 2 tz,
\end{equation}
where $\phimin$ is any element in $M(t,z)$, and we have exploited \eqref{load} to take
the limit as $n \to \infty$. Combining \eqref{liminf-nk} and \eqref{comparison-n},
 we ultimately have
 $\enei2 {t_n}{z_n} \to \enei2tz$, and  the weak limit $\phi$ of the sequence $(\phiminn)$
 is in fact an element in $M(t,z)$, which we will hereafter denote with $\varphi$.

Now, it follows from \eqref{subdiff-marginal}, \eqref{frsub1}, and
\eqref{gateau2} that the sequence $\xi_n \in \diff{t_n}{z_n}$ in
\eqref{contin-prop} is given, for every $n \in \N$, by $\xi_n =
-\Delta_q z_n +\zeta_n + \rmD_z W(\nabla\phiminn, z_n), $ for some
$\zeta_n \in
\partial I_{\mathrm{K}}(z_n)$ and $\phiminn \in M(t_n,z_n)$.
  Arguing by
comparison and  relying on   the aforementioned
\cite[Prop.\,2.17]{Brezis73}, from the boundedness of $\xi_n $
 in $L^2(\Omega;\R^m)$  we infer that
\begin{equation}
\label{key} \sup_n \left(\| \Delta_q z_n\|_{L^2 (\Omega;\R^m)} + \|
\zeta_n  \|_{L^2 (\Omega;\R^m)} \right) <+\infty.
\end{equation}
Relying on~\cite{Savare98},  we find that
for every $\nu \in [1,1+1/q)$ there holds $\sup_n \|
z_n\|_{W^{\nu,q} (\Omega;\R^m)}<+\infty$.  Since $W^{\nu,q}
(\Omega;\R^m) \Subset W^{1,q}(\Omega;\R^m)$ for all $\nu\in
(1,1+1/q)$, we conclude that, indeed, the weak convergence of
$(z_n)$ in $L^2(\Omega;\R^m)$ improves to
\begin{equation}
\label{strong-w1q} z_n \to z \quad \text{ in
$W^{1,q}(\Omega;\R^m)$.}
\end{equation}
Therefore, $\cE^1(z_n) \to \cE^1(z)$. On account of
 the previously proved convergence of $\enei 2{t_n}{z_n}$,  we obtain $\ene {t_n}{z_n} \to \ene tz$.
Finally, combining estimate \eqref{key} with \eqref{strong-w1q}, and exploiting the monotonicity
of the operator $-\Delta_q$ (cf.\ \cite{Brezis73}),   we
find
\begin{equation}
\label{e:q-ident}
 -\Delta_q  z_n \weakto -\Delta_q z \ \  \text{ in $L^2
(\Omega;\R^m)$.}
\end{equation}
 Furthermore, from \eqref{key} we also deduce that,
up to a not relabeled subsequence,
\begin{equation}
\label{e:zeta-ident} \zeta_n \weakto \zeta \ \  \text{ in $L^2
(\Omega;\R^m)$, with } \ \zeta \in
\partial I_{\mathrm{K}}(z)
\end{equation}
(the latter fact follows from the strong-weak closedness of the
graph of $\partial I_{\mathrm{K}}$ in $L^2 (\Omega;\R^m) \times L^2
(\Omega;\R^m)$).

  Estimate
\eqref{even-better2} yields
\begin{equation}
\label{aux-lab2} \sup_n \| \rmD_z W(\nabla \phiminn, z_n) \|_{L^2
(\Omega;\R^m)} \leq c_8.
\end{equation}
 Then, along some
  (not relabeled) subsequence, the sequence  $ (\rmD_z W(\nabla \phiminn,
z_n) )_n$ is weakly converging  in $L^2(\Omega;\R^m)$.
It remains to prove that
\begin{equation}
\label{desired} \rmD_z W(\nabla \phiminn,
z_n) \weakto \rmD_z W(\nabla \varphi, z)
\ \ \text{ in $L^2 (\Omega;\R^m)$.}
\end{equation}
To this aim, we mimick the argument in the proof of
\cite[Prop.\,3.3]{FraMie06ERCR}. We fix $\eta \in
W^{1,q}(\Omega;\R^m)$ and $h>0$, and apply \eqref{w3}(ii) with the
choices $z_1=z_n$ and $z_2=z_n + h\eta$. Indeed, $z_n (x) \in
\mathrm{K}$ and  \eqref{strong-w1q} ensure that, for sufficiently
large $n$ and sufficiently small $h$, we have  $z_n,\, z_n + h\eta
\in \mathrm{K}'$. Arguing like in
  \eqref{ingredient1}--\eqref{ingredient2} and exploiting estimate~\eqref{even-better},
   there exist constants $C>0$ and
$\alpha \in (0,1]$ such that for every $n \in \N$
\begin{equation}
\label{framie1}
\begin{aligned}
\Big| \frac1h \int_{\Omega}\big( W(\nabla \phiminn, z_n \pm h \eta) & - W(\nabla \phiminn, z_n)
 \mp h \rmD_z W(\nabla \phiminn, z_n)\eta \big) \dd x\Big| \\ &  \leq C h^{\alpha}
\|\eta \|_{L^\infty (\Omega;\R^m)}^{\alpha} \|\eta \|_{L^2
(\Omega;\R^m)}\doteq \omega(h).
\end{aligned}
\end{equation}
 On the other hand, again combining \eqref{strong-w1q} and the weak
convergence of $\phiminn$ with Ioffe's theorem, we conclude that for
(sufficiently small) $h>0$ there holds
\[
\begin{aligned}
 \liminf_{n \to \infty}
\int_{\Omega}\big( W(\nabla \phiminn, z_n
\pm h \eta)   - W(\nabla \phiminn, z_n)
\big) \dd x \geq \frac1h \int_{\Omega}\big( W(\nabla
\varphi, z \pm h \eta)   -
W(\nabla \varphi, z) \big) \dd x.
\end{aligned}
\]
Estimate \eqref{framie1} and the above inequality yield
\begin{equation}
\label{modulus}
\begin{aligned}
 & \limsup_{n \to \infty}   \int_{\Omega} \rmD_z W(\nabla
\phiminn, z_n)\eta \dd x\\   &\quad  \leq
\limsup_{n \to \infty} \frac1h \int_{\Omega}\big( W(\nabla
\phiminn, z_n) - W(\nabla
\phiminn, z_n - h \eta) \big) \dd x +
\omega(h) \\ & \quad \leq -\frac1h \int_{\Omega}\big(
W(\nabla \varphi, z - h \eta) - W(\nabla
 \varphi, z) \big) \dd x + \omega(h)
\leq \int_{\Omega} \rmD_z W(\nabla
\varphi, z)\eta \dd x + 2\omega(h),
\end{aligned}
\end{equation}
where  the last inequality follows from \eqref{framie1} written for
$(\nabla \varphi,z)$. Analogously, we infer  that
\[
\liminf_{n \to \infty} \int_{\Omega} \rmD_z
W(\nabla\phiminn, z_n)\eta \dd x \geq
\int_{\Omega} \rmD_z W(\nabla \varphi,
z)\eta \dd x - 2\omega(h).
\]
 Since $h>0$ is arbitrary, we conclude that
 \[
 \lim_{n
\to \infty} \int_{\Omega} \rmD_z W(\nabla \phiminn, z_n)\eta \dd x= \int_{\Omega} \rmD_z W(\nabla
 \varphi, z)\eta \dd x \quad \text{for every $\eta
\in W^{1,q}(\Omega;\R^m)$.}
\]
In view of \eqref{aux-lab2},
\eqref{desired} follows. Thus,
\eqref{e:q-ident}, \eqref{e:zeta-ident} and \eqref{desired}
 entail that the weak limit $\xi$ of $(\xi_n)$
fulfills
$\xi \in \diff tz$, and \eqref{contin-prop} ensues.

Finally, let us observe that $\mathrm{graph}(\diffname)= \cup_{m\in \N}
G^m$, with
 \[   G^m = \left\{(t,u,\xi) \in [0,T] \times L^2(\Omega;\R^m) \times L^2(\Omega;\R^m) \, : \ \xi \in \diff tu,  \
|\ene{t}{u}| \leq m\right\}.
\]
Now, it follows from the  closedness property~\eqref{contin-prop}
that every $G_m$ is a closed, hence Borelian, set. Hence, $\mathrm{graph}(\diffname)$
is a Borel set.
 \QED 
\begin{corollary}[Variational sum rule]
\label{var-sum-rule}
 Assume~\eqref{load},  \eqref{w2}--\eqref{w3},
 \eqref{d1}--\eqref{d4}, and \eqref{hyp-phi-dir}.
Then, the dissipation potentials $(\Psi_z)_{z \in \domainenergy}$
and the reduced energy functional $\cE$ comply with the variational sum
rule~\eqref{eq:42-bis}.
\end{corollary}
\PROOF This follows from Lemmas~\ref{lemma:diss},
\ref{l:subdifferentials}, and \ref{l:subdifferentials}, combined
with Proposition~\ref{l:variational-sum-rule}.
 \QED
\paragraph{\textbf{Chain rule.}}
\begin{lemma}[Chain rule]
\label{l:chain-rule} Assume~\eqref{load}, \eqref{w2}--\eqref{w3},
\eqref{d1}--\eqref{d4}, and \eqref{hyp-phi-dir}. Then, the function
$\Ptname: \mathrm{graph}(\diffname) \to \R$ defined in \eqref{def-P}
complies with \eqref{e:ass_p_a}. Moreover, the system
$(V,\cE,\Psi,\diffname,\Ptname)$ fulfills the closedness
condition~\eqref{eq:468}, and the chain-rule
inequality~\eqref{eq:48strong}.
\end{lemma}
\PROOF We first observe that
\begin{equation}
\label{weakly-compact} R(t,z,\xi) \ \text{ is weakly sequentially
compact in $W^{1,p}(\Omega;\R^d)$ for every } (t,z,\xi) \in
\mathrm{graph}(\diffname).
\end{equation}
Indeed,  every sequence  $(\phiminn)_n \subset  R(t,z,\xi)$ is
bounded in $W^{1,p}(\Omega;\R^d)$ thanks to \eqref{e2-lower-order}.
Hence, up to a subsequence it converges to some $\varphi$. From the
arguments in Lemma \ref{closedness} it follows that $\varphi \in
R(t,z,\xi)$. Thus, it is immediate to see that the maximum in
formula~\eqref{def-P} is attained.

For every $(t,z) \in [0,T]\times W^{1,q}(\Omega)$, $h\in (0,T-t]$,
and $\phimin(t)\in M(t,z)$ there holds
\[
\frac{\ene{t+h}z -\ene tz}{h}= \frac{\enei 2{t+h}z -\enei
2tz}{h}\leq \frac{1}h \pairing{}{W^{1,p}}{-\ell(t+h)
+\ell(t)}{\phimin(t)},
\]
whence $\limsup_{h \down 0}\frac{\ene{t+h}z -\ene tz}{h}\leq \Pt
tz\xi$. On the other hand, it follows from~\eqref{load}
and~\eqref{e2-lower-order} that $ |\Pt tz\xi|   \leq  \|\ell'(t)
\|_{W^{1,p}(\Omega;\R^d)^*} \cdot \sup_{\phimin \in M(t,z)} \|
 \phimin \|_{W^{1,p}(\Omega;\R^d)} \leq C$.
Therefore, \eqref{e:ass_p_a} is fulfilled.

Combining  the previously proved
closedness property \eqref{contin-prop} with arguments analogous to
those developed for \eqref{weakly-compact}, it is possible to check  that
\eqref{eq:468} holds in
a slightly stronger form, viz.
\[
\begin{aligned}
&
\big(
t_n \to t, \
  u_n\to u  \text{ in $\V$,}
  \ \diff{t_n}{u_n} \ni \xi_n\weakto \xi   \text{ in
     $\V^*$, } \
    \Pt {t_n}{u_n}{\xi_n} \to p,  \ \ene {t_n}{u_n}\to \EE
    \big)
\\
& \Longrightarrow \
      (t,u)\in \dom(F),\quad
      \xi\in \diff tu,\quad
      p\leq \Pt{t}{u}{\xi},\quad
      \EE=\ene tu.
      \end{aligned}
      \]
Hence, mimicking the argument at the end of the proof of Lemma \ref{closedness},
  it is possible to check
   that for every $\lambda \in \R$, the set $\Ptname^{-1}([\lambda,+\infty))$ is a Borel
   set of $[0,T] \times \V \times \V^*$. Therefore,
   $\Ptname: \mathrm{graph}(\diffname) \to
\R$ defined by~\eqref{def-P} is a Borel function.

Finally, in order to prove that the chain rule~\eqref{eq:48strong}
is fulfilled, let us fix a curve $ z\in \AC ([0,T];L^2(\Omega))$ and
a function $\xi \in L^1 (0,T;L^2(\Omega))$
fulfilling~\eqref{e:basic-conditions} and \eqref{conditions-1}.
Taking into account \eqref{stima2} and \eqref{replacement}, we have
a fortiori that
\begin{equation}
\label{additional-regu}
\begin{gathered}
 z \in L^\infty (0,T; W^{1,q}(\Omega;\R^m))
\cap H^1(0,T; L^2 (\Omega;\R^m)) \subset \mathrm{C}^0
([0,T];L^\infty(\Omega;\R^m)), \\ \xi \in L^2
(0,T;L^2(\Omega;\R^m)).
\end{gathered}
\end{equation}
Furthermore,  there exist  measurable selections $t \mapsto \zeta(t)
\in \partial I_{\mathrm{K}}(z(t))$ and
 $t \mapsto \phimin(t) \in M(t,z(t))$ such that
\begin{equation}
\label{selection-xi} \xi(t) = -\Delta_q z(t) + \zeta(t) + \rmD_z
W(\nabla \phimin(t), z(t)) \quad \foraa\,
t \in (0,T). \end{equation} Arguing as in the proof of
Lemma~\ref{closedness}, from $\xi \in L^2
(0,T;L^2(\Omega;\R^m))$ we deduce that
\begin{equation}
\label{quoted-later}
\begin{aligned}
\!\!\!\!
 \| \Delta_q z(t)\|_{L^2(0,T;L^2(\Omega;\R^m))}
 + \| \zeta(t) \|_{L^2(0,T;L^2(\Omega;\R^m))}  +  \|\rmD_z W(\nabla
\phimin(t), z(t))
\|_{L^\infty(0,T;L^2(\Omega;\R^m))} \leq C,
\end{aligned}
\end{equation}
where the latter estimate follows from \eqref{even-better2}.
  Thus,  the chain rule for the
\emph{convex} functional $\cE^1$ (see \cite{Brezis73}) yields that
\begin{equation}
\label{chain-rule-1}
\begin{aligned}
   & \text{the map $t\mapsto \cE^1 (z(t)) $ is absolutely continuous,
and } \\ &    \frac{\d}{\d t} \cE^1 (z(t)) =\int_{\Omega} (-\Delta_q
z(t) + \zeta(t)) z'(t) \dd x \quad \foraa \, t \in (0,T).
\end{aligned}
\end{equation}
As for the map $t \mapsto \enei 2t{z(t)}$, there exist constants
$C>0$, $\alpha \in (0,1]$ such that
 for every $0 \leq s
\leq t \leq T$ we have
\[
\begin{aligned}
&
 \enei 2t{z(t)}  - \enei 2s{z(s)}
 \\ &  = \enei 2t{z(t)}  - \enei 2t{z(s)} + \enei 2t{z(s)}  - \enei
2s{z(s)}
 \\ &
  \leq \enei 2t{z(t)}  - \enei 2t{z(s)} +\marginali 2t {z(s)}{\phimin(s)} -\marginali 2s {z(s)}{\phimin(s)}
 \\  &
  \leq  C \|z(t) -z(s)\|_{L^\infty
(\Omega;\R^m)}^{\alpha}  \|z(t) -z(s)\|_{L^2 (\Omega;\R^m)}  \\ &
\quad  + \int_{\Omega} \rmD_z W(\nabla \phimin(s)),z(s)) (z(t)-z(s))
\dd x
 - \pairing{}{W^{1,p}}{\ell(t) - \ell(s)}{\phimin(s)}.
\end{aligned}
\]
where the second inequality follows from estimate
\eqref{useful-preliminary} with $z_1 = z(s)$ and $z_2 =z(t)$, also
taking into account \eqref{even-better}.
 Exchanging the
role of $s$ and $t$, we thus conclude for every $0 \leq s \leq t
\leq T$
\[
\begin{aligned}
 & | \enei 2t{z(t)}   - \enei 2s{z(s)} |  \\ &
  \leq   C \|z(t) - z(s)
\|_{L^2(\Omega;\R^m)}  \big(  \|z(t) -z(s)\|_{L^\infty
(\Omega;\R^m)}^{\alpha}     + 2\sup_{t \in (0,T)} \| \rmD_z W(\nabla
\phimin(t),z(t))\|_{L^2(\Omega;\R^m)} \big)
 \\ & \quad  +  \| \ell(t) -
\ell(s) \|_{W^{1,p}(\Omega;\R^d)^*} (\|
\phimin(t)\|_{W^{1,p}(\Omega;\R^d)} + \|\phimin(s)
\|_{W^{1,p}(\Omega;\R^d)} )
\\ & \leq C' \left(  \|z(t) - z(s)
\|_{L^2(\Omega;\R^m)} + |t-s| \right),
\end{aligned}
\]
where the second inequality  follows from \eqref{additional-regu},
\eqref{load},  and estimates
\eqref{e2-lower-order} and \eqref{quoted-later}. Thus, the map $t
\mapsto \enei 2t{z(t)}$ is absolutely continuous.

Finally, let us fix $t \in (0,T)$, such that formula
\eqref{chain-rule-1} for $\frac{\d}{\d t}\cE^1(z(t))$ holds,
$\frac{z(t+h) -z(t)} h \to z'(t)$ in $L^2 (\Omega)$,  $\ell(t+h)
-\ell(t) \to \ell'(t)$ in $W^{1,p}(\Omega;\R^d)^*$, and
$\frac{\d}{\d t}\enei 2t{z(t)}$ exists (the set of such $t$'s has
full measure). Now, in view of \eqref{useful-preliminary}, and again
taking into account \eqref{additional-regu} and \eqref{even-better},
for all $h \in (-t,0]$ and $\phiminvv t \in R(t,z(t),\xi(t))$ there
holds
\[
\begin{aligned}
 & \frac1h   (\enei 2{t+h}{z(t+h)} -\enei2 t{z(t)}) \\
 &  \geq    C
\|z(t+h) -z(t)\|_{L^\infty (\Omega;\R^m)}^{\alpha}  \frac{\|z(t+h)
-z(t)\|_{L^2 (\Omega;\R^m)}}{h}  \\
& \quad + \frac1h \int_{\Omega} \rmD_z W(\nabla \phiminvv t,z(t))
(z(t+h)-z(t)) \dd x -\frac1h \pairing{}{W^{1,p}}{\ell(t+h)
-\ell(t)}{\phiminvv t}.
\end{aligned}
\]
Taking the $\lim_{h \up 0}$ in the above inequality and using
that $z \in \mathrm{C}^0 ([0,T];L^\infty (\Omega;\R^m))$ by
\eqref{additional-regu}, we conclude that for every $\phiminvv t \in
R(t,z(t),\xi(t))$
\begin{equation}
\label{5.58} \frac{\d}{\d t}\enei 2t{z(t)} \geq \int_{\Omega} \rmD_z
W(\nabla \phiminvv t,z(t)) z'(t) \dd x -
\pairing{}{W^{1,p}}{\ell'(t) }{\phiminvv
t}. \end{equation}
 Now,
from the definition of $R(t,z(t),\xi(t))$ it follows that, in
correspondence   to the map $t\mapsto  \phiminvv t$, there exists a
selection $t \mapsto \tilde{\zeta}(t) \in \partial
I_{\mathrm{K}}(z(t))$ such that $\tilde{\zeta}(t) + \rmD_z
W(\nabla \phiminvv t,z(t)) = \zeta(t) +
\rmD_z W(\nabla \phimin (t),z(t))$ for
almost all $t \in (0,T)$ (where $\zeta$ and $\phimin$ are the
selections in \eqref{selection-xi}). Thus, using the chain rule for
$I_{\mathrm{K}}$, we have
\[
\begin{aligned}
   \int_{\Omega}  \rmD_z W(\nabla \phiminvv t,z(t)) z'(t)
\dd x = \int_{\Omega}  \rmD_z W(\nabla \phimin(t),z(t)) z'(t) \dd x  +\!\! \!\! \!\! \ddd{\int_{\Omega}
\zeta(t) z'(t) \dd x - \int_{\Omega} \tilde{\zeta(t)} z'(t) \dd
x}{=0}{}\!\! \!\!.
\end{aligned}
\]
Since the selection $t\mapsto \phiminvv t \in R(t,z(t),\xi(t))$ in
\eqref{5.58} is arbitrary, from  the above equality we
ultimately conclude
\begin{equation}
\label{chain-rule-2}  \frac{\d}{\d t} \enei2 t{z(t)} \geq
\int_{\Omega} \rmD_z W(\nabla\phimin(t),z(t)) z'(t) \dd x +  \Pt {t}{z(t)}{\xi(t)}  \quad
\foraa \, t \in (0,T).
\end{equation}
Combining \eqref{chain-rule-1} and \eqref{chain-rule-2}, we obtain
\medskip \eqref{eq:48strong}.
 \QED

\noindent Thus, we have shown that all the abstract assumptions of
Section \ref{ss:3.0} are fulfilled, which implies that Theorem
\ref{th:exist-marginal-elast} follows from Theorem
\ref{thm:viscous2}.

\section{Proofs}
\label{s:5} \paragraph{\textbf{Plan of the proof of
Theorem~\ref{thm:viscous2}.}} First, in Section \ref{ss:4.1}  we
provide some ``stationary estimates'' on every single step of the
incremental minimization scheme. In particular, in Lemma
\ref{lemma:degiorgi} we prove the crucial energy inequality
\eqref{e:interp-ene-ineq}, which is the starting point for the a
priori estimates on the approximate solutions. We prove the latter
estimates in Proposition \ref{prop:a-prio-est}.  Hence  we deduce in
Proposition \ref{prop:compactness-young} that, along some
subsequence, the approximate solutions converge to a curve $u \in
\AC ([0,T];\V)$. In Section \ref{ss:4.3} we conclude the proof of
Theorem~\ref{thm:viscous2}, showing that $u$ is in fact a solution
of the Cauchy problem for \eqref{eq:1-bis}. In doing so, we rely on
some
technical results proved in the Appendix. 
\subsection{Discrete energy inequality}
\label{ss:4.1} In the following, we gain further insight into
problems \eqref{e:min-scheme} and \eqref{interpmin} (which give rise
to  approximate solutions), by fixing some crucial properties of
the general minimization problem
\begin{equation}
\label{e:gen-min-prob} \IMIN{\mathsf{t}}{r}{u}:=\inf_{v \in
\domainenergy} \left\{ r \Psi_u \left(\frac{v-u}{r} \right) +
\ene{\mathsf{t}+r}{v} \right\} \qquad \text{for given $\mathsf{t}
\in [0,T]$, $u \in \domainenergy$, $0<r<T-\mathsf{t}$.}
\end{equation}
 The following
result is the Banach-space counterpart
to~\cite[Lemmas 4.4,4.5]{Rossi-Mielke-Savare08} (see
also~\cite{Ambrosio95, Ambrosio-Gigli-Savare08, Rossi-Savare06}).
\begin{lemma}
\label{lemma:degiorgi}
 Assume \eqref{e:4.1}, and \eqref{Ezero-bis}--\eqref{e:ass_p_a}.
 Then, for every $\mathsf{t} \in [0,T]$, $u \in \domainenergy$, and
  for all $0<r<T-\mathsf{t}$
\begin{equation}
\label{nonempty} \text{the set} \ \
\MIN{\mathsf{t}}{r}{u}:=\argmin_{v \in \domainenergy} \left\{ r
\Psi_u \left(\frac{v-u}{r} \right) + \ene{\mathsf{t}+r}{v} \right\}
\quad \text{is nonempty.}
\end{equation}
 Moreover,
  for all $\mathsf{t} \in [0,T]$ there exists
  a measurable
selection $r\mapsto u_r \in \MIN{\mathsf{t}}{r}{u}$ such that
\begin{equation}
\label{selection}
  0 \in \partial\Psi_u \left( \frac{u_r-u}{r} \right) +
  \diff{\mathsf{t}+r}{u_r}.
\end{equation}
Further,
 there holds
\begin{align}
 &
\label{e:est-a-1}
 \forall\,\mathsf{t} \in [0,T], \ u \in \domainenergy, \ 0<r<T-\mathsf{t},   \  u_r \in
 \MIN{\mathsf{t}}{r}{u}\, :  \quad
 \cg{u_r} \leq C_3 \cg{u},
 \\
 &
\label{e:est-a-2} \lim_{r \down 0} \sup_{u_r \in
\MIN{\mathsf{t}}{r}{u}} \| u_r - u \|=0, \qquad \lim_{r \down
0}\IMIN{\mathsf{t}}{r}{u} = \ene{\mathsf{t}}{u} \quad \text{for all
$\mathsf{t} \in [0,T], \ u \in \domainenergy$}\,,
 \end{align}
 with $C_3$ the constant in~\eqref{e:useful-later}.
Finally,
\begin{equation}
\label{e:differentiability} \text{the map $(0,T-\mathsf{t})\ni r
\mapsto \IMIN{\mathsf{t}}{r}{u}$ is a.e.\ differentiable in $(0,
T-\mathsf{t})$}
\end{equation}
and for every $r_0 \in (0,T-\mathsf{t})$ and every  measurable selection $r \in
(0,r_0] \mapsto u_r \in \MIN{\mathsf{t}}{r}{u} $ there holds
\begin{equation}
\label{e:interp-ene-ineq}
\begin{aligned}
 r_0 \Psi_u \left(\frac{u_{r_0}-u}{r_0}
 \right)+ \int_0^{r_0} \Psi_u^* \left(-\xi_r \right) \dd r +
 \ene{\mathsf{t}+r_0}{u_{r_0}} \leq  \ene{\mathsf{t}}{u} + \int_0^{r_0} \Pt{\mathsf{t}+r}{u_{r}}{\xi_r} \dd r
 \end{aligned}
\end{equation}
where $\xi_r $  is \emph{any} selection  in $
\diff{\mathsf{t}+r}{u_{r}} \cap \left(- \partial \Psi_u
(\frac{u_r-u}{r})\right)$.
\end{lemma}
\noindent \PROOF The \emph{direct method of the calculus of
variations} gives \eqref{nonempty} because of the coercivity
condition~\eqref{eq:17-bis}.
 Further, \cite[Cor.\,III.3, Thm.\,III.6]{Castaing-Valadier77} guarantee
 the existence of a measurable selection $ (0,+\infty)\ni r  \mapsto u_r \in \MIN{\mathsf{t}}{r}{u}$,
  which complies with the Euler equation~\eqref{selection} thanks
to~\eqref{eq:42-bis}.

 Estimate \eqref{e:est-a-1}  follows from the chain of inequalities
\begin{equation}
\label{2nd-ingredient}
\begin{aligned}
\cg{u} \geq
 r\Psi_u \left(\frac{v-u}{r} \right) +
\ene{\mathsf{t}+r}{u_r} \geq \ene{\mathsf{t}+r}{u_r}    \geq
\frac1{C_3}\cg{u_r},
\end{aligned}
\end{equation}
where the first one is  due to the minimality of $u_r$, and the
third one to~\eqref{e:useful-later}. We refer
to~\cite[Lemma\,4.4]{Rossi-Mielke-Savare08} for the proof
of~\eqref{e:est-a-2}, only pointing out that the first limit
in~\eqref{e:est-a-2} follows from the superlinear growth of
$\Psi_u$.

Then, to check~\eqref{e:differentiability} we fix $0 <r_1 <r_2$ and
remark that
\begin{align}
  \label{eq:comparison1}
  \begin{aligned}
 \IMIN{\mathsf{t}}{r_2}{u}  & -\IMIN{\mathsf{t}}{r_1}{u} -  \left(
\ene{\mathsf{t}+r_2}{u_{r_1}} {-} \ene{\mathsf{t}+r_1}{u_{r_1}}
\right)  \\ & \leq
  r_2\Psi_u\left(\frac{u_{r_1}-u}{r_2}\right)
  -r_1 \Psi_u\left(\frac{u_{r_1}-u}{r_1}\right)
 \\ &  \leq
(r_2-r_1)\Psi_u\left(\frac{u_{r_1}-u}{r_2}\right)+ r_1 \left(
\Psi_u\left(\frac{u_{r_1}-u}{r_2}\right)
-\Psi_u\left(\frac{u_{r_1}-u}{r_1}\right)\right)
 \\\ & \leq
(r_2-r_1) \left( \Psi_u\left(\frac{u_{r_1}-u}{r_2}\right) -\left
\langle {w_2^1},{\frac{u_{r_1}-u}{r_2}} \right \rangle\right) =
-(r_2-r_1)\Psi_u^* (w_2^1)
\end{aligned}
\end{align}
where the first inequality follows from \eqref{e:gen-min-prob}, the
second one from algebraic manipulations, the third one by  choosing
some $ w_2^1 \in \partial \Psi_u(({u_{r_1}-u})/{r_2}) $ (which is
nonempty, cf.\ \eqref{e:bounded-operator}), and the
 last passage from
an elementary convex analysis identity. Since $-(r_2-r_1)\Psi_u^*
(w_2^1) \leq 0$ by \eqref{psipos}, we conclude that
\begin{equation}
\label{e:key-point}
\begin{aligned}
\IMIN{\mathsf{t}}{r_2}{u}   & \leq \IMIN{\mathsf{t}}{r_1}{u}+
\ene{\mathsf{t}+r_2}{u_{r_1}} - \ene{\mathsf{t}+r_1}{u_{r_1}}
\\ & \leq \IMIN{\mathsf{t}}{r_1}{u}  + C_1(r_2 - r_1)
\cg{u_{r_1}} \leq \IMIN{\mathsf{t}}{r_1}{u} +  C_1(r_2 - r_1)
C_3\cg{u}     \,,
\end{aligned}
\end{equation}
the second inequality thanks to~\eqref{eq:diffclass_a},  and the
third one to~\eqref{e:est-a-1}. Therefore,
  the map $r \mapsto
\IMIN{\mathsf{t}}{r}{u}$ is given by the sum of a nonincreasing and
of an absolutely continuous function, whence  we deduce that it is almost everywhere
differentiable, viz.\ \eqref{e:differentiability}. In order to conclude
\eqref{e:interp-ene-ineq},  we fix ${r} \in (0,T-t)$,
 outside a negligible set, such that ${r}$ is a
differentiability point   of  the map $r \mapsto
\IMIN{\mathsf{t}}{r}{u}$,
 and we consider
 a selection $w_{h}^{{r}} \in \partial \Psi_u(({u_{{r}} -u})/({{r}+h}))
 $ for $h>0$  sufficiently small. We also fix a sequence
 $h_k \down 0$ such that
\begin{equation}
\label{e:trick} \liminf_{h_k \down 0}
\frac{\ene{\mathsf{t}+{r}+h_k}{u_{{r}}} -
  \ene{\mathsf{t}+{r}}{u_{{r}}}}{h_k} =
\liminf_{h \down 0} \frac{\ene{\mathsf{t}+{r}+h}{u_{{r}}} -
  \ene{\mathsf{t}+{r}}{u_{{r}}}}{h}\,.
\end{equation}
Since $\partial \Psi_u: \V \rightrightarrows \V^*$  is a bounded
operator, from \eqref{e:est-a-2}
 we easily deduce that $\|w_{h_k}^{{r}}\|_* \leq C$, so that
 there exist $w_{{r}} \in \partial \Psi( (u_{{r}} -u)/{{r}})$
 and  a   subsequence
 such that $w_{h_j}^{{r}}\weakto w_{{r}}$ in $\V^*$.
 Then,
we find that
\[
\begin{aligned}  & \Psi_u^*( w_{{r}}) \leq \liminf_{h_j \down
0}\Psi_u^*(w_{h_j}^{{r}})
  \leq \limsup_{h_j \down
0}\Psi_u^*(w_{h_j}^{{r}})
\\ &
 \leq \lim_{h_j \down 0} \left\langle w_{h_j}^{{r}},
\frac{u_{{r}} -u}{{r}+h_j} \right\rangle - \liminf_{h_j \down 0}
\Psi_u \left( \frac{u_{{r}} -u}{{r}+h_j}\right)
 \leq
\left \langle w_{{r}}, \frac{u_{{r}} -u}{{r}} \right\rangle -
\Psi_u\left( \frac{u_{{r}} -u}{{r}}\right) =\Psi_u^*( w_{{r}})\,,
\end{aligned}
\]
using  an elementary convex analysis identity and  that both
$\Psi_u^*$ is
  weakly lower semicontinuous on $\V^*$ and $\Psi_u$ weakly lower semicontinuous on
  $\V$. Therefore $\lim_{j}\Psi_u^*(w_{h_j}^{{r}})=\Psi_u^*(
  w_{{r}})$.
 Since the limit
is independent of the subsequence, we conclude that, for the whole
sequence $w_{h_k}^{{r}}$ there holds
 \[
\lim_{h_k \down 0}\Psi_u^*(w_{h_k}^{{r}})=\Psi_u^*( w_{{r}})=
\Psi_u^*( -\xi_{{r}}) \quad \text{for all }\xi_{{r}} \in
\diff{\mathsf{t}+ {r}}{u_{{r}}} \cap \left(-\partial \Psi_u \left(
\frac{u_{{r}} -u}{{r}}\right) \right)\,, \] the last identity thanks
to condition~\eqref{eq:psi-sum-2}. Then, from
 \eqref{eq:comparison1} we deduce
 \[
 \begin{aligned}
\tfrac{\rmd}{\rmd r} \IMIN{\mathsf{t}}{r}{u}\Restr{r={r}}   +
 \Psi_u^*( -\xi_{{r}}) & =\lim_{h_k \down 0}
\left(\frac{\IMIN{\mathsf{t}}{{r}+h_k}{u}
-\IMIN{\mathsf{t}}{{r}}{u}}{h_k} +\Psi_u^*(w_{h_k}^{{r}})
 \right) \\ & \leq \liminf_{h_k \down 0} \frac{\ene{\mathsf{t}+{r}+h_k}{u_{{r}}} -
  \ene{\mathsf{t}+{r}}{u_{{r}}}}{h_k} \leq \Pt{\mathsf{t+ {r}}}{u_{{r}}}{\xi_{{r}}}\,,
  \end{aligned}
\]
the latter inequality due to~\eqref{e:trick} and~\eqref{e:ass_p_a}.
 Since ${r}$  is arbitrary, we ultimately find
\begin{equation}
\label{e:agognata} \tfrac{\rmd}{\rmd r}
\IMIN{\mathsf{t}}{r}{u}\Restr{r={r}} + \Psi_u^*( -\xi_{{r}})
 \leq  \Pt{\mathsf{t}+ {r}}{u_{{r}}}{\xi_{{r}}}   \qquad
\forae\, {r} \in (0,T-t)\,.
\end{equation}
 Hence, \eqref{e:interp-ene-ineq}
follows from integrating \eqref{e:agognata}  on  the interval
$(0,r_0)$, also using the second of \eqref{e:est-a-2}.
 \QED
 \begin{remark}
\label{rem:alternative} Under assumption \eqref{structure-condition}
as a replacement of~\eqref{eq:psi-sum-2}, it is possible to prove
inequality~\eqref{e:agognata} in the following way. We obtain the
differentiability property~\eqref{e:differentiability} in the same
way as throughout~\eqref{eq:comparison1}--\eqref{e:trick} and then
we observe that, for a
 fixed ${r} \in (0,T-t)$
 outside a negligible set, such that ${r}$ is a
differentiability point   of  the map $r \mapsto
\IMIN{\mathsf{t}}{r}{u}$, we have the following chain of
inequalities for all $h
>0$ (in which we have set $\tilde{u}_{{r},h} = u + \frac{{r}+h}{{r}}
(u_{{r}} - u)$):
\[
\begin{aligned}
\IMIN{\mathsf{t}}{{r}+h}{u} - \IMIN{\mathsf{t}}{{r}}{u} & \leq
\ene{\mathsf{t} + {r}+h}{\tilde{u}_{{r},h}} -\ene{\mathsf{t} +
{r}}{u_{{r}}}   + ({r}+h) \Psi_u \left(\frac{\tilde{u}_{{r},h}
-u}{{r}+h} \right) - {r}\Psi_u \left( \frac{u_{{r}}-u}{{r}}\right)
\\ &
= \ene{\mathsf{t} + {r}+h}{\tilde{u}_{{r},h}} - \ene{\mathsf{t} +
{r}}{\tilde{u}_{{r},h}} + \ene{\mathsf{t} + {r}}{\tilde{u}_{{r},h}}
-\ene{\mathsf{t} + {r}}{u_{{r}}} + h \Psi_u \left(
\frac{u_{{r}}-u}{{r}}\right),
\end{aligned}
 \]
where in the second passage we have used that $\Psi_u
(\frac{\tilde{u}_{{r},h} -u}{{r}+h}) = \Psi_u (
\frac{u_{{r}}-u}{{r}})$. Then,  upon dividing the above inequality
by $h>0$ and taking the $\limsup$ as $h \down 0$,
\eqref{structure-condition} and  \eqref{strengthened} yield (recall
that  $\dire\ene{t}{u;v} = \lim_{h \down 0} \frac1h(
\ene{t}{u+hv}{-}\ene{t}{u})$)
\[
\begin{aligned} \tfrac{\rmd}{\rmd r} \IMIN{\mathsf{t}}{r}{u}\Restr{r={r}}
 & \leq \liminf_{h \down 0} \frac1h \left(\ene{\mathsf{t}+
{r}+h}{\tilde{u}_{{r},h}} {-} \ene{\mathsf{t} +
{r}}{\tilde{u}_{{r},h}}\right) + \dire\cE_{\mathsf{t} +
{r}} \left( u_{{r}}; \frac{u_{{r}}-u}{{r}} \right) + \Psi_u
\left(\frac{u_{{r}}-u}{{r}} \right)
\\
& \leq  \Pt{\mathsf{t+ {r}}}{u_{{r}}}{\xi_{{r}}} + \left\langle
-\xi_{{r}}, \frac{u_{{r}}-u}{{r}} \right\rangle  + \Psi_u
\left(\frac{u_{{r}}-u}{{r}} \right) \\ & =
 \Pt{\mathsf{t+ {r}}}{u_{{r}}}{\xi_{{r}}}-\Psi_u^* \left(-\xi_{{r}} \right) \quad \text{for all \
} \xi_{{r}} \in \diff{\mathsf{t}+ {r}}{u_{{r}}} \cap \left(-\partial
\Psi_u \left( \frac{u_{{r}} -u}{{r}}\right) \right).
\end{aligned}
\]
\end{remark}
\subsection{A priori estimates and compactness for  the approximate solutions}
\label{ss:4.2}
\begin{proposition}[A priori estimates]
  \label{prop:a-prio-est}
 Assume
\eqref{e:4.1}--\eqref{eq:41.1}, and
\eqref{Ezero-bis}--\eqref{e:ass_p_a} for the generalized gradient
system $(\V,\cE,\Psi,\diffname,\Ptname)$. Let $ \pwC \UU{\tau}$, $
\upwC \UU{\:\tau} $, $ \pwL \UU{\tau}$, $ \pwM \UU{\tau}$,  and
$\pwM \xi{\tau}$
 be the interpolants defined by~\eqref{e:pwc}--\eqref{interpmin} and \eqref{interpxi}.
  Then,
  the \emph{discrete upper energy estimate}
 \begin{equation}
 \label{eq:discr-en-ineq-1}
 \begin{aligned}
\int_{\pwC {\mathsf{t}}{\tau}(s)}^{\pwC {\mathsf{t}}{\tau}(t)}
 \Psi_{\upwC{\UU}{\:\tau}(r)}\left( \pwL {\UU'}{\tau}(r)\right) \dd r  &
 + \int_{\pwC {\mathsf{t}}{\tau}(s)}^{\pwC {\mathsf{t}}{\tau}(t)}
 \Psi_{\upwC{\UU}{\:\tau}(r)}^*(-\pwM {\xi}{\tau}(r)) \dd r
+ \ene{\pwC {\mathsf{t}}{\tau}(t)}{\pwC{\UU}{\tau}(t)} \\ &  \leq
\ene{\pwC {\mathsf{t}}{\tau}(s)}{\pwC{\UU}{\tau}(s)} +
\int_{\pwC{\mathsf{t}}{\tau}(s)}^{\pwC {\mathsf{t}}{\tau}(t)}
\Pt{r}{\pwM {\UU}{\tau}(r)}{\pwM {\xi}{\tau}(r)} \dd r
\end{aligned}
\end{equation}
holds for every
  $0\leq s \leq t \leq T$.
  Moreover, there exists a positive constant $S$  such that
  the  following
  estimates are valid for every $\tau>0$:
  \begin{align}
& \label{aprio2-bis} \sup_{t \in (0,T)} |\ene{t}{\pwC \UU{\tau}(t)}|
\leq S, \qquad \sup_{t \in (0,T)} |\ene{t}{\pwM \UU{\tau}(t)}| \leq
S, \qquad
 \sup_{t \in (0,T)}\left|
\Pt{t}{\pwM \UU{\tau}(t)}{\pwM \xi{\tau}(t)}\right|\leq S,
\\
& \label{aprio2} \int_0^T  \Psi_{\upwC{\UU}{\:\tau}(s)}\left( \pwL
{\UU'}{\tau}(s)\right) \dd s   \leq S, \qquad \int_0^T
\Psi_{\upwC{\UU}{\:\tau}(s)}^*(-\pwM {\xi}{\tau}(s)) \dd s \leq S,
\\
& \label{aprio-3-bis} \text{the families $ (\pwL
{\UU'}{\tau})\subset L^1(0,T;\V)$ and $(\pwM {\xi}{\tau})\subset
L^1(0,T;\V^*)$ are uniformly
integrable, and}\\
 & \label{aprio3}
 \sup_{t \in (0,T)} \| \pwC U\tau (t) {-} \upwC
U{\:\tau} (t)\|+ \sup_{t \in (0,T)} \| \pwL U\tau (t) {-} \pwC U\tau
(t)\|+
 \sup_{t \in (0,T)} \| \pwM U\tau (t) {-}
\upwC U{\:\tau} (t)\|=o(1)
\end{align}
as $\tau \down 0.$
\end{proposition}
\noindent \PROOF The proof of Proposition~\ref{prop:a-prio-est}
closely follows  the argument
for~\cite[Prop.\,4.7]{Rossi-Mielke-Savare08}.  For the reader's
convenience we just outline   its main steps here, referring
to~\cite{Rossi-Mielke-Savare08} for the details.

 Let $t_{n-1}, $ $t_{n}$ be two consecutive nodes of
the partition $ \mathscr{P}_\tau $ and let $t \in (t_{n-1},t_n]:$
applying inequality \eqref{e:interp-ene-ineq}
 with the choices
$ \mathsf{t}=t_{n-1}, $ $u= \Utau^{n-1},$
  $r_0= t
-t_{n-1}, $ $u_{r_0}=\pwM {U}{\tau}(t), $ $ u_r=\pwM {U}{\tau}(r) $
and $\xi_r =\pwM {\xi}{\tau}(r)  $ for  $r \in (t_{n-1},t) $ (where
$\pwM{\UU}\tau$ and $\pwM {\xi}{\tau}$ are defined
by~\eqref{interpmin} and \eqref{interpxi}, respectively), we easily
obtain
\begin{equation}
 \label{eq:finer-discr-en-ineq-1}
 \begin{aligned}
(t-t_{n-1})\Psi_{\upwC {\UU}{\:\tau}(t)}& \left(\frac{\pwM
{\UU}{\tau}(t) - \upwC {\UU}{\:\tau}(t)}{t-t_{n-1}} \right)   +
\int_{t_{n-1}}^{t} \Psi_{\upwC {\UU}{\:\tau}(r)}^* ( - \pwM
{\xi}{\tau}(r)) \dd r
\\ & + \ene{t}{\pwM {\UU}{\tau}(t)} \leq  \ene{t_{n-1}}{\pwC
{\UU}{\tau}(t_{n-1})} + \int_{t_{n-1}}^{t} \Pt{r}{\pwM
{\UU}{\tau}(r)}{\pwM {\xi}{\tau}(r)} \dd r\,.
\end{aligned}
\end{equation}
Writing \eqref{eq:finer-discr-en-ineq-1} for $t=t_n$ yields
\begin{equation}
\label{e:useful-again}
\begin{aligned}
\int_{t_{n-1}}^{t_n}
 \Psi_{\upwC {\UU}{\:\tau}(r)} \left( \pwL {\UU'}{\tau}(r)\right) \dd r  & + \int_{t_{n-1}}^{t_n}
 \Psi_{\upwC {\UU}{\:\tau}(r)}^*(-\pwM {\xi}{\tau}(r)) \dd r
+ \ene{t_n}{\pwC{\UU}{\tau}(t_n)} \\ &  \leq
\ene{t_{n-1}}{\pwC{\UU}{\tau}(t_{n-1})} + \int_{t_{n-1}}^{t_n}
\Pt{r}{\pwM {\UU}{\tau}(r)}{\pwM {\xi}{\tau}(r)} \dd r\,.
\end{aligned}
\end{equation}
Upon summing up on the subintervals of the partition, we obtain
\eqref{eq:discr-en-ineq-1}.
 Now, we estimate the
right-hand side of~\eqref{e:useful-again}  via
 \[ \begin{aligned}
\ene{t_{n-1}}{\pwC {\UU}{\tau}(t_{n-1})} +
 \int_{t_{n-1}}^{t_n}  \Pt{s}{\pwM
{\UU}{\tau}(s)}{\pwM {\xi}{\tau}(s)}
 \dd s   &  \leq \ene{t_{n-1}}{\pwC {\UU}{\tau}(t_{n-1})}+  C_2 \int_{t_{n-1}}^{t_n} \cg{\pwM {\UU}{\tau}(s)}
 \dd s \\ & \leq \ene{t_{n-1}}{\pwC {\UU}{\tau}(t_{n-1})}+  C_2 C_3
  \int_{t_{n-1}}^{t_n} \cg{\upwC {\UU}{\:\tau}(s)}
 \dd s\,,
 \end{aligned}
 \]
the first inequality due to~\eqref{e:ass_p_a} and the second one
to~\eqref{e:est-a-1}. On the other hand, condition
\eqref{e:useful-later} yields $\ene{t_n}{\pwC{\UU}{\tau}(t_n)} \geq
C_3^{-1} \cg{\pwC {\UU}{\tau}(t_n)}$. Taking into account the
positivity of the two other integral terms on the left-hand side of
\eqref{e:useful-again} (cf.\ \eqref{psipos}), and  summing it  up on
the intervals of the partition, we obtain the following inequality
\begin{equation}
\label{inequalityA} \cg{\pwC{\UU}{\tau}(t_k)} \leq C \left(
\ene{0}{u_0} + \int_0^{t_k} \cg{\upwC {\UU}{\:\tau}(s)}
 \dd s  +1\right)\,.
\end{equation}
 Then, the first
estimate  in \eqref{aprio2-bis} follows from applying
to~\eqref{inequalityA} a  discrete version of the
   Gronwall lemma (see, e.g.,\ \cite[Lemma\,4.5]{Rossi-Savare06}),
   and the second of~\eqref{aprio2-bis} is  a consequence of \eqref{e:est-a-1}.
The bound in \eqref{aprio2-bis} for the sequence $\{ \Pt{t}{\pwM
\UU{\tau}(t)}{\pwM \xi{\tau}(t)}\}$ again follows from the estimate
for $\ene t {\pwM U{\tau}(t)}$, via~\eqref{e:ass_p_a}.

 Ultimately,  the right-hand side in the discrete energy
inequality \eqref{eq:discr-en-ineq-1} is bounded. Thus, we conclude
\eqref{aprio2}. From \eqref{eq:finer-discr-en-ineq-1} we also deduce
\begin{equation}
\label{e:again-cited} \sup_{t \in [0,T]}
(t-\upwC{\mathsf{t}}{\:\tau}(t)) \Psi_{\upwC{\UU}{\:\tau}(t)}
\left(\frac{\pwM {\UU}{\tau}(t) - \upwC
{\UU}{\:\tau}(t)}{t-\upwC{\mathsf{t}}{\:\tau}(t)} \right) \leq C
\end{equation}
Now, combining this information with \eqref{aprio2-bis} and
\eqref{eq:41.1} (cf.\  \eqref{superlinear}), we infer that
\[
\begin{aligned}
 \forall\, M>0 \ \exists\, S>0 \ &  \forall\,\tau>0, \ t \in [0,T]\, :
\\ &  \| \pwM {\UU}{\tau}(t) - \upwC {\UU}{\:\tau}(t)\| \leq
(t-\upwC{\mathsf{t}}{\:\tau}(t)) S + \frac1M
(t-\upwC{\mathsf{t}}{\:\tau}(t)) \Psi_{\upwC{\UU}{\:\tau}(t)}
\left(\frac{\pwM {\UU}{\tau}(t) - \upwC
{\UU}{\:\tau}(t)}{t-\upwC{\mathsf{t}}{\:\tau}(t)} \right).
\end{aligned}
\]
Estimates \eqref{aprio2} and, again, the superlinear growth
condition~\eqref{eq:41.1}, yield the  uniform integrability of
$(\pwM {\xi}{\tau})$ and $ (\pwL {\UU'}{\tau})$, and the latter in
turn implies~\eqref{aprio3}.
 \QED
\noindent Hereafter, we will  use the short-hand notation
 \begin{equation}
 \label{not-ptau} P_{\tau}(t):= \Pt{t}{\pwM
\UU{\tau}(t)}{\pwM \xi{\tau}(t)}.
\end{equation}
The following result subsumes all compactness information on the
approximate solutions. Some of the convergences below  are stated in
terms of a (limit) \emph{Young measure} associated with the family
$(\pwL U{\tau}', \pwM {\xi}{\tau},
       P_{\tau})_\tau \subset \V \times \V^* \times \R$, the latter space
       endowed with the \emph{weak} topology. The definition of  Young measure, and
       related results,  are recalled in  Appendix A.
       Without going into details,  we may just mention here
 that the aforementioned limit  Young measure allows to
 express the limit as $\tau \downarrow 0 $ of  the sequence $(\mathcal{J}(\pwL U{\tau}', \pwM {\xi}{\tau},
       P_{\tau}))_\tau$ for any \emph{weakly continuous} functional $\mathcal{J}$ on $\V \times \V^*\times \R$
       (and the   $\liminf $ as $\tau \downarrow 0 $ of  the sequence $(\mathcal{H}(\pwL U{\tau}', \pwM {\xi}{\tau},
       P_{\tau}))_\tau$ for any \emph{weakly lower semicontinuous}
        functional $\mathcal{H}$ on $\V \times \V^*\times \R$).
\begin{proposition}[Compactness]
\label{prop:compactness-young}  Assume \eqref{e:4.1}--\eqref{mosco},
and \eqref{Ezero-bis}--\eqref{eq:468}. Then, for every vanishing
sequence $(\tau_k) $ of time-steps there exist a (not relabeled)
subsequence, a curve $u \in \AC ([0,T];\V)$, a function
$\energy:[0,T] \to \R$ of bounded variation, and a time-dependent
Young measure
 $\mmu=\{\mu_t\}_{t\in (0,T)} \in \mathscr{Y}(0,T;\V \times \V^* \times \R),$
 such that as $k\up+\infty$
 \begin{align}
&
   \label{c1} \pwC U{\tau_k}, \, \upwC U{\:\tau_k}, \, \pwL U{\tau_k}, \,  \pwM
   U{\tau_k}
   \rightarrow u \quad \mbox{in $
     L^{\infty}(0,T;\V)   $},
   \\
   &
   \label{c2} \pwL U{\tau_k}' \rightharpoonup u' \quad \mbox{weakly in $
     L^1(0,T;\V) $},
   \\&
   \label{c3}
   \left\{\begin{aligned}
 & \ene{t}{\pwC \UU{\tau}(t)} \to \energy(t)\quad
\text{for all $t \in [0,T]$,} \qquad \energy(0)=\ene{0}{u_0},
\\
& \energy(t) \geq   \ene{t}{u(t)} \quad \text{for all $t \in
[0,T]$,}
\\
& \energy(t)=\ene{t}{u(t)} \quad \foraa\, t \in (0,T),
   \end{aligned}
   \right.
   \end{align}
   and, moreover, $\mmu$  is the limit Young measure
associated with $(\pwL U{\tau_k}', \pwM {\xi}{\tau_k},
       P_{\tau_k})$
      in  the space $\V \times \V^* \times \R$
     (endowed with the weak topology), which implies
   \begin{subequations}
   \begin{align}
  & \label{c4bis}
 u'(t) = \int_{\V \times \V^* \times \R} v \,
\mathrm{d}\mu_t(v,\zeta,p) \quad \foraa\, t \in (0,T),
\\
& \label{c5bis}
 \pwM {\xi}{\tau_k} \weakto \tilde{\xi} \ \ \text{ in $L^1
(0,T;\V^*)$} \quad \text{with }\ \tilde{\xi}(t) := \int_{\V \times
\V^* \times \R} \zeta \, \mathrm{d}\mu_t(v,\zeta,p)  \quad \foraa\,
t \in (0,T),
\\
& \label{c6bis}
\begin{aligned}
& P_{\tau_k} \weaksto P \ \ \text{ in $L^\infty (0,T)$} \text{ with
}\\  & P(t) := \int_{\V \times \V^* \times \R} p \,
\mathrm{d}\mu_t(v,\zeta,p) \leq  \int_{\V \times  \V^* \times \R}
\Pt{t}{u(t)}{\zeta } \, \mathrm{d}\mu_t(v,\zeta,p) \quad \foraa\, t
\in (0,T). \end{aligned}
\end{align}
\end{subequations}
Finally, the following energy inequality  holds for all $0 \leq s
\leq t \leq T$:
\begin{align}
   \label{c6}
   \begin{aligned}
     &  \int_s^t\int_{\V \times \V^* \times \R}  \left( \Psi_{u(r)}(v) +  \Psi_{u(r)}^*(-\zeta)\right) \, \d
     \mu_r (v,\zeta,p)\, \d r+\energy(t) \\
& \leq \energy(s) +
     \int_s^t P(r)\, \d r
      \leq \energy(s) +
     \int_s^t \int_{\V \times \V^* \times \R}
      \Pt{r}{u(r)}{\zeta} \, \d
     \mu_r (v,\zeta,p)\, \d r\,.
     \end{aligned}
   \end{align}
\end{proposition}
\PROOF Let $(\tau_k) $  be a vanishing sequence of time-steps. It
follows from the uniform integrability~\eqref{aprio-3-bis} of  the
sequence $(\pwL U{\tau_k}')$ that $(\pwL U{\tau_k})$ is
equicontinuous on $\V$. Furthermore, \eqref{aprio2-bis} and
 assumption \eqref{eq:17-bis}
give that $\pwM U{\tau_k}$  is contained in some compact subset of
$\V$. Hence $\pwL U{\tau_k}$ is contained in its convex hull, which
is also compact.  Therefore, with the Arzel\`a-Ascoli theorem we
conclude that there exists $u \in \mathrm{C}^0 ([0,T];\V)$  such
that, up to a subsequence,
\begin{equation}
\label{conv1}
 \pwL U{\tau_k} \to u  \quad \text{in $\mathrm{C}^0
([0,T];\V)$.} \end{equation}
 Combining this with \eqref{aprio3}, we
conclude convergences~\eqref{c1}. Next, \eqref{c2} ensues from the
aforementioned uniform integrability of $(\pwL U{\tau_k}')$ via the
Dunford-Pettis criterion (see, e.g.,\
\cite[Cor.\,IV.8.11]{Dunford-Schwartz58}).

 Secondly, from  the third of~\eqref{aprio2-bis} we have that, up to
 a further subsequence,
\begin{equation}
\label{conv2} \text{
 $P_{\tau_k}$ converges weakly$^*$ in $L^\infty
 (0,T)$ to some $P \in L^\infty (0,T)$.}
\end{equation}
  Thus, to prove \eqref{c3} we
 proceed in the same way as for  \cite[Prop.\,4.7]{Rossi-Savare06},
 viz.\ we deduce from the discrete energy inequality
 \eqref{eq:discr-en-ineq-1} that the map
 \[
t \mapsto \eta_\tau(t):= \ene{\pwC
{\mathsf{t}}{\tau}(t)}{\pwC{\UU}{\tau}(t)} - \int_{0}^{\pwC
{\mathsf{t}}{\tau}(t)}P_\tau(r) \dd r \qquad \text{is nonincreasing
on $[0,T]$.}
 \]
Therefore by Helly's theorem there exists $\eta:[0,T]\to \R$,
nonincreasing, such  that, up to a subsequence, $\eta_{\tau_k}(t)
\to \eta(t)$ for all $t \in \R$. In view of \eqref{conv2}, we
conclude that
\begin{equation}
\label{passage1}
 \ene{\pwC
{\mathsf{t}}{\tau_k}(t)}{\pwC{\UU}{\tau_k}(t)} \to \energy(t):=
\eta(t) + \int_0^t P(r) \dd r \qquad \forall\, t \in [0,T].
\end{equation}
This ultimately yields  the first of \eqref{c3}  via
\eqref{eq:diffclass_a}  and \eqref{aprio2-bis}, which give
\begin{equation}
\label{passage2} \left|\ene{\pwC
{\mathsf{t}}{\tau}(t)}{\pwC{\UU}{\tau}(t)}{-}
\ene{t}{\pwC{\UU}{\tau}(t)} \right| \leq  C_1 |\pwC
{\mathsf{t}}{\tau}(t){-}t|\, \cg{\pwC{\UU}{\tau}(t)} \leq S
C_1\,|\pwC {\mathsf{t}}{\tau}(t){-}t| \to 0 \quad \text{as
$\tau\down 0$.}
\end{equation}
Then, the second of \eqref{c3} is a straightforward consequence of
the lower semicontinuity of $\ene{t}{\cdot}$, while the third of
\eqref{c3} follows from assumption \eqref{eq:468}.

In view of  estimates \eqref{aprio2-bis} and \eqref{aprio2} (which
imply the uniform integrability of the sequence $\{\pwM \xi{\tau}
\}$ in $L^1 (0,T;\V^*)$ as well), we  are in the position of applying the Young measure result
in Theorem~\ref{thm.balder-gamma-conv} to the sequence $(\pwL
U{\tau_k}', \pwM {\xi}{\tau_k},
       P_{\tau_k})$, with which we associate a limit Young measure
       $\mmu= \{\mu_t\}_{t \in (0,T)}$ such that for a.a. $t \in (0,T)$
\begin{equation}
\label{e:concentration-bis}
\begin{gathered}
 \mu_t \ \ \text{is concentrated on the set $L(t)$ of the
limit points of $(\pwL U{\tau_k}'(t), \pwM {\xi}{\tau_k}(t),
P_{\tau_k}(t))$}
\\
\text{with respect to the  weak-weak-strong topology of $\V \times
\V^* \times \R$,}
\end{gathered}
\end{equation}
(cf.\ \eqref{e:concentration}), and  there hold
\eqref{heq:gamma-liminf} and  \eqref{eq:35}. Note that the latter
relations imply~\eqref{c4bis}, \eqref{c5bis}, and \eqref{c6bis}.
Then,  from Jensen's inequality we have
\begin{align}
 \label{c4ter}
 \left\{
\begin{array}{l}
  \int_{\V \times \V^* \times \R} \Psi_{u(t)}(v)\,
\mathrm{d}\mu_t(v,\zeta,p)  \geq \Psi_{u(t)}(u'(t)),
\\
 \int_{\V \times
\V^* \times \R} \Psi_{u(t)}^*(-\zeta)\, \mathrm{d}\mu_t(v,\zeta,p)
\geq \Psi_{u(t)}^*(-\tilde{\xi}(t))
\end{array}
\right.\qquad \foraa\, t \in (0,T)\,.
\end{align}
 Passing to
the limit in the Euler equation \eqref{interpxi}, we deduce from
\eqref{eq:468}, from convergence \eqref{c1} for $(\pwM
{U}{\tau_k})$, and  from the first of \eqref{c3},  that for a.a. $t
\in (0,T)$ the set $L(t)$ has the following property
\begin{equation}
\label{e:limit} \text{for all $(v,\zeta,p) \in L(t)$ there holds} \
\ \zeta \in \diff{t}{u(t)}, \ \  p \leq \Pt{t}{u(t)}{\zeta}.
\end{equation}
Hence,  from the latter inequality  and \eqref{eq:35} we also deduce
the  inequality in \eqref{c6bis}. Furthermore, we apply the
$\Gamma$-$\liminf$ inequality \eqref{heq:gamma-liminf} with the
choice $ \mathcal{H}_k (t,v,\zeta,p)=
\Psi_{\upwC{\UU}{\:\tau_k}(t)}(v)$ (notice that
\eqref{hyp:gamma-liminf} is fulfilled in view of assumption
\eqref{mosco},  of \eqref{aprio2-bis}, and of \eqref{c1}). Thus, we
obtain for all $0 \leq s \leq t \leq T$
\begin{equation}
\label{liminf1}
\begin{aligned}
 \liminf_{k \to \infty}
 \int_{\pwC {\mathsf{t}}{\tau_k}(s)}^{\pwC
{\mathsf{t}}{\tau_k}(t)}
 \Psi_{\upwC{\UU}{\:\tau_k}(r)}\left( \pwL {\UU'}{\tau_k}(r)\right)
 \dd
 r  \geq \int_s^t \int_{\V \times \V^* \times \R}  \Psi_{u(r)}(v)
\dd \mu_r(v,\zeta,p) \dd r
 \end{aligned}
 \end{equation}
The choice $\mathcal{H}_k (t,v,\zeta,p)=
\Psi_{\upwC{\UU}{\:\tau}(t)}^*(\zeta)$ (which complies with
\eqref{hyp:gamma-liminf} thanks to~\eqref{e:impo-consequence} and
again \eqref{aprio2-bis})
 obviously gives
\begin{equation}
\label{liminf2}
 \liminf_{k \to \infty}
 \int_{\pwC {\mathsf{t}}{\tau_k}(s)}^{\pwC
{\mathsf{t}}{\tau_k}(t)} \Psi_{\upwC{\UU}{\:\tau_k}(r)}^*(-\pwM
{\xi}{\tau_k}(r))\dd r \geq \int_s^t \int_{\V \times \V^* \times \R}
 \Psi_{u(r)}^*(-\zeta) \dd \mu_r(v,\zeta,p) \dd r.
\end{equation}

Therefore, we pass to the limit in the discrete energy inequality
\eqref{eq:discr-en-ineq-1}. Using \eqref{c3}, \eqref{c6bis},
\eqref{conv2}, \eqref{liminf1}, and \eqref{liminf2},  we conclude
inequality \eqref{c6}.
This completes the proof.
 \QED
\subsection{Proof of  Theorem \ref{thm:viscous2}.}
\label{ss:4.3}
\paragraph{\textbf{Step $1$: a Young measure argument.}}
 It follows from  the a priori estimates
\eqref{aprio2-bis}--\eqref{aprio2} and from \eqref{c3},
\eqref{e:concentration-bis}--\eqref{liminf1} that the curve $u \in
\AC([0,T;\V)$ and the Young measure $\{ \mu_t\}_{t \in (0,T)}$
comply with assumptions \eqref{ass-chr1}--\eqref{ass-chr2-bis} of
Theorem \ref{th:app-chain-rule}. Therefore, the map $t \mapsto
\ene{t}{u(t)}$ is absolutely continuous and we have the following
chain of inequalities
\begin{equation}
\label{one}
\begin{aligned}
 \int_0^t & \int_{\V \times \V^* \times \R}\left( \Psi_{u(r)}(v){+}\Psi_{u(r)}^*(-\zeta)
 \right) \, \d
     \mu_r (v,\zeta,p)\, \d r + \ene{t}{u(t)}
\\   & \leq \ene{0}{u(0)} +  \int_0^t P(r)\, \d r
\leq
   \ene{t}{u(t)} + \int_0^t \int_{\V \times \V^* \times \R}  \langle -\zeta,u'(r) \rangle  \, \d
     \mu_r (v,\zeta,p)\, \d r
\end{aligned}
\end{equation}
where the first inequality follows from \eqref{c6} (written for $t
\in (0,T] $ and $s=0$) and from the second of \eqref{c3}, while the
second inequality is a consequence of the \emph{Young measure}
chain-rule inequality \eqref{ymeasure-chain}. Taking into account
inequality~\eqref{c4ter}, we thus conclude
\[
\int_0^t  \int_{\V \times \V^* \times \R}\left(
\Psi_{u(r)}(u'(r)){+}\Psi_{u(r)}^*(-\zeta) {-}  \langle -\zeta,u'(r)
\rangle
 \right)
\, \d
     \mu_r (v,\zeta,p)\, \d r \leq 0.
\]
Since the integrand is nonnegative, we find
\begin{equation}
\label{e:ultimately}
 \int_{\V \times \V^*\times \R}
\left(\Psi_{u(t)}(u'(t)){+}\Psi^*_{u(t)}(-\zeta) {-} \langle -\zeta,
u'(t)\rangle\right) \, \d
     \mu_t (v,\zeta,p)=0 \quad \foraa\, t \in (0,T).
\end{equation}
Now, it follows from the above discussion that all inequalities in
\eqref{one} indeed hold as equalities. Again using the chain-rule
inequality~\eqref{ymeasure-chain}, it is easy to deduce that, for
almost all $t \in (0,T)$, we have
\begin{equation}
\label{chainequality}
\begin{aligned}
 \int_{\V \times \V^*\times \R}
\Big(\Psi_{u(t)}(u'(t)) & {+}\Psi^*_{u(t)}(-\zeta) {-} p \Big) \, \d
     \mu_t (v,\zeta,p)
\\  & = \int_{\V \times \V^*\times \R} \left(\langle -\zeta,
u'(t)\rangle {-} p \right) \, \d
     \mu_t (v,\zeta,p)  = -\frac \d{\d t}\ene t{u(t)}
     \end{aligned}
\end{equation}
Thus,  we conclude for every $0\leq s \leq t \leq T$ the energy
identity
\begin{equation}
\label{enid1-new}
 \int_s^t  \int_{\V \times \V^* \times \R}\left(
 \Psi_{u(r)}(v){+}\Psi_{u(r)}^*(-\zeta) {-}p
 \right) \, \d
     \mu_r (v,\zeta,p)\, \d r
  = \ene{s}{u(s)}-\ene{t}{u(t)}\,.
\end{equation}

\paragraph{\textbf{Step $2$: a measurable selection.}}
Let us consider the measure $\nu_t:= (\pi_{2,3})_{\#}(\mu_t)$, i.e.\
the marginal of $\mu_t$ with respect to the $(\zeta,p)$-component,
defined by $(\pi_{2,3})_{\#}(\mu_t) (B) = \mu_t(\pi_{2,3}^{-1}(B)) $
for all  $B \in \mathscr{B}(\V^* \times \R)$ (the Borel
$\sigma$-algebra of $\V^* \times \R$).  As a consequence of
\eqref{e:limit} and \eqref{e:ultimately}, for almost all $t \in
(0,T)$  the measure $\nu_t$ is concentrated on the set
\begin{equation}
\label{new-set} \mathcal{S}(t,u(t),u'(t)):= \left\{ (\zeta,p)\in
\V^*\times \R\, : \ \zeta \in \diff t{u(t)} \cap (-\partial
\Psi_{u(t)}(u'(t)), \ p \leq \Pt{t}{u(t)}{\zeta} \right\},
\end{equation}
namely
\begin{equation}
\label{concentrated}
 \nu_t(( \V^*\times \R) \setminus
\mathcal{S}(t,u(t),u'(t)))=0 \quad \foraa\, t \in (0,T).
\end{equation}
In particular, for almost all $t \in (0,T)$ the set
$\mathcal{S}(t,u(t),u'(t))$ is nonempty. Then, Lemma \ref{l:helpful}
in the appendix below guarantees that there exists
 a measurable selection $t \in (0,T) \mapsto (\xi(t),p(t)) \in
 \mathcal{S}(t,u(t),u'(t))$
 such that
 \begin{equation}
 \label{psi-star-min}
\Psi_{u(t)}^*(-\xi(t)) -p(t) = \min_{(\zeta,p) \in
\mathcal{S}(t,u(t),u'(t))} \{ \Psi_{u(t)}^*(-\zeta) -p\} \doteq
\mathcal{M}^*(t) \quad \foraa\, t \in (0,T).  \end{equation} In
particular, $\xi$ satisfies equation \eqref{xi-equation}, hence we
conclude that $u$ solves the Cauchy problem for \eqref{eq:1-bis}.
 In fact, we have
\begin{equation}
\label{diss-bounded} \int_0^T \Psi_{u(t)}^* (-\xi(t))\, \d t
<+\infty,
\end{equation}
which in particular yields $\xi \in L^1 (0,T; \V^*)$ via
\eqref{eq:41.1}. To check \eqref{diss-bounded}, it is sufficient to
observe that
\[
\begin{aligned}
\int_0^T  \Psi_{u(t)}^* (-\xi(t))\, \d t &   \leq  \int_0^T
\mathcal{M}^*(t) + p(t)\, \d t  \\ & \leq \int_0^T  \int_{\V^*
\times \R}\left( \Psi_{u(r)}^*(-\zeta) - p
 \right) \, \d
     \nu_r (\zeta,p)\, \d r + \int_0^T \Pt {r}{u(r)}{\xi(r)} \, \d r
     \\ & \leq \int_0^T  \int_{\V \times \V^*
\times \R}\left( \Psi_{u(r)}^*(-\zeta) - p
 \right) \, \d
     \mu_r (v,\zeta,p)\, \d r +C_2 \int_0^T \cg {u(r)} \, \d r  \leq C
     \end{aligned}
\]
where the second inequality ensues from \eqref{concentrated} and the
fact that $p(t) \leq  \Pt {t}{u(t)}{\xi(t)}$ for almost all $t\in
(0,T)$, the third inequality from~\eqref{e:ass_p_a}, and the last
one from  \eqref{enid1-new} and the fact that $\sup_{t \in (0,T)}
\cg{u(t)}<+\infty$.

\paragraph{\textbf{Step $3$: proof of the energy identity~\eqref{enid1}.}}
On the one hand, we observe that for every $0 \leq s \leq t \leq T$
there holds
\begin{equation}
\label{left-energy}
\begin{aligned}
 \int_s^t  & \left(\Psi_{u(r)}(u'(r))  {+} \Psi_{u(r)}^*(-\xi(r)) {-}\Pt
 {r}{u(r)}{\xi(r)}\right)
\, \d r \\ & \leq \int_s^t \Psi_{u(r)}(u'(r)) {+}
\Psi_{u(r)}^*(-\xi(r)) {-}p(r) \, \d r \leq
 \ene{s}{u(s)}-\ene{t}{u(t)},
 \end{aligned}
\end{equation}
where the first estimate follows from  $p(t) \leq \Pt
 {t}{u(t)}{\xi(t)}$ for almost all $t \in (0,T)$ by definition of
 $\mathcal{S}(t,u(t),u'(t))$, and the second estimate is due to
\eqref{enid1-new}, combined with \eqref{concentrated} and
\eqref{psi-star-min}.
 On the other hand,
applying the chain-rule inequality \eqref{eq:48strong} to the pair
$(u,\xi)$ we have
\begin{equation}
\label{right-inequality} \ene{t}{u(t)} - \ene{s}{u(s)}\geq \int_s^t
\left( \langle \xi(r), u'(r) \rangle  {+}\Pt
 {r}{u(r)}{\xi(r)} \right) \, \d r \quad \text{for every $0 \leq s \leq t \leq T$.}
\end{equation}
Combining \eqref{left-energy} and \eqref{right-inequality} and
arguing in the same way as throughout
\eqref{one}--\eqref{e:ultimately}, we obtain that all inequalities
in~\eqref{left-energy} ultimately hold as equalities; in particular,
$ p(t) = \Pt {t}{u(t)}{\xi(t)}$ for almost all $t \in (0,T).$
 We
have thus proved that the pair $(u,\xi)$ satisfies the energy
identity~\eqref{enid1}. A comparison between the latter and the
\emph{Young-measure energy identity}~\eqref{enid1-new} also reveals
that, for almost all $t \in (0,T)$,
\begin{align}
\label{almost-all-0}
  & \Psi_{u(t)}(u'(t))= \int_{V \times \V^* \times \R}
\Psi_{u(t)}(v) \, \d \mu_t (v,\zeta,p),
\\ &
 \label{almost-all}
\begin{aligned}
 \Psi_{u(t)}^*(-\xi(t))  - \Pt {t}{u(t)}{\xi(t)}
& = \min_{(\zeta,p) \in \mathcal{S}(t,u(t),u'(t))} \{
\Psi_{u(t)}^*(-\zeta) -p\} \\ & = \Psi_{u(t)}^*(-\zeta) {-}p \quad
\text{for $\nu_t$-almost all $(\zeta,p) \in \V^*\times \R$.}
\end{aligned}
\end{align}
Taking into account that $\Psi_{u(t)}^*(-\zeta)=
\Psi_{u(t)}^*(-\xi(t))$ due to condition \eqref{eq:psi-sum-2}, we
thus conclude the \emph{maximum selection principle}
\begin{equation}
\label{max-sel} \Pt {t}{u(t)}{\xi(t)} = \max \{p \, : \ (\zeta,p)
\in \mathcal{S}(t,u(t),u'(t)) \}.
\end{equation}

\paragraph{\textbf{Step $4$: enhanced convergences.}}
Convergences \eqref{convu4}--\eqref{convu5} and
\eqref{enhanced-middle} are proved by passing to the limit in
\eqref{eq:discr-en-ineq-1}, written for $s=0$ and $t \in [0,T]$. We
use the short-hand notation \eqref{not-ptau}, as well as
\[
 A_k(t) =\int_{0}^{\pwC {\mathsf{t}}{\tau_k}(t)}
\Psi_{\upwC{\UU}{\:\tau_k}(r)}\left( \pwL {\UU'}{\tau_k}(r)\right)
\dd r, \ \   B_k(t)  = \int_{0}^{\pwC {\mathsf{t}}{\tau_k}(t)}
\Psi_{\upwC{\UU}{\:\tau_k}(r)}^*(-\pwM {\xi}{\tau_k}(r))\dd r,\ \
 C_k (t) =\ene{\pwC {\mathsf{t}}{\tau_k}(t)}{\pwC U{\tau_k}(t)}.
\]
For all $t \in [0,T]$ we find
\begin{equation}
\label{liminf-limsup}
\begin{array}{ll}
  & \int_{0}^{t} \Psi_{u(r)}(u'(r))\,\d r
  +\int_{0}^{t} \Psi_{u(r)}^*(-\xi(r))\, \d r\ +   \ene {t}{u(t)}
 \\
  & =\int_{0}^{t} \Psi_{u(r)}(u'(r))\,\d r
  +\int_{0}^{t} \int_{\V \times \V^*\times \R}\Psi^*_{u(r)}(-\zeta)\,  d\nu_r (v,\zeta,p) \, \d r\ +   \ene {t}{u(t)}
\\ & \leq \liminf_{k \to \infty} A_k(t) +
 \liminf_{k \to \infty}B_k(t)   + \liminf_{k \to \infty}C_k(t)
\\ &  \leq \limsup_{k \to \infty}\left( A_k(t)  + B_k(t)   +C_k(t) \right)
\\ &   \leq \limsup_{k \to \infty} \ene{0}{\pwC{\UU}{\tau_k}(0)} +
\limsup_{k \to \infty}  \int_{0}^{\pwC {\mathsf{t}}{\tau_k}(t)}
P_{\tau}(r) \dd r
\\ & \leq \ene {0}{u(0)} + \int_0^t  \Pt {r}{u(r)}{\xi(r)}\, \d r\\
 & = \int_{0}^{t} \Psi_{u(r)}(u'(r))\,\d r
  +\int_{0}^{t} \Psi_{u(r)}^*(-\xi(r))\, \d r\ +   \ene {t}{u(t)}
\end{array}
\end{equation}
where the first identity follows from \eqref{almost-all}, the second
estimate
 from \eqref{c3},
\eqref{passage1}--\eqref{passage2}, and
\eqref{liminf1}--\eqref{liminf2}, the third estimate is trivial and
the fourth one  ensues from  inequality \eqref{eq:discr-en-ineq-1},
whereas the fifth estimate is a consequence of \eqref{c6bis}, and
the sixth identity is due to \eqref{enid1}. Altogether, all
inequalities in~\eqref{liminf-limsup} turn out to be equalities, and
with an elementary argument we conclude \eqref{convu4},
\eqref{convu5}, as well as \eqref{enhanced-middle}.

\paragraph{\textbf{Step $5$: the strictly convex case.}}
Finally, if we further assume \eqref{strict-convexity}, from
\eqref{almost-all-0}, \eqref{almost-all},
 and the
strict convexity of $\Psi^*_{u(t)}(\cdot)$ we infer
\begin{equation}
\label{deltas} 
(\pi_{2})_{\#}(\mu_t)= \delta_{\xi(t)} \qquad \forae \ t \in (0,T).
\end{equation}
Hence, from~\eqref{c5bis} we deduce convergence~\eqref{convu3}. This
concludes the proof of Theorem~\ref{thm:viscous2}. \QED
\begin{remark}
\upshape
\label{rmk:deltas}
Notice that, if in addition we assume $\Psi_u$ to be strictly convex
for all $u$, then we also have $(\pi_{1})_{\#}(\mu_t)=
\delta_{u'(t)}$. The latter relation, joint
with~\eqref{deltas},  yields
\[
\pwL {\UU'}{\tau_k}(t) \weakto u'(t), \qquad \pwM {\xi}{\tau_k}(t)
\weakto \xi(t) \qquad \forae \  t \in (0,T).
\]
\end{remark}
\paragraph{\textbf{Sketch of the proof of Theorem \ref{thm:viscous-stability}.}}
For every $n \in \N$, the solution pair $(u_n,\xi_n)$ fulfills the
energy identity associated with the Cauchy problem \eqref{cauchy-n},
namely  there holds
\begin{equation}
\label{e:starting-point}
 \int_{0}^{t} \big( \Psi_{u_n(r)}^n(u_n'(r))
 {+}(\Psi^n_{u_n(r)})^*(-\xi_n(r)) \big)
 \,\d r +   \enei n {t}{u_n(t)} =
\enei n {0}{u_0^n} + \int_0^t  \Pti n{r}{u_n(r)}{\xi_n(r)}\dd r
\end{equation}
for all $t \in [0,T]$. From \eqref{e:starting-point} we deduce all the a priori estimates
on the sequence $(u_n,\xi_n)$,
with the very same arguments as in the proof of Proposition \ref{prop:a-prio-est}.
 Indeed, we exploit condition \eqref{data-convergence}
on $\enei n {0}{u_0^n}$,
and use \eqref{e:ass_p_a} (for a constant uniform
with respect to $n \in \N$), to estimate the terms on the right-hand side of \eqref{e:starting-point}.
Then, all of the terms on the left-hand side are  estimated as well.
 Combining this  with
the coercivity properties of the potentials $(\Psi_u^n)$, viz.
\[
\begin{aligned}
& \! \! \! \! \! \!  \forall\, R>0, \,  M>0\,
\\ & \! \! \! \! \! \!  \left\{
\begin{array}{lllllll}
\exists K>0  &  \forall\, u\in \domainenergy \text{ with } \sup_{n
\in \N}\cgn u\leq R  & \forall\, v \in \V &  :  &  \|v \|\geq K  &
\Rightarrow &  \Psi_u^n(v) \geq M\|v\|,
\\
 \exists K^*>0  &  \forall\, u\in \domainenergy \text{ with }
\sup_{n \in \N}\cgn u\leq R  &  \forall\, \xi \in \V^* &  :
 &  \|\xi \|_{*} \geq K^* &  \Rightarrow & (\Psi_u^n)^*(\xi) \geq
M\|\xi\|_{*},
\end{array}
\right.
\end{aligned}
\]
we have that the sequence $(u_n') \subset L^1 (0,T;\V)$ is
 uniformly
integrable. Furthermore, the estimate $\sup_{n \in \N} \sup_{t \in [0,T]} \enei {n}{t}{u_n(t)}$
yields compactness which, combined with uniform integrability, ensures convergences \eqref{convu1-n}, along a
subsequence, to some curve $u \in \AC ([0,T];\V)$.
Like in  Proposition~\ref{prop:compactness-young}, up to a subsequence
we also find some limit Young measure for the
sequence
$(u_n',\xi_n,P_n)$, with $P_n(t):=  \Pti n{t}{u_n(t)}{\xi_n(t)} $.

Finally, in order to pass to the limit as $n \to \infty$, we reproduce on the time-continuous level
the arguments developed in Steps 1--4 of the proof of Theorem \ref{thm:viscous2}. Namely,
combining
 semicontinuity arguments with properties \eqref{closure-n} and \eqref{ppsi-n}, we  take the limit as $n
 \to \infty$ of \eqref{e:starting-point}, and deduce that the curve $u$ fulfills  the upper energy estimate.
 We obtain the lower energy estimate from the chain rule, and in this way we conclude that $u$ is a
 solution to the Cauchy problem for \eqref{eq:1-bis}.
\QED



\appendix
\section{Young measure tools}
\label{s:a-1}
 \noindent In this section, we collect some results on
 parametrized (or Young) measures
with values in infinite-dimensional  spaces, see   e.g.\
\cite{Balder84, Balder85, Ball89, Balder_lecture_notes98,
castaing-fitte-valadier04, Valadier90}. In particular, we shall
focus on Young measures with values in a reflexive Banach space
$\mathcal{V}$. The definitions and results we are going to recall
below, apply in Section \ref{ss:4.2} (cf. Proposition
\ref{prop:compactness-young}), to the space \medskip $\mathcal{V}=
\V \times \V^* \times \R$.

\paragraph{\textbf{Notation.}}
Given an interval $I \subset \R$, we  denote by $\mathscr{L}_{I}$
the $\sigma$-algebra of the Lebesgue measurable subsets of $I$ and,
given
   a reflexive Banach space $\mathcal{V}$,
 by $\mathscr B(\mathcal{V})$ its Borel $\sigma$-algebra.
We  use the symbol $\otimes$ for product $\sigma$-algebrae. We
recall that a $\mathscr{L}_{I} \otimes \mathscr
B(\mathcal{V})$-measurable function $h : I \times \mathcal{V} \to
(-\infty,+\infty]$ is a \emph{normal integrand} if for a.a. $t \in
(0,T)$ the map $x \mapsto h_t(x)= h(t,x)$ is lower semicontinuous on
$\mathcal{V}$.

We consider the space $\mathcal{V}$ endowed   with the \emph{weak}
topology, and say that a $\mathscr{L}_{(0,T)} \otimes \mathscr
B(\mathcal{V})$--measurable functional $\mathcal{H}: (0,T) \times
\mathcal{V} \to (-\infty,+\infty]$ is a \emph{weakly-normal
integrand}  if for a.a. $t \in (0,T)$ the map
\begin{equation}
\label{def:wws}
\begin{gathered}
 \text{$w \mapsto h(t,w)$
is sequentially lower semicontinuous on $\mathcal{V}$ w.r.t.\ the
weak topology.}
\end{gathered}
\end{equation}
We   denote by $\mathscr{M} (0,T; \mathcal{V})$ the set of all
$\mathscr{L}_{(0,T)}$-measurable functions $y: (0,T) \to
\mathcal{V}$. A sequence $(w_n) \subset\mathscr{M} (0,T;
\mathcal{V}) $ is said to be \emph{weakly-tight}  if there exists a
weakly-normal integrand $\mathcal{H}:(0,T) \times \mathcal{V}
\rightarrow (-\infty,+\infty]  $ such that the map
\[
\begin{gathered}
\text{$w \mapsto  \mathcal{H}_t(w)$ has compact sublevels w.r.t.\
 the weak topology of $\mathcal{V}$, and}
\\
\sup_n \int_0^T \mathcal{H}(t,w_n(t))  \dd t <\infty.
\end{gathered}
\]
 \begin{definition}[\bf (Time-dependent) Young measures]
  \label{parametrized_measures}
  A \emph{Young measure} in the space $\mathcal{V} $
  is a family
  $\mmu:=\{\mu_t\}_{t \in (0,T)} $ of Borel probability measures
  on $ \mathcal{V} $
  such that the map on $(0,T)$
\begin{equation}
\label{cond:mea} t \mapsto \mu_{t}(B) \quad \mbox{is}\quad
{\mathscr{L}_{(0,T)}}\mbox{-measurable} \quad \text{for all } B \in
\mathscr{B}(\mathcal{V}).
\end{equation}
We denote by $\mathscr{Y}(0,T; \mathcal{V})$ the set of all Young
measures in $\mathcal{V} $.
\end{definition}
The following $\Gamma$-$\liminf$ result is a straightforward
consequence of \cite[Thm.\,4.2]{Stef08?BEPD}.
\begin{theorem}
\label{thm.balder-gamma-conv}
 Let $\{\mathcal{H}_n\}$, $\mathcal H
: (0,T) \times \mathcal{V}\to (-\infty,+\infty]$ be weakly-normal
integrands such that for all $ w \in \mathcal{V}$ and $ \forae\, t
\in (0,T)$
\begin{equation}
\label{hyp:gamma-liminf}
\begin{aligned}
\mathcal{H} (t,w) \leq \inf \Big\{\liminf_{n\to
    \infty}\mathcal{H}_n (t,w_n)\, : & \
w_n \weakto w \ \text{in  $\mathcal{V}$}
 \Big\}
\,.
\end{aligned}
\end{equation}
Let $(w_n) \subset \mathscr{M} (0,T;\mathcal{V}) $ be a weakly-tight
sequence. Then, there exist a subsequence $(w_{n_k})$ and a Young
measure $\mmu=\{\mu_t\}_{t \in (0,T)}$ such that for a.a. $t \in
(0,T)$
\begin{equation}
\label{e:concentration}
\begin{gathered}
  \mbox{$ \mu_{t} $ is
      concentrated on
      the set
      $ L(t):=
      \bigcap_{p=1}^{\infty}\overline{\big\{w_{n_k}(t)\,: \ k\ge p\big\}}^{\mathrm{w}}$}
      \end{gathered}
  \end{equation}
of the limit points of the sequence $(w_{n_k}(t))$ with respect to
the weak topology of $\mathcal{V}$ and,
 if the
sequence $t \mapsto \mathcal{H}_{n_k}^- (t,w_{n_k}(t))$ is uniformly
integrable, there holds
\begin{equation}
\label{heq:gamma-liminf} \liminf_{k \to +\infty} \int_0^T
\mathcal{H}_{n_k} (t,w_{n_k}(t))\, \dd t \geq
\int_0^T\int_{\mathcal{V}} \mathcal{H} (t,w)\, \dd \mu_t(w) \, \dd
t\,.
\end{equation}
\end{theorem}

 As a corollary (the reader is referred to the discussion in
\cite{Stef08?BEPD} for more details), we have a generalization of
the so-called {\em
  Fundamental Theorem of Young measures}, see \cite[Thm.\,3.2]{Rossi-Savare06}
for the case of the weak topology in Hilbert spaces, and the
classical results~\cite[Thm.\,1]{Balder84},
\cite[Thm.\,2.2]{Balder85},
\cite[Thm.\,4.2]{Balder_lecture_notes98},
\cite[Thm.\,16]{Valadier90}.
 \begin{theorem}[\textbf{The Fundamental Theorem for strong-weak-weak topologies}]
  \label{cor:app-1}
Let $1 \leq p \leq \infty$ and let  $(w_n) \subset L^p
(0,T;\mathcal{V})$ be a bounded sequence. If $p=1$, suppose further
that
 $(w_n)  $ is uniformly integrable in $  L^1 (0,T;\mathcal{V})$.   Then,
  there exists a subsequence $(w_{n_k})$ and
  a Young  measure
  $ \mmu=\{ \mu_{t} \}_{t \in (0,T)}  \in \mathscr{Y}(0,T;\mathcal{V})$
  such that for a.a.\ $t\in (0,T)$ relation \eqref{e:concentration} holds
  and, setting
  \[
\rmw(t):=\int_{\mathcal{V}} w \, \dd \mu_t (w)  \qquad \foraa\, t
\in (0,T)\,,
  \]
there holds
\begin{equation}
  \label{eq:35}
w_{n_k} \weakto \rmw \ \ \text{ in $L^p (0,T;\V)$},
\end{equation}
with $\weakto$ replaced by $\weaksto$ if $p=\infty$.
\end{theorem}
\noindent In fact, in Section \ref{ss:4.2} (cf.\ Proposition
\ref{prop:compactness-young}), Theorem \ref{cor:app-1} applies to
the sequence $w_k (t):= (\pwL U{\tau_k}'(t), \pwM {\xi}{\tau_k}(t),
P_{\tau_k}(t))\subset \V \times \V^*\times \R$.
\section{Extension of the chain rule to Young measures}
\noindent From now on, we work with Young measures valued in the
reflexive space $\mathcal{V}:= \V \times \V^*\times \R$, whose
elements are denoted by $(v,\zeta,p)$. The main result of this
section is a \emph{Young measure version} of the chain-rule
inequality \eqref{eq:48strong}.
\begin{proposition}
\label{th:app-chain-rule} In the framework of
\eqref{e:4.1}--\eqref{eq:41.1}, suppose that $\cE$ complies with
\eqref{Ezero-bis}--\eqref{eq:468}. Let $u \in \AC ([0,T];\V)$ be an
absolutely continuous curve
 such that
 \begin{equation}
 \label{ass-chr1}
(t,u(t)) \in \domaindiff \quad \forae\, t \in (0,T), \ \ \text{and}
\ \ \sup_{t \in (0,T)} \ene{t}{u(t)}<+\infty.
\end{equation}
 Let $\mmu=\{ \mu_t
\}_{t \in (0,T)}$ be a Young measure in $\V \times \V^* \times \R$
such that
\begin{align}
& \label{ass-chr2} \int_0^T \int_{\V \times \V^* \times \R}\left(
\Psi_{u(t)}(v)+ \Psi^*_{u(t)}(-\zeta) \right) \dd \mu_t
(v,\zeta,p)\dd t<+\infty,
\\
& \label{ass-chr3}
 \foraa\, t \in (0,T) \quad u'(t) = \int_{\V\times \V^* \times \R} v \, \d \mu_t(v,\zeta,p),
 \\
 &
 \label{ass-chr2-bis} \foraa\, t \in (0,T) \quad \text{for all $(v,\xi,p) \in \mathrm{supp}(\mu_t)$ there holds
 } \ \xi \in \diff t{u(t)},  \ p \leq \Pt{t}{u(t)}{\xi}.
\end{align}
Then,
\begin{equation}
\label{ymeasure-chain}
\begin{gathered}
\text{the map} \ \ t \mapsto \ene{t}{u(t)} \   \ \text{is absolutely
continuous on $(0,T)$, and}
\\
   \frac \d{\d t}\ene t{u(t)} \geq \int_{\V \times \V^* \times \R}
  \left( \la u'(t),\zeta \ra+
    p \right) \dd \mu_t(v,\zeta,p)
    \quad
    \text{for a.a.\ }t\in (0,T).
\end{gathered}
\end{equation}
\end{proposition}
\noindent The proof of this result closely follows the argument for
\cite[Thm.\,3.3]{Rossi-Savare06}, to which we constantly refer the
reader.
\\
\PROOF We split the argument in three claims.
\par\noindent
\textbf{Claim $1$:} \textit{let us set for almost all $t \in (0,T)$}
\begin{equation}
\label{Kt} \mathcal{K}(t,u(t)):= \left\{ (\xi,p)\in \V^*\times \R\,
: \ \xi \in \diff t{u(t)}, \ p \leq \Pt{t}{u(t)}{\xi} \right\}.
\end{equation}
 \textit{There exists a sequence  of strongly
measurable maps $(\xi_n, p_n) :(0,T) \to \V^* \times \R$ such that}
\begin{equation}
\label{measurability} \{(\xi_n(t),p_n(t)) \, : \ n \in \N \} \subset
\mathcal{K}(t,u(t)) \subset \overline{ \{(\xi_n(t),p_n(t)) \, : \ n
\in \N \}},
\end{equation}
\textit{($\overline{A}$ denoting the closure, with respect to the topology of
$\V^*\times \R$, of a subset $A \subset\V^*\times \R$)}.
\\
 First of
all, let us observe  that  the set
\begin{equation} \label{gt}
\begin{gathered}
 K: = \left\{(t,u,\xi,p) \in [0,T]\times \V \times
\V^*\times \R \, : \ \xi \in \diff tu,  \ p \leq \Pt{t}{u}{\xi} \right\}
\\
\text{is a Borel set of $[0,T]\times \V \times \V^* \times \R$.}
\end{gathered}
 \end{equation}
This follows from the fact that
$\mathrm{graph}(\diffname)$ is a Borel set of
$[0,T]\times \V \times \V^*$, and $\Ptname: \mathrm{graph}(\diffname) \to \R$ a Borel function.
Now, it follows
from~\eqref{ass-chr1} and~\eqref{ass-chr3}--\eqref{ass-chr2-bis}
that there exists a subset $\mathcal{T} \subset (0,T)$ of full
measure such that $\mathcal{K}(t,u(t)) \neq \emptyset$ for all $t
\in [0,T]$. Let us then consider the graph of the multivalued
function $t \in \mathcal{T} \mapsto \mathcal{K}(t,u(t)) \subset
\V^*\times \R$, i.e.\ the set
\[
\mathscr{K} = \left\{ (t,\xi,p)\in \mathcal{T} \times \V^*\times
\R\, : \ (\xi,p) \in \mathcal{K}(t,u(t))   \right\} = \left\{
(t,\xi,p)\in \mathcal{T} \times \V^*\times \R\, : \ (t,u(t),\xi,p) \in
K\right\}.
\]
Due to the latter representation, to~\eqref{gt} and to the fact that
the function $u :[0,T] \to \V$ is Borelian, we can conclude that
$\mathscr{K}$ is  a Borel set of $\mathcal{T} \times \V^*\times
\R$.
 Therefore, \eqref{measurability} ensues from
\cite[Thm.\,III.22]{Castaing-Valadier77}.
\par\noindent
\textbf{Claim $2$:} \textit{it is possible to choose the  measurable
maps $\xi_n :(0,T) \to \V^*$ fulfilling \eqref{measurability} in
such a way that}
\begin{equation}
\label{measurability-bound} \xi_n \in L^1 (0,T;\V^*) \ \text{ for
every $n\in \N$ and } \ \sup_n \int_0^T \Psi^*_{u(t)}(-\xi_n(t))\,
\d t <+\infty.
\end{equation}
We set
\[
\mathcal{M}_*(t): = \min_{(\xi,p) \in \mathcal{K}(t,u(t))}
\Psi_{u(t)}^* (-\xi) \quad \text{for almost all $t \in (0,T)$.}
\]
It follows from \eqref{measurability} that
\begin{equation}
\label{e:technical-measurability} \text{the map  \ } t \mapsto
\mathcal{M}_*(t) = \min_{n} \Psi_{u(t)}^* (-\xi_n(t)) \ \text{ is
measurable on $(0,T)$.}
\end{equation}
Moreover, \eqref{ass-chr2} and \eqref{ass-chr3} yield that
\begin{equation}
\label{consequence-bis} \int_0^T \mathcal{M}_*(t)\, \d t \leq
\int_0^T \int_{\V^*}\Psi_{u(t)}^* (-\zeta)\, \d \mu_t (v,\zeta,p) \,
\d t<+\infty.
\end{equation}
 Arguing as the proof of \cite[Lemma\,3.4]{Rossi-Savare06},  we
 recursively define the following family of subsets of
 $\mathcal{T}$,  i.e.
 \[
A_0:= \emptyset, \quad A_k:= \left\{t \in \mathcal{T}\,  : \
\Psi_{u(t)}^* (-\xi_k(t)) \leq \mathcal{M}_*(t) + 1 \right\}
\setminus \bigcup_{j=0}^{k-1}A_j.
 \]
 Due to~\eqref{e:technical-measurability}, for every $k \in \N$ the
 set $A_k$ is measurable  and, by construction,
 the family $\{A_k\}_{k\in\N}$
is disjoint with $\bigcup_{k\in\N} A_k=\mathcal{T}$. Hence, we set
\[
  \bar{\xi}(t):=\sum_{k=1}^{+\infty}\xi_k(t)\nchi_{A_k}(t) \quad
   \bar{p}(t):=\sum_{k=1}^{+\infty}p_k(t)\nchi_{A_k}(t)\quad
  \text{for all $t \in \mathcal{T}$.}\]
Notice that  the map $t \mapsto (\bar\xi(t),\bar{p}(t))$
 is a measurable selection of $\mathcal{K}(t,u(t))$, and
that
\begin{equation}
\label{passage} \Psi_{u(t)}^* (-\bar{\xi}(t)) \leq \mathcal{M}_*(t)
+ 1 \quad \text{ for every $t\in\mathcal T$.}
\end{equation}
 In particular, it
follows from \eqref{eq:41.1} that $\bar{\xi} \in L^1 (0,T;\V^*)$.
Then, we use $(\bar \xi,\bar{p})$ to construct a new countable
family of functions, setting
\[
  (\xi_{n,k}(t), p_{n_k}(t)):=
  \left\{
  \begin{array}{ll}
    (\xi_n(t),p_n(t)) &\text{if } \Psi_{u(t)}^* (-\xi_n(t))\leq k,\\
    (\bar{\xi}(t),\bar{p}(t))&\text{otherwise,}
    \end{array}
  \right.
\]
such that the functions $\xi_{n,k} $ belong to $ L^1(0,T;\V^*)$, and
the pairs $(\xi_{n,k}, p_{n_k})$  satisfy
\[
  (\xi_{n,k}(t), p_{n_k}(t)) \in \mathcal{K}(t,u(t)),\quad
  \{(\xi_{n}(t),p_n(t)):n\in\N\}\subset
  \{(\xi_{n,k}(t), p_{n_k}(t)):n,k\in \N\}\quad
  \text{if }t\in\mathcal \mathcal{T}, \quad\square
\]
as well as estimate \eqref{measurability-bound}, in view of
\eqref{consequence-bis}.
\par\noindent
\textbf{Claim $3$:} \textit{inequality \eqref{ymeasure-chain}
holds.}
\\
Indeed,  it follows from Claims $1$ and $2$ that we can apply the
chain rule \eqref{eq:48strong} to the pairs $(u,\xi_n)$
 for all $n \in \N$ (indeed, estimate \eqref{ass-chr2} and
the first of \eqref{ass-chr3} yield $\int_0^T \Psi_{u(t)}(u'(t))\,
\d t<+\infty$). Therefore, we conclude  for all $n \in \N$ that
there exists a set of full measure $\mathcal{T}_n \subset
\mathcal{T} $
\[
\begin{gathered}
     \text{the map  $t\mapsto \ene t{u(t)}$ is absolutely continuous
     and}
     \\
    \frac \d{\d t}\ene t{u(t)} \geq \la \xi_n (t),u'(t)\ra+
    \Pt{t}{u(t)}{\xi_n(t)} \geq \la \xi_n (t),u'(t)\ra+ p_n(t)
    \quad
    \text{for all }t\in \mathcal{T}_n.
    \end{gathered}
\]
Thus, we infer
\begin{equation}
\label{inequality-convexified} \frac \d{\d t}\ene t{u(t)} \geq \la
\xi ,u'(t)\ra+
    p
    \quad
    \text{for all }
(\xi,p) \in \overline{\mathrm{conv}}(\mathcal{K}(t,u(t))), \
\text{for all }
    t\in \mathcal{T}_\infty,
\end{equation}
with $\overline{\mathrm{conv}}(\mathcal{K}(t,u(t)))$  the closed
convex hull of $\mathcal{K}(t,u(t))$ and
 $\mathcal{T}_\infty = \bigcap_{n\in \N} \mathcal{T}_n$ (note that $\mathcal{T}_\infty$ is  a
 subset of $(0,T)$ of full measure, too).
Then, \eqref{ymeasure-chain} follows upon integrating
\eqref{inequality-convexified} with respect to the measure $\mu_t$,
again taking into account \eqref{ass-chr2-bis} and \eqref{ass-chr3}.
  \QED
\noindent We conclude with the following
\begin{lemma}[Measurable selection]
\label{l:helpful}
 In the framework of
\eqref{e:4.1}--\eqref{eq:41.1}, suppose that $\cE$ complies with
\eqref{Ezero-bis}--\eqref{eq:468}. Let $u \in \AC ([0,T];\V)$ be an
absolutely continuous curve
 complying with~\eqref{ass-chr1}, and suppose that the set
 \begin{equation}
 \label{non-emptiness}
 \begin{aligned}
 \mathcal{S}(t,u(t),u'(t)): =   \left\{ (\zeta,p)\in
\V^*\times \R\, : \ \zeta \in \diff t{u(t)} \cap (-\partial
\Psi_{u(t)}(u'(t))), \ p \leq \Pt{t}{u(t)}{\zeta} \right\}, \\
\text{ is nonempty for almost all $t \in (0,T)$.}
\end{aligned}
\end{equation}
Then, there exists measurable functions $\xi :(0,T) \to \V^*$,
$p:(0,T) \to \R$ such that
\begin{equation}
\label{properties-selection}
\begin{aligned}
   (\xi(t),p(t)) \in
\argmin \{ \Psi_{u(t)}^*(-\zeta) -p\, : \ (\zeta,p) \in
\mathcal{S}(t,u(t),u'(t))\} \quad  \foraa\, t \in (0,T)\,.
\end{aligned}
\end{equation}
\end{lemma}
\PROOF First of all, let us observe that
\begin{equation}
\label{argmin-not-empty} \argmin \{ \Psi_{u(t)}^*(-\zeta) -p\, : \
(\zeta,p) \in \mathcal{S}(t,u(t),u'(t))\} \neq \emptyset \quad
\foraa\, t \in (0,T)\,.
\end{equation}
 For, let $(\xi_n,p_n) \subset \mathcal{S}(t,u(t),u'(t))$
be a infimizing sequence: then, there exist some positive constants
$C$ and $C'$ such that
\begin{equation}
\label{chain-bis} \Psi_{u(t)}^*(-\xi_n) \leq C  + p_n \leq C +
\Pt{t}{u(t)}{\xi_n} \leq C + C_2 \cg{u(t)}  \leq C' \quad \text{for
every $n \in \N$,}
\end{equation}
where the first inequality trivially  follows from the fact that
$(\xi_n,p_n)$ is infimizing, the second one from the definition of
$\mathcal{S}(t,u(t),u'(t))$, the third one from \eqref{e:ass_p_a},
and the last one from \eqref{e:useful-later} and
assumption~\eqref{ass-chr1}. In view of the latter, and of the
superlinear growth condition~\eqref{eq:41.1}, from the bound for
$\Psi_{u(t)}^*(-\xi_n)$ we infer that $\sup_n \|\xi_n \|_*
<+\infty$. It is also clear from \eqref{chain-bis} that $\sup_n
|p_n|<+\infty$, therefore there exist $(\xi_*,p_*)$ such that, up to
a subsequence, $\xi_n \weakto \xi_*$ in $\V^*$ and $p_n \to p_*$.
Exploiting the closedness condition~\eqref{eq:468}, we infer that
$(\xi_*,p_*) \in \mathcal{S}(t,u(t),u'(t))$, and by lower
semicontinuity
\[
\Psi_{u(t)}^*(-\xi_*) -p_* \leq \liminf_{n \to \infty}\left(
\Psi_{u(t)}^*(-\xi_n)-p_n\right) = \inf_{(\zeta,p) \in
\mathcal{S}(t,u(t),u'(t))} \{ \Psi_{u(t)}^*(-\zeta) -p\},
\]
whence \eqref{argmin-not-empty}.

 In order to prove
\eqref{properties-selection}, we first observe that for all $t \in
[0,T]$ the set
\begin{equation}
 \label{st}
\begin{aligned}
  S :  & = \left\{(t,u,v,\xi,p) \in [0,T] \times \V \times V \times
\V^*\times \R \, : \ \xi \in \diff tu \cap (-\partial\Psi_u (v)),  \ p \leq
\Pt{t}{u}{\xi} \right\} \\ &  \text{is a Borel set of $ [0,T] \times \V  \times \V \times
\V^* \times \R$.} \end{aligned}
 \end{equation}
 This follows from the same arguments as for~\eqref{gt} in
 Proposition~\ref{th:app-chain-rule}. Now, due to
 \eqref{non-emptiness} there exists a subset $\mathcal{T}' \subset
 (0,T)$ of full measure such that $\mathcal{S}(t,u(t),u'(t)) \neq
 \emptyset$ for all $t \in \mathcal{T}'$. We thus consider the graph
 of the multivalued function $t \in \mathcal{T}' \mapsto \mathcal{S}(t,u(t),u'(t)) \subset
\V^*\times \R$, i.e.\ the set
\[
\begin{aligned}
\mathscr{S}  & :  = \left\{ (t,\xi,p)\in \mathcal{T}' \times
\V^*\times \R\, : \ (\xi,p) \in \mathcal{S}(t,u(t),u'(t)) \right\}\\
&  = \left\{ (t,\xi,p)\in \mathcal{T}' \times \V^*\times \R\, : \
(t,u(t),u'(t),\xi,p) \in S\right\}.
\end{aligned}
\]
Then, we combine the latter representation of $\mathscr{S}$ with~\eqref{st}, and
 the fact that the functions $u :(0,T) \to
\V$ and  $u' :(0,T) \to \V$ are Borelian up to choosing a suitable representative for $u'$.
Thus, we conclude that $\mathscr{S}$ is a Borel
set of $\mathcal{T}' \times \V^*\times \R$. Hence, the existence
of a measurable selection $(\xi,p) $ as in
\eqref{properties-selection} is a consequence  of \cite[Cor.\,III.3, Thm.\,III.6]{Castaing-Valadier77},
cf.\ also  Filippov's theorem, see e.g.\ \cite[Thm.~8.2.11]{AubFra90SVA}.
 \QED 

\bibliographystyle{siam}
\bibliography{bibliografia2009,refs_doublynonlinear09,bibliography_rate-independent10}
\end{document}